\documentclass[11pt]{amsart}

\input xy             
\xyoption{all}    

\newcommand{\vMax}{{\underline{\mbox{Max} }  } \ }

\newcommand{\vword}{{\mbox{w-ord}}}

\newcommand{\Gdi}{{\mathcal G}^{(d-i)}}

\newcommand{\calo}{{\mathcal {O}}}
\newcommand{\Sing}{\mbox{Sing\ }}
\newcommand{\ord}{\mbox{ord}}

\usepackage{latexsym, amsfonts, amsmath, amssymb, amscd, epsfig, xypic}

\newtheorem{Theorem}{Theorem}[section]
\newtheorem{Lemma}[Theorem]{Lemma}
\newtheorem{Proposition}[Theorem]{Proposition}
\newtheorem{Corollary}[Theorem]{Corollary}
\newtheorem{Definition}[Theorem]{Definition}
\newtheorem{Remark}[Theorem]{Remark}
\newtheorem{Paragraph}[Theorem]{}
\newtheorem{Example}[Theorem]{Example}

\title[]{Singularities in positive characteristic, stratification and simplification of the singular locus}

\author{A. Bravo  \and O. Villamayor U.}

\address{Instituto de Ciencias Matematicas CSIC-UAM-UC3M-UCM and Dpto. Matem\'aticas,
Facultad de Ciencias, Universidad Aut\'onoma
de Madrid, Canto Blanco 28049 Madrid, Spain}
\email{ana.bravo@uam.es, villamayor@uam.es}
\thanks{The authors were partially supported by MTM 2006-10548.}


 \date{}
\dedicatory{} \commby{}

\begin{document}
\maketitle
\begin{abstract}
We introduce an upper semi-continuous function that stratifies the
highest multiplicity locus of a hypersurface in arbitrary
characteristic  (over a perfect field). The blow-up along the maximum stratum defined by
this function leads to a form of simplification of the
singularities,  known as the reduction to the monomial case. 
\end{abstract}

\tableofcontents

\part*{Introduction}

Resolution of singularities is a classical and central problem in
algebraic geometry. Using a non-constructive argument,   Hironaka proved in the mid sixties (cf.
\cite{Hironaka64})  that the singularities of varieties over fields of characteristic
zero  could always be resolved. Several constructive
(algorithmic) proofs of resolution of singularities have been
published since the late eighties  (\cite{BM}, \cite{EH},
\cite{Hauser}, \cite{kollar}, \cite{Villa89}, \cite{Villa92},
\cite{WLL}).

\

There are some results on low dimensional varieties over arbitrary fields   by Abhyankar,  and more recently by
some  other authors  (see also \cite{CJS}, \cite{CP1}), \cite{CP2}), \cite{Cut}) but the general question
of resolution 
 remains open.  Another important  contribution is  de Jong's work on alterations,
 which  provides a weaker statement, but it is strong  enough  for certain applications. 

\

In the next paragraphs  we  describe some of the main ideas  of  the
proof of algorithmic resolution in characteristic zero, paying  special attention to the part of the argument  that fails in positive
characteristic. After this exposition,   we explain the results obtained in
this paper.

 \

Suppose that  $X$ is a  reduced scheme over a field of characteristic
zero. An algorithmic  desingularization of  $X$ can be obtained 
in two steps, say A and B. In Step A,  a suitable sequence of monoidal
transformations on smooth centers is defined so as to produce a
{\em simplification} of the singularities of $X$, meaning that they can be assumed to be contained in  some smooth lower dimensional scheme, where they can be  described  in terms of an  ideal  of a divisor with normal crossings support.   This step  is accomplished
by using  an inductive argument. Once this process is finished, it is
said that $X$ is within  the {\em monomial case}. In Step B, the
monomial case is treated: a combinatorial argument
 leads to a resolution of singularities of $X$.

\

 {\bf Step A. Simplification of singularities.}   The
goal  here  is to define a stratification of any reduced scheme $X$ by
means of an upper semi-continuous function  $\Gamma_X:X\to (\Lambda,\leq)$,  where $(\Lambda,\leq)$  is a  fixed well ordered set,   such that:
\begin{itemize}
\item [1.] The maximum value of $\Gamma_X$ is attained  on a smooth closed
subscheme, \underline{Max}$\Gamma_X$, which describes the {\em
 worst singularities of} $X$, and  the minimum value of $\Gamma_X$ is attained on the non-singular locus of $X$. 
\item [2.] The blow-up of $X$ along \underline{Max}$\Gamma_X$, $\pi_1: X_1\to
X$,  {\em improves   the singularities of $X$} in the
following sense:  as in (1.), a new   upper semi-continuous function 
$\Gamma_{X_1}:X_1\to (\Lambda,\leq)$ is defined with the following properties:
\begin{itemize}
\item[(a)] If $x_1\in X_1\setminus\pi^{-1}_1(\mbox{\underline{Max}}\Gamma_X)
\simeq X\setminus\mbox{\underline{Max}}\Gamma_X$ maps to $x\in
X$, then
 $\Gamma_X(x)=\Gamma_{X_1}(x_1)$.
 \item[(b)] The maximum value of $\Gamma_X$, $\mbox{Max} \Gamma_X$, drops,
 i.e., $\mbox{ Max} \Gamma_{X_1}  <\mbox{Max}\Gamma_X$.
 \end{itemize}
\item[3.]  A  {\em simplification } of
the singularities of $X$ is obtained after a finite number of monoidal transformations defined by $\Gamma_{X_0}, \Gamma_{X_1}, \ldots, \Gamma_{X_{n-1}}$, which can have the following form:
 \begin{equation} \label{simplification}
 X=X_0 \leftarrow X_1 \leftarrow
\ldots \leftarrow X_n. \end{equation} This is usually referred to
as a reduction to the  {\em monomial case}. Its meaning will
be explored in more detail below. 
 \end{itemize}

 \

The question now is
how to define the functions $\Gamma_{X_i}: X_i\to (\Lambda,\leq)$,
for $i=0,1,\ldots,n-1$. Suppose that $X$ is embedded in a smooth $d$-dimensional scheme
$V^{(d)}$. A first approximation  to $\Gamma_X$ is to consider the
order of the ideal of definition of $X$, ${\mathcal
I}^{(d)}\subset {\mathcal O}_{V^{(d)}}$, at the closed points
$x\in V^{(d)}$. Notice that this defines an upper semi-continuous
function, say,
$$\Gamma^{(d)}: X\to {\mathbb Z}_{\geq 0}.$$  However this function
is {\em too coarse}, since in general it does not satisfy
properties (1.) and (2.b). This becomes  clear if we assume, for
instance, that $X$ is a hypersurface: then the maximum order of
${\mathcal I}(X)$ is located  in the set of points  of $X$ with maximum
multiplicity, which may not  be a smooth
 subscheme of $V^{(d)}$.

\

An important point in characteristic zero is how the
previous function $\Gamma^{(d)}$ can be {\em refined}, thanks to
the existence of {\em hypersurfaces of maximal contact}: the
closed set $\mbox{\underline{Max}} \Gamma^{(d)}$ is locally
contained in a smooth $(d-1)$-dimensional scheme $V^{(d-1)}$, and
it can be described by means of an ideal ${\mathcal
I}^{(d-1)}\subset {\mathcal O}_{V^{(d-1)}}$. Then a new function
$\Gamma^{(d-1)}:\mbox{\underline{Max}} \Gamma^{(d)}\subset
V^{(d-1)} \to {\mathbb Z}_{\geq 0}$ is defined on  
$\mbox{\underline{Max}} \Gamma^{(d)}\subset V^{(d-1)}$, now using
the order at closed points of the ideal ${\mathcal I}^{(d-1)}$.
The construction  of an upper semi-continuous function $\Gamma:
X\to (\Lambda,\leq)$ satisfying all properties (1.),  (2.) and (3.)  from
Step A is  made by collecting the information from
$\Gamma^{(d)}, \Gamma^{(d-1)}, \ldots$, and employing an inductive argument.

\

 {\bf Step B. The monomial case.} {\rm Once Step A is
 accomplished as a composition of a finite number of monoidal
 transformations, $V^{(d)}\leftarrow V_n^{(d)}$,
it can be assumed that, locally,  the  {\em
worst singularities of $X_n\subset V_n^{(d)}$}  are contained in some 
   smooth-$(d-e)$-dimensional  closed subscheme
$V^{(d-e)}_n\subset  V^{(d)}_n$,  and that they can be defined in terms of an ideal
of a divisor with normal crossings support. Here $e\geq 1$, and  the
ideal of a divisor with normal crossings is called a monomial
ideal. In this case (monomial case) it is relatively easy to
enlarge sequence (\ref{simplification})  to resolve the singularities of $X$.

\

{\bf Hypersurfaces of maximal contact and the problems in
positive characteristic.} Hypersurfaces of maximal contact play a
central role in constructive resolution in characteristic zero. This topic
is related to  Abhyankar's notion of Tschirnhausen transform (or
Tschirnhausen substitution): given the  equation of a singular
embedded hypersurface, Tschirnhausen provides the equation of a
smooth hypersurface, in the ambient space, that locally contains
the highest multiplicity locus of the singular hypersurface
(maximal contact). Moreover, this containment  is preserved by
monoidal transformations with centers included  in the locus with the highest
multiplicity (see \cite{Aab}, \cite{Ab}, \cite{AbM}).

\

It is the work of  J. Giraud where  hypersurfaces of maximal
contact arise by means of techniques that  involve differential
operators on smooth schemes, in a first attempt to address the
problem of embedded desingularization in arbitrary characteristic
(cf. \cite{Giraud1975}). This approach, which uses differential
operators, played a central role in  the development of algorithmic
resolution of singularities in characteristic zero.

\

However, in positive characteristic hypersurfaces of maximal
contact may not exist (see for instance  \cite{HHauser} and
\cite{Nar}). For this reason, the argument explained in Step A
cannot be extended to this setting.

\

{\bf The aim of this paper.} In this paper we  show that the stratifying functions with the
prescribed properties of  Step A,  are defined in any characteristic.  This extension is made possible by the introduction of two tools:  the characteristic free  techniques
introduced in \cite{hpositive} that avoid maximal contact;  and the  function defined thanks to  Main Theorem
\ref{ordertau} (see Definition \ref{deford}). 

\

 When applied to characteristic zero, this approach leads to  the same  upper semi-continuous functions, and to the reduction to the monomial case, as, for instance, in \cite{EncVil97:Tirol}. This  coincidence in characteristic zero is fully proved in  \cite{mariluz}. 
 
\

 Thus, 
 Main Theorem \ref{ordertau} provides a characteristic free form of induction; and as a consequence,  a simplifications of singularities  over arbitrary characteristic (in the spirit of the reduction to  monomial case) can   be obtained  (this is shown in  \cite{positive}).   In other words,  the so called ``reduction to
the monomial case" is possible in positive characteristic. 
This means that the 
monomial case arises in some lower dimension via
induction (cf. \cite{positive}). 

\

{\bf About Step B.}  In characteristic zero, once  $X$ is within the monomial case,    
 it is possible to define another upper semi-continuous function, of  a combinatorial nature. 
 This function again stratifies $X$ in smooth strata;  then the iteration of a finite sequence of  blow-ups at its  maximum stratum    easily leads  to a resolution of singularities of $X$.
  Therefore,  the  extension of Step A to arbitrary
characteristic, i.e., the reduction to the monomial case in any characteristic, opens a
door to new invariants. Over fields of positive characteristic,
hypersurfaces whose highest multiplicity locus is  in the ``monomial
case",  turn out to have very particular properties. The treatment of
this specific case, which we hope to address in the future, would
imply resolution of singularities over arbitrary fields (see
\cite{BVV}).

\

{\bf Other approaches to resolution in arbitrary characteristic.}
The form of induction dealt with in this paper is different from the one used in
\cite{Hironaka08}, and also different from those used in the
Kawanoue program (\cite{kaw} and \cite{MK}), and in \cite{WL}. All
these approaches   strongly rely on techniques of differential
operators, but differ in their approach to induction. In these
works,  restriction to smooth hypersurfaces of maximal contact is
replaced by a notion of restriction to singular hypersurfaces,
which are also, in some generalized sense, of maximal contact.
Some questions concerning an approach to stratification of
singularities in positive  characteristic have also been addressed
in \cite{EH}. There are other invariants for singularities in
positive characteristic,  studied in works
of Cossart, Hauser, and Moh, which are also related to the
problem of embedded resolution of singularities. We include an example to illustrate
the effect of Step A on a particular singularity. We chose here
one of Hauser's kangaroo points (see  Example \ref{canguro}); these are
singularities where  pathologies specific to  positive
characteristic arise. For instance,  the definition of a resolution invariant which works in characteristic zero, but  over positive characteristic   increases after a finite 
number of blow-ups at closed points (cf.  \cite{HHHauser} and \cite{HauserBul} for full details).

\

\noindent {\bf Elimination: a strategy for overcoming the failure
of maximal contact in positive characteristic}

\

{\bf Maximal contact vs. elimination. }  In \cite{hpositive},   the
concepts of {\em hypersurfaces of maximal contact} and {\em
restriction to hypersurfaces of maximal contact} are replaced by the
notion of {\em transversal projections} and {\em elimination
algebras} (respectively). This allows us to use  induction in
any characteristic. We illustrate this procedure next.

\

Let $X $ be an algebraic variety embedded in a $d$-dimensional smooth scheme $V^{(d)}$ over a 
field $k$, and let $x\in X$ be an $n$-fold 
closed point. Then, locally, in a neighborhood of $x$, Weierstrass Preparation Theorem 
provides a projection to a $d-1$-dimensional smooth scheme and the following situation 
can be assumed to hold. 

\
   
Let  $A$  be a smooth $k$-algebra,  and  let  $X$ be a  hypersurface in $\mbox{Spec}(A[Z])$ defined by
$$f(Z)=Z^n+a_1Z^{n-1}+\ldots+a_n\in A[Z].$$

Suppose  $\Upsilon_n$ is  the set of $n$-fold points of $X$ (i.e., the
points of multiplicity $n=\mbox{deg}f(Z))$, and let
$B=A[Z]/\langle f(Z)\rangle$. Then the natural projection
$$\beta: \mbox{Spec} (B) \to \mbox{Spec}(A)$$ is a finite morphism, and Zariski's
multiplicity formula for projections ensures that the map induces
a bijection between $\Upsilon_n$ and $\beta(\Upsilon_n)$ (see
\ref{Motivation} for more details).

\

Now,  in this setting 

\begin{enumerate} \item [i.]  A suitable  $A[Z]$-Rees algebra ${\mathcal
G}$, with singular locus  $\mbox{Sing
}{\mathcal G} = \Upsilon_n$, is  associated to $f(Z)$  (we refer to  Definition
\ref{singularlocus} for the notion of singular locus of a Rees
algebra).

\item[ii.]   An {\em elimination
algebra} ${\mathcal
R}_{\mathcal G}$ is associated to ${\mathcal G}$.  This elimination algebra 
 is a Rees algebra over the ring $A$
(independent of the variable $Z$), and has the property   that
$\beta(\Upsilon_n)=\beta(\mbox{Sing }{\mathcal G})\subset
\mbox{Sing }{\mathcal R}_{\mathcal G}$.
\end{enumerate}

\

Following this approach, the highest multiplicity locus of $X$ is {\em
projected bijectively} to the  smooth scheme $\mbox{Spec}(A)$. This projection replaces  the restriction, 
used in characteristic zero, of the highest multiplicity locus of $X$ to a hypersurface of maximal contact.

\

More specifically, and parallel to the  arguments given in Step A, the ambient space
$V^{(d)}$ here is $\mbox{Spec}(A[Z])$,  and  the ideal ${\mathcal
I}^{(d)}$ is replaced by ${\mathcal G}$. Then the restriction to
the hypersurface of maximal contact $V^{(d-1)}$ is replaced by the
projection to $\mbox{Spec}(A)$,  and  the information encoded by
${\mathcal I}^{(d-1)}$ is now encoded by the Rees algebra ${\mathcal
R}_{\mathcal G}$, over the ring  $A$, so it is independent of the
variable $Z$.

\

We now proceed in the same manner as  in characteristic zero.
 Consider the order function for Rees algebras (see \ref{orderRees} for the precise definition):
$$\begin{array}{rrcl}
\Gamma^{(d)}: &  \mbox{Sing }{\mathcal G} & \to & {\mathbb
Q}_{\geq 0}\\
 & x & \to & \mbox{ord}_x{\mathcal G}.
 \end{array}$$
As it happens with ideals, this order function turns out to be too coarse to satisfy
 properties (1.) and (2.b) of  Step A. Thus, we refine this function   by considering
$$\begin{array}{rrcl}
\Gamma^{(d-1)}: & \mbox{\underline{Max}}\Gamma^{(d)} & \to &
{\mathbb
Q}_{\geq 0}\\
 & x & \to & \mbox{ord}_{x_1}{\mathcal R}_{\mathcal
 G},
 \end{array}$$
where $x_1=\beta(x)$,  and where we use the fact that
$\beta(\mbox{\underline{Max}}\Gamma^{(d)}) \subset \mbox{Sing
}{\mathcal R}_{\mathcal
 G}$.  After this,  the construction  proceeds by induction.

 \

 In the present  article,  we show that this procedure can be iterated  and that the functions
 we construct are independent of the choice of the projection.
  By  induction  on the dimension,  we can
construct an upper semi-continuous function to some well ordered set,
$\Gamma: X\to (\Lambda,\geq)$,   that stratifies $\Upsilon_n$ in
smooth strata,  and also fulfills  conditions (1.), (2.) and (3.) in step A  (see Main Theorem \ref{ordertau} and Theorem
\ref{stratification}).   Therefore,  this   part of the inductive argument, used in characteristic zero,   can be extended to positive 
characteristic.  The iteration of blow-ups at the maximum strata of this function,  leads to a simplification of the singularities in some lower dimension (see \cite{positive}).

\

The idea of projecting the maximum multiplicity locus to a smooth
subscheme, follows along the lines   of Jung's procedure for
resolving  hypersurface singularities,  via a simplification of the discriminant.

\

The paper is organized  in five parts. Part \ref{Char0} contains a
brief exposition of the main ideas behind algorithmic resolution
of singularities over fields of characteristic zero. As indicated
above, algorithmic resolution is achieved in two steps A and B:
 a reduction to the monomial case, and then a treatment of
the monomial case. We indicate the conditions required to extend step A to arbitrary fields.

\

The problem of resolution over arbitrary fields will be formulated in
terms of Rees algebras, so Part \ref{partRees} is devoted to
recalling  some notions of  the theory. In Sections \ref{Rees} and
\ref{DiffRees} we present a brief introduction.   Special
attention will be  paid to  Rees algebras enriched with the  action of
differential
 operators. We will see how these objects   provide a suitable framework to
 define our invariants.
 In Section \ref{transformationsweak}  we discuss  transformations of Rees algebras, and Hironaka's notion of weak equivalence. The least number of variables needed to express
 the initial form of a hypersurface at a singular point is a central
 invariant in the theory. In fact, these are the variables
 that can be eliminated from the problem. This is what we call the $\tau$-invariant
 and its study, addressed in Section \ref{tangentcones},
 will play a central role in the construction of our stratifying
 function.

 \

 Part \ref{partElimination} is dedicated to presenting elimination algebras and to
 reviewing
 some of their properties: Section \ref{evui}
 contains a detailed study of {\em universal  elimination algebras}, and the specialization
 to usual elimination algebras via change of  base rings. In  Section \ref{localprojection}
 we study conditions under which elimination can be defined, with special attention to the fact
 that these conditions are open.

 \

  Part \ref{partMain} contains the main results: in Section
  \ref{eliminationmonoidal} we study the behavior of elimination
  algebras under monoidal transformations;
 Main Theorem \ref{ordertau} is stated in   Section \ref{sectionorder}, and
 the proof is given in Section \ref{sectionorder1}.
  Theorem \ref{ordertau} makes it possible to construct
   the upper semi-continuous functions introduced in Definition \ref{deford},
 and   Sections \ref{nonsimple} and  \ref{stratificationtheorem}
 are devoted to the stratification that results from this
 function.

 \

 Finally in Part \ref{partEpilogue} we explain how to use
 our results to reach the {\em monomial case}. We refer to \cite{positive} for full details.

\

{\em Acknowledgements:} We are indebted to A. Benito Sualdea and
M. L. Garc\'{\i}a Escamilla for their careful review of this
manuscript, and the numerous useful  suggestions they provided. We are also grateful
to the referee's work. His/her valuable comments that have helped us  to
improve the presentation of this paper. Prof. J.M. Aldaz also gave us some useful suggestions.



\part{Algorithmic resolution of singularities over fields of characteristic
zero} \label{Char0}

Here we briefly present the main ideas underlying the algorithmic resolution of singularities  in characteristic
zero. This is done here, as in \cite{Villa89} and \cite{Villa92},  in terms of  {\em pairs} and {\em
basic objects}.   We conclude this first part with an example to illustrate how
the algorithm works (cf. Example
\ref{examplerespairs}). For more details we
refer the reader to the introductory presentation in \cite{EncVil97:Tirol}.

\

We will show how the language of Rees
algebras, which is required for this new  approach  in arbitrary characteristic,
parallels the one of pairs (see \ref{paraleloragadprs},
\ref{ragadprs3} and \ref{rgn}).

\section{The language of pairs and basic objects}\label{ragadprs}

{\em Pairs}  provide a suitable language to formulate resolution 
problems.  However, once we start the process of resolution,  exceptional divisors
  appear, and we  need to keep track on this information too.
The information provided by pairs and  the ambient space where they are defined,   
 as well as the set of exceptional divisors that the resolution
process produces,  are  codified  in terms of {\em basic objects}.

\

Suppose  that $X$ is a hypersurface embedded in some
smooth $d$-dimensional space $V$. If our goal is to resolve the
singularities of $X$, then    we will start by paying attention to the
worst singularities of $X$, namely, those points where the
multiplicity of $X$ is the highest, say $b$. We will see that the
natural pair associated to this closed set is $({\mathcal
I}(X),b)$. If there are no exceptional divisors to take care of,
the basic object we will be interested in is $(V,({\mathcal
I}(X),b), \{\emptyset\})$.

\begin{Paragraph}\label{unouno}{\bf Pairs.}
 {\rm A {\em pair}
 $(J,b)$ on a smooth scheme $V$ is
defined by a non-zero sheaf of ideals $J\subset \calo_V$ and a
positive integer $b$. The    {\em singular locus}  of a  pair
$(J,b)$ consists of  the set of points in $V$ where $J$ has order
at least $b$, i.e.,
$$\Sing(J,b):=\{ x\in V | \nu_x(J)\geq b \},$$ where  $\nu_x$
denotes the order function at the regular local ring $
\calo_{V,x}$. The set $\Sing(J,b)$ is closed in $V$.

\

As indicated above,  it is typical  to take $J$ as the
sheaf of ideals defining a hypersurface $X\subset V$, and $b$ as
the maximum of the multiplicities at points of $X$.

\

Hironaka defines the function:
\begin{equation}\label{defeord}
\begin{array}{rrcl}
\mbox{ord}^{}_{(J,b)}: & \Sing(J,b)& \to  & \mathbb{Q}_{\ \geq 1}
\\
 & x & \to & \mbox{ord}^{}_{(J,b)}(x)=\frac{\nu_x(J)}{b}.
\end{array}
\end{equation}

If $J$ is the defining ideal of a hypersurface $X\subset V$ as
before, we can think that the worst singularities of $X$ are located
at the points where  Hironaka's function is
maximum. Thus,  our aim is to  lower  the maximum value of this
function by defining a suitable sequence of blow-ups at smooth
centers in an effort to improve  the
singularities of  the strict transform of $X$.  }
\end{Paragraph}

\begin{Paragraph} {\bf Basic objects and resolution.}
{\rm Since our resolution problem has been codified in terms of
pairs, the next step is to understand how pairs transform under
blow-ups. Given a pair  $(J,b)$, a smooth closed subscheme
$Y\subset V$ is said to be {\em permissible} if  $Y \subset
\Sing(J,b)$. If  $ V \stackrel{\pi}{\longleftarrow} V_{1} \supset
H=\pi^{-1}(Y) $ denotes the monoidal transformation at a
permissible center $Y$ then the {\em total transform of $J$ in
${\mathcal O}_{V_1}$}, $J\calo_{V_1}$, can be expressed as a
product, $$J\calo_{V_1}= I(H)^b J_1$$ for a uniquely defined $J_1
$ in $ \calo_{V_1}$. The new couple $(J_1, b)$ is called the {\em
transform} of $(J, b)$, say:
\begin{equation}\label{ecnwc}
\begin{array}{ccccc}
 & V & \stackrel{\pi}{\longleftarrow} & V_{1} ,& \\ & (J,b) & & (J_1,b) &\\
\end{array}
\end{equation}
Observe that, in general, $J_1$ is strictly contained in the 
{\em weak transform of $J$ in ${\mathcal O}_{V_1}$} (the weak transform is 
defined as $\overline{J}_1=I(H)^{-\nu_Y(J)}\cdot J{\mathcal O}_{V_1}$). 
We refer to Example \ref{motivtrans} below for the motivation of this definition.

\ 

However,  some geometric conditions have to be imposed  in order to
define a {\em sequence of transformations} of a pair. Every monoidal
transformation introduces an exceptional divisor and we require that
these divisors have normal crossings. To keep track of this
additional information  we define a {\em couple} $(V,E)$ to be  a
smooth scheme $V$ together with a set of smooth hypersurfaces
$E=\{H_1, \dots , H_r\}$  so that their union has normal crossings.
If $Y$ is closed and smooth in $V$,  and has normal crossings with
$E$ (i.e., with the union of hypersurfaces of $E$),  then we define a {\em
transform of the couple}, say
$$(V,E)\leftarrow (V_1,E_1),$$  where $V\leftarrow V_1$ is the
blow-up at $Y$; and $E_1=\{H_1, \dots , H_r, H_{r+1}\}$, where
$H_{r+1}$ is the exceptional locus, and each $H_i$ denotes again
the strict transform of $H_i$ in $V_1$, for $1\leq i \leq r$.

\

We finally define a {\em basic object} to be a couple $(V,
E=\{H_1,\ldots, H_r\})$ together with a pair $(J,b)$, and we
denote it by
$$ (V,(J,b),E).$$ 
With this notation, $J$ comes with a factorization $J={\mathcal I}(H_1)^{\alpha_1}\cdots {\mathcal I}(H_r)^{\alpha_r} \overline{J}$ for suitable $\alpha_1,\ldots,\alpha_r\in {\mathbb N}$, and $\overline{J}\subset {\mathcal O}_V$.  We  say that $ (V,(J,b),E)$ is a {\em $d$-dimensional basic object} if
the dimension of $V$ is $d$. If a smooth center $Y$ defines a
transformation of $(V,E)$, and in addition $Y \subset \Sing(J,b)$,
then a transform of the couple $(J,b)$, say $(J_1,b)$, is defined
as above.
 In this case we say that $$ (V,(J,b),E) \longleftarrow
(V_1,(J_1,b),E_1)$$ is a {\em transformation} of the basic object.
So we will ask permissible centers to satisfy this normal crossings
condition.

\

 A {\em sequence of permissible transformations} is denoted by
\begin{equation}\label{transfuno}
(V,(J,b),E) \longleftarrow (V_1,(J_1,b),E_1)\longleftarrow \cdots
\longleftarrow (V_s,(J_s,b),E_s);
\end{equation}
 and such sequence is said to be a {\em resolution} of the basic object
 if
$\Sing(J_s,b)=\emptyset.$  }
\end{Paragraph}

\begin{Example} \label{examplerespairs}{\rm A resolution of $(V,(J,b),E)=(V,({\mathcal
I}(X),b),\{\emptyset\})$, with $b$  the maximum multiplicity of a
hypersurface $X$, lowers  the maximum multiplicity of
the strict transform of $X$ in $V_s$.}
\end{Example}

\begin{Example}  The  monomial case.  {\rm Let $(V,(J,b),E)$ be a
basic  object with $E=\{H_1,\ldots,H_l\}$. Notice that if
$J={\mathcal I}(H_1)^{a_1}\cdots {\mathcal I}(H_l)^{a_l}$ with
$a_i\in {\mathbb N}$  for  $i=1,\ldots,l$,  then it is relatively
easy to find a resolution of $(V,(J,b),E)$, which  can be achieved
using a combinatorial argument.}
\end{Example}

\begin{Example} \label{motivtrans}
{\rm Consider the surface in a three dimensional affine space, $X:=\{z^2+(x^2-y^3)^2=0\}\subset {\mathbb
A}_k^3$. Since its maximum order is $2$, we will want to resolve the basic object 
$({\mathbb A}_k^3, ({\mathcal I}(X),2),E^{(3)}=\{\emptyset\})$. 

Notice that $\mbox{Sing}(({\mathcal I}(X),2))$ is contained in the smooth surface $Z:=\{z=0\}\simeq {\mathbb A}^2$,  and that  there,  we can describe it as the singular locus of the pair $(\langle (x^2-y^3)^2\rangle, 2)=(J,2)$. 

Consider the blow-up at the origin of ${\mathbb A}_k^3$, $\pi: V^{(3)}\to  {\mathbb A}_k^3$ and denote the exceptional divisor by  $H$.  This also induces a blow-up $\overline{\pi}: V^{(2)}\to  {\mathbb A}_k^2$ with exceptional divisor $\overline{H}$.

The  strict transform of $X$, $X_1\subset V_1^{(3)}$, still has points of order 2, and  moreover  $$\mbox{Sing }({\mathcal I}(X_1),2)\subset Z_1,$$ where $Z_1\subset V_1^{(3)}$ is the strict transform of $Z$. A quick computation shows that this set is   $\mbox{Sing }({\mathcal I}(\overline{H})^{-2}J{\mathcal O}_{V^{(2)}}, 2)=(J_1,2)$. Notice that  $J_1$  is strictly  contained in the weak transform of $J$ in $V_1^{(3)}$, but from the way it is defined  there is a commutative diagram of 
restriction and transformations of basic objects: 

$$\begin{array}{ccc}
({\mathbb A}_k^3, ({\mathcal I}(X),2), E^{(3)}=\{\emptyset\}) & \leftarrow & (V^{(3)}, ({\mathcal I}(X_1),2), E^{(3)}=\{H\})  \\
\downarrow  & & \downarrow  \\
(Z, (\langle (x^2-y^3)^2\rangle, 2), E^{(2)}=\{\emptyset\})   & \leftarrow & (V^{(2)}, ({\mathcal I}(\overline{H})^{-2}J{\mathcal O}_{V^{(2)}}, 2), E^{(2)}=\{\overline{H}\}).\end{array} $$ 
With this law of transformation of pairs, it follows  that a resolution of the 2-dimensional basic object induces a resolution of the original in 3-dimensional space. This illustrates the fact  that resolution of basic objects is obtained by induction on  the dimension. }
\end{Example}

\section{Algorithmic resolution of basic
objects}\label{algoritmicpairs}

Given   a basic object $(V,(J,b),E)$,  algorithms for resolving singularities  provide a resolution as in (\ref{transfuno}), where
the choice of the centers of the monoidal transformations is given
by  the ``worst stratum''  which is defined by a suitable
upper semi-continuous function.

\

We distinguish two steps in   algorithmic resolution:
\begin{itemize}
\item {\bf Step A. Reduction to the mononial case. } In this step a sequence of permissible
transformations is  defined  to  {\em simplify } the structure
of the basic object.
\item {\bf Step B. Treatment of the monomial case.}  This step involves the resolution of a
basic object that is assumed to be within the monomial case.
\end{itemize}

\

Step A is accomplished by both defining a suitable
upper semi-continuous function constructed from the so called {\em
satellite functions}, and using  an inductive argument. In
Step B the monomial case is treated using an upper semi-continuous
function of combinatorial nature.

\

In the following paragraphs we sketch  how to accomplish  Step A;  
see \cite{EncVil97:Tirol} for
more details on this matter and a treatment of Step B.


\

\noindent {\bf Step A. Satellite functions}

\

Let $V$ be a $d$-dimensional smooth
scheme and consider a sequence of transformations of basic objects
which is not necessarily a resolution,
\begin{equation}\label{Atransfuno}
(V_0,(J_0,b),E_0)= (V,(J,b),E) \longleftarrow
(V_1,(J_1,b),E_1)\longleftarrow \cdots \longleftarrow
(V_s,(J_s,b),E_s).
\end{equation}
Let $\{H_{r+1}, \dots , H_{r+s}\}(\subset E_s)$ denote the
exceptional (irreducible) hypersurfaces introduced by the sequence of blow-ups. Satellite functions are upper-semi continuous functions defined at each step of   a sequence   like (\ref{Atransfuno}), that derive from Hironaka's order function. 

\begin{Paragraph} {\bf The first satellite function. }
\label{lasat1}   {\rm Given a sequence like \ref{Atransfuno} above, 
 there is a
well defined factorization of the sheaf of ideals $J_s \subset
\calo_{V_s}$:
\begin{equation}\label{eqfcword}
J_s=I(H_{r+1})^{b_{s,r+1}}I(H_{r+2})^{b_{s,r+2}}\cdots
I(H_{r+s})^{b_{s,r+s}}\cdot \overline{J}_s
\end{equation} such  that $\overline{J}_s$ does not vanish along
$H_{r+i}$ for  $0\leq i \leq s$.

\

Define  $\vword^{(d)}_{(J_s,b)}$ (or simply $\vword^{(d)}_{s}$):
\begin{equation}\label{eqdword}
\begin{array}{rrcl}
\vword^{(d)}_s: & \Sing(J_s,b)& \to & \mathbb{Q}
 \\  & x & \to & \vword^{(d)}_s(x)=\frac{\nu_x({\overline{J}_s})}{b},
\end{array}
\end{equation}
where $\nu_x(\overline{J}_s)$ denotes the order of
$\overline{J}_s$ at $\calo_{V_s,x}$.  This function has the
following properties:

 1) It  is
upper semi-continuous. In particular the set of points where it reaches its maximum value, 
$\mbox{max w-ord}^{(d)}_s$, is closed.
This set is denoted by $\vMax \vword^{(d)}_s$.

2) For any index $ i \leq s$, there is an expression
$$J_i=I(H_{r+1})^{b_{i,r+1}}\cdots I(H_{r+i})^{b_{i,r+i}}\cdot
\overline{J}_i,$$ and hence the function $\vword_i^{(d)}:
\Sing(J_i,b) \to {\mathbb Q}$ can also be defined.

3)  If each transformation of basic objects  $(V_i,(J_i,b_i), E_i)
\leftarrow (V_{i+1},(J_{i+1},b_{i+1} )E_{i+1})$ in
(\ref{Atransfuno}) is defined with center $Y_i \subset \vMax
\vword_i^{(d)}$, then \begin{equation}\label{satelitedecrece} \max
\vword^{(d)} \geq \max \vword_1^{(d)}\geq \dots \geq \max
\vword_s^{(d)}.\end{equation}

\

Observe that at the beginning  of a resolution process    $\vword^{(d)}_{0} : \Sing(J_0,b) \to \mathbb{Q}$  is the same as   $\mbox{ord}_{(J_0,b)}:  \Sing(J_0,b) \to \mathbb{Q}$
in (\ref{defeord}). For $i>0$,  the functions $\vword^{(d)}_{i}:  \Sing(J_i,b)  \to \mathbb{Q}$ 
 are defined on the weak transform of the ideal $J_{i-1}$ in $V_i$, and therefore for indices $i\geq 1$ they differ from Hironaka's order function, although they strongly depend on it. 
These  are the {\em first  satellite functions} of the function
introduced in (\ref{defeord}). They represent small variations of the
original function, and satisfy the inequalities stated
 in (\ref{satelitedecrece}).

\

Observe that $ \max \vword_s^{(d)}=0$ when
\begin{equation}\label{eqfcword1}
J_s=I(H_{r+1})^{b_{s,r+1}}I(H_{r+2})^{b_{s,r+2}}\cdots
I(H_{r+s})^{b_{s,r+s}},
\end{equation}
i.e., when $ \overline{J}_s=\calo_{V_s}$ in (\ref{eqfcword}); in
this  case we say that $(V_s, (J_s,b), E_s)$ is in the {\em monomial
case}. In the monomial case   it is easy to enlarge sequence
(\ref{Atransfuno})  to obtain a resolution.  Therefore  the
functions $\vword_i^{(d)}: \Sing(J_i,b) \to {\mathbb Q}$ measure
how far $J_i$ is from being a locally monomial sheaf of ideals   supported on  the exceptional locus. }
\end{Paragraph}

\begin{Paragraph}\label{lasat2}{\bf The second satellite functions.} {\rm
Consider a sequence of permissible transformations of
$d$-dimensional basic objects,
\begin{equation}\label{Atransfdos}
(V,(J,b),E) \longleftarrow (V_1,(J_1,b),E_1)\longleftarrow \cdots
\longleftarrow (V_s,(J_s,b),E_s),
\end{equation}
with
$$\mbox{max w-ord}^{(d)}\geq \mbox{max w-ord}_1^{(d)}\geq \ldots \mbox{max
w-ord}_s^{(d)}.$$ Then if $\mbox{max w-ord}_s^{(d)}>0$, the  function
$t_s^{(d)}$ is defined in the following way: let $s_0\leq s$ be the smallest index such
that
$$\mbox{max w-ord}^{(d)}\geq \mbox{max w-ord}_1^{(d)}\geq \ldots  \mbox{max
w-ord}_{s_0}^{(d)}=\mbox{max w-ord}_{s_0+1}^{(d)}=\ldots =\mbox{max
w-ord}_s^{(d)},$$ and set
$$E_s=E_s^+\sqcup E_s^-$$
where $E_s^-$ are the strict transforms of the hypersurfaces in
$E_{s_0}$. Then  we can define:
\begin{equation}\label{ect}
\begin{array}{rccl}
t_s^{(d)}: & \mbox{Sing }(J_s,b) & \longrightarrow & \left({\mathbb Q}\cup \{\infty\}\right) 
\times
{\mathbb N} \\
 & x & \to & (\mbox{w-ord}^{(d)}_s(x),n^{(d)}_s(x))
 \end{array}
 \end{equation}
where
$$n^{(d)}_s(x)=\sharp\{H_i\in E_s^{-}: x\in H_i\}$$ and  ${\mathbb Q}
\times
{\mathbb N}$ is lexicographically ordered.
The function $t_s^{(d)}$  is upper semi-continuous, and it is
designed to ensure the normal crossings condition
of the permissible centers with the smooth hypersurfaces in $E_s$.


}
\end{Paragraph}

\begin{Remark}{\rm The
 functions $t_i^{(d)}$ depend  on Hironaka's order function (see
(\ref{defeord}) in \ref{unouno}). These functions are referred to as 
{\em inductive functions} since they play a key role in the inductive arguments 
of Step A.   }
\end{Remark}

\

\noindent {\bf Step A. Induction and maximal contact}

\

In general satellite functions
 are too coarse to provide, just by themselves, an
upper semi-continuous function that leads  to resolution, or  to
the monomial case. For instance, it can be easily seen that for
$({\mathbb A}_k^3, (\langle z^2+(x^2-y^3)^2 \rangle,
2),E^{(3)}=\{\emptyset\})$ the maximum of $t^{(3)}$ is not
smooth so it does not define a permissible center.  Induction is thus used to solve this problem: the information provided by the satellite
functions is refined using an inductive argument as explained in the
next paragraphs.

\begin{Paragraph}{\bf Simple basic objects and induction.} \label{racsim} {\rm
Let $V$ be a smooth scheme. A pair $(J,b)$ is said to be {\em
simple } if $\mbox{ord}^{}_{(J,b)}:\Sing(J,b) \to \mathbb{Q}$ \ is
  constant with $\mbox{ord}_{(J,b)}(x)=1$ for all $x\in \mbox{Sing }(J,b)$; namely when the order of $J$ is exactly
$b$ at the local ring $ \calo_{V,x}$ for any  $ x\in\Sing(J,b)$. A basic
object $(V,(J,b),E)$ is said to be a {\em simple basic object}
when $(J,b)$ is simple.

\

In the case of characteristic zero,  the  resolution of simple basic
objects can be defined if we assume, by induction,  the resolution of
basic objects on lower dimensional ambient spaces.  This   is
guaranteed by the notion of {\em maximal contact}: a
$d$-dimensional simple basic object $(V^{(d)},(J,b),E^{(d)})$  can
be restricted, locally, to a smooth hypersurface,  defining a
$(d-1)$-dimensional basic object on this smooth lower dimensional
space, $(V^{(d-1)},(J^{\prime},b^{\prime}),E^{(d-1)})$.
Furthermore, the link between the original basic object and the
restricted one is sufficiently  strong so that a resolution of the
latter induces a resolution of the former, since  there are
commutative diagrams of transformations and restrictions:
$$\begin{array}{ccccccc}
(V^{(d)},(J,b),E^{(d)}) & \leftarrow &
(V_1^{(d)},(J_1,b),E_1^{(d)}) & \leftarrow
& \cdots & \leftarrow & (V_s^{(d)},(J_s,b),E_s^{(d)})\\
\downarrow & & \downarrow & & \cdots & & \downarrow \\
(V^{(d-1)},(J^{\prime},b^{\prime}),E^{(d-1)}) & \leftarrow &
(V_1^{(d-1)},(J_1^{\prime},b^{\prime}),E_1^{(d-1)}) & \leftarrow &
\cdots & \leftarrow &
(V_s^{(d-1)},(J_s^{\prime},b^{\prime}),E_s^{(d-1)}).
\end{array}$$

In other words, {\bf simple basic objects can be resolved by induction}:
starting from a simple basic object, a new basic object can be 
defined in lower dimension. A resolution of the latter  can also be defined in terms 
of satellite functions in lower dimension. 
This in turns defines a resolution of the 
original simple basic object. }
\end{Paragraph}

\begin{Paragraph} {\bf The non-simple case.}  {\rm  The previous discussion shows how induction comes in when dealing with simple basic objects, but {\bf only} when dealing with simple basic objects.  Even this form of induction will lead us to the non-simple case.   For instance, a resolution of the 3-dimensional simple basic object  $({\mathbb A}_k^3, (\langle z^2+(x^2-y^3)^2 \rangle,
2),E^{(3)}=\{\emptyset\})$ can be obtained 
by finding a resolution of the 2-dimensional basic object $({\mathbb A}_k^2, (\langle (x^2-y^3)^2 \rangle,
2),E^{(2)}=\{\emptyset\})$ or equivalently of   $({\mathbb A}_k^2, (\langle (x^2-y^3) \rangle,
1),E^{(2)}=\{\emptyset\})$, but  the latter 
 is not simple. Notice that $\mbox{Sing} (\langle (x^2-y^3) \rangle,
1)$ is not contained in a smooth one dimensional scheme. However,  there is a 2-dimensional 
simple basic object naturally   attached to  it: $({\mathbb A}_k^2, (\langle (x^2-y^3) \rangle,
2),E^{(2)}=\{\emptyset\})$. Observe that 
$$\mbox{Sing }(\langle (x^2-y^3) \rangle,
2)=\underline{\mbox{Max}}-\vword^{(2)}(\langle (x^2-y^3) \rangle,
1),$$
and that:  \begin{enumerate} 
\item[i.] A resolution of $({\mathbb A}_k^2, (\langle (x^2-y^3) \rangle,
2),E^{(2)}=\{\emptyset\})$ induces a lowering of the maximum  value of 
$\vword^{(2)}$ on $({\mathbb A}_k^2, (\langle (x^2-y^3) \rangle,
1),E^{(2)}=\{\emptyset\})$; 
\item[ii.] Since  $({\mathbb A}_k^2, (\langle (x^2-y^3) \rangle,
2),E^{(2)}=\{\emptyset\})$ is simple, 
a resolution can be found by means of an inductive argument, now in dimension 
one. 
\end{enumerate}

So, one property of the first satellite function is that it is
naturally attached to a simple basic object. Given a
non-necessarily simple basic object, $(V^{(d)},(J,b),E^{(d)})$,
there is a simple basic object attached to it, say 
$(V^{(d)},(\widetilde{J},\widetilde{b}),E^{(d)})$,  such that
$$\mbox{Sing }(\widetilde{J},\widetilde{b})=\underline{\mbox{Max}}-\vword^{(d)} (J,b).$$
And moreover,    a resolution of $(V^{(d)},(\widetilde{J},\widetilde{b}),E^{(d)})$, say
$$(V^{(d)},(\widetilde{J},\widetilde{b}),E^{(d)}) \leftarrow (V^{(d)}_1,(\widetilde{J}_1,\widetilde{b}),E^{(d)}_1) \leftarrow \ldots \leftarrow 
(V^{(d)}_s,(\widetilde{J}_s,\widetilde{b}),E^{(d)}_s)$$
where $\mbox{Sing }(\widetilde{J}_s,\widetilde{b})=\{\emptyset\}$,  
induces a sequence of
 permissible transformations of $(V^{(d)},(J,b),E^{(d)})$,  say 
$$(V^{(d)},(J,b),E^{(d)})   \leftarrow 
(V_1^{(d)},(J_1,b),E_1^{(d)})   \leftarrow
  \cdots   \leftarrow   (V_s^{(d)},(J_s,b),E_s^{(d)}),$$
that lowers the maximum value 
  $\max
 \vword^{(d)}(J,b)$, i.e., 
$$ \max \vword^{(d)} = \max \vword_1^{(d)}= \dots > \max
\vword_s^{(d)}.$$
As a consequence,  by  successively resolving  the simple
basic objects
 attached to the functions $\vword^{(d)}_{i}$ we produce a sequence of transformations 
where  $$ \max \vword^{(d)} \geq \max \vword_1^{(d)}\geq \dots \geq \max
\vword_s^{(d)}=0;$$ which  implies  that $J_s$
is monomial. 

\

We should estress here the role played by the function $t^{(d)}$.  
If $(J,b)$ is simple at a point $x$, and if  the codimension of $\mbox{Sing }(J,b)$ 
is one in a neighborhood of $x$, then it can be shown that 
$\underline{\mbox{Max}}\vword^{(d)}$ is smooth locally at $x$, and 
hence it is a 
canonical center  to blow-up. As a consequence, the function $t^{(d)}$ assigns the value 
``$\infty$" to these points. If $(J,b)$ is simple at $x$, and if the codimension of $\mbox{Sing }(J,b)$ is greater than one, then via $t^{(d)}$ a $(d-1)$-dimensional basic object can be  attached to $(V,(J,b),E)$ as explained above.}
\end{Paragraph}
\begin{Paragraph}\label{technical} {\bf Step A.  Technical problems.}  {\rm There are three main sub-steps in step A: the local restriction to hypersurfaces of maximal contact, commutative diagrams of restrictions and blow-ups, and  the attachment of simple basic
objects to non-simple basic objects.

$\bullet$ \underline{Maximal contact.} For a fixed simple basic
object there may be different choices of hypersurfaces of maximal
contact. However it can be shown  that they all lead to the same
resolution. This result is the main outcome of the so called
Hironaka's trick; an alternative and enlightening proof of this result
is given  by  J. W\l odarczyk (see \cite{WLL}).

$\bullet$ \underline{Commutative diagrams of restrictions and
permissible transformations.} Step A is accomplished  by an
inductive argument;  for this argument to hold  it is required  that restrictions and permissible transformations commute  as
indicated in \ref{racsim}.

$\bullet$ \underline{Association of simple basic objects to
non-simple basic objects.} For a fixed basic object there may be
different choices of simple basic objects that can be associated
to it. It can be shown that all different choices lead to the same
resolution.}

\end{Paragraph}

 \begin{Paragraph}
 {\bf Summarizing.} {\rm Given a basic object $(V,(J,b),E)$, an upper semi-continuous
function is defined by means of the satellite
 functions in dimension $d$ and lower dimensions: if $x\in
 \mbox{Sing}(J,b)$, then the upper continuous-function associates
 to it a set of values
 $$t^{(d)}(x), t^{(d-1)}(x), \ldots, t^{(d-r)}(x) \in ({\mathbb Q}\cup \{\infty\})\times {\mathbb N}$$  
for some $r< d$.  This is done by identifying $x$ with the points in  
 successive restrictions to hypersurfaces of maximal contact. Then a resolution of 
$(V,(J,b),E)$  is achieved by blowing up the centers defined by
these functions. }
\end{Paragraph}

\begin{Example} {\rm  To find a resolution of singularities of
 $X:=\{z^2+(x^2-y^3)^2=0\}\subset {\mathbb
A}_k^3$, we start by finding a resolution of  the basic object
$({\mathbb A}_k^3, ({\mathcal I}(X),2),E^{(3)}=\{\emptyset\})$.
Since
$$\mbox{Sing }({\mathcal I}(X),2) = \{z=0, x^2-y^3=0\}=C,$$ we
can take $\{z=0\}\simeq {\mathbb A}_k^2$ as a hypersurface of
maximal contact. We then   associate to the original 3-dimensional basic
object a 2-dimensional one: $$({\mathbb A}_k^2,(\langle
(x^2-y^3)^2\rangle,2), E^{(2)}=\{\emptyset\}).$$ This basic
object  is not simple, so we attach to $\mbox{{Max}
w-ord}^{(2)}$ a simple basic object (s.b.o.), $({\mathbb A}^2_k,
(\langle x^2-y^3\rangle, 2),\{\emptyset\})$, and  then find another
hypersurface of maximal contact, $\{x=0\}\simeq {\mathbb A}_k^1$,
and a 1-dimensional basic object:

$$\begin{array}{cccclcc}
({\mathbb A}^3_k, (\langle z^2+(x^2-y^3)^2\rangle,2),
\{\emptyset\})  &  & &  & & & \\
 \mbox{\tiny{\textsf{Restriction}}} \downarrow &  \mbox{\tiny{\textsf{s. b. o.}}} & &  & & & \\
({\mathbb A}^2_k, (\langle (x^2-y^3)^2\rangle, 2), \{\emptyset\})
 & \leftrightarrow  & ({\mathbb A}^2_k, (\langle
x^2-y^3\rangle, 2),\{\emptyset\}) & & & & \\
 & &  \downarrow \mbox{\tiny{\textsf{Restriction}}} & \mbox{\tiny{\textsf{s. b. o.}}} &  & &    \\
& &  ({\mathbb A}^1_k, (\langle y^3\rangle, 2),\{\emptyset\}) &
\leftrightarrow  & ({\mathbb A}^1_k, (\langle y^3\rangle,
3),\{\emptyset\}) & &
\end{array}$$

\

The information can be interpreted  in the following way:

\

- A resolution of the simple basic object $({\mathbb A}^3_k, (\langle
z^2+(x^2-y^3)^2\rangle,2), \{\emptyset\})$ can be constructed  by finding
a resolution of $({\mathbb A}^2_k, (\langle (x^2-y^3)^2\rangle,
2), \{\emptyset\})$.

\

- Lowering the maximum order of $\langle (x^2-y^3)^2\rangle$ in
${\mathbb A}^2_k$ (reached at  $(0,0)$) is equivalent to
finding a resolution of the simple basic object $({\mathbb A}^2_k, (\langle
x^2-y^3\rangle, 2),\{\emptyset\})$.

\

-  A resolution of $({\mathbb A}^2_k, (\langle x^2-y^3\rangle,
2),\{\emptyset\})$   can be constructed  by resolving $({\mathbb A}^1_k,
(\langle y^3\rangle, 2),\{\emptyset\})$.

\

- The  maximum order of  $\langle y^3\rangle $ in ${\mathbb
A}^1_k$ is forced to drop by resolving $({\mathbb A}^1_k, (\langle
y^3\rangle, 3),\{\emptyset\})$.

\

By collecting the information provided by the order function in
different dimensions,  the following upper semi-continuous
function is defined:
$$\Gamma_X(p)=\left\{\begin{array}{lcl}
((1,0),(2,0), (\frac{3}{2},0)) &  \mbox{ if } & p=(0,0,0)\\
((1,0),(1,0), (\infty,0) &  \mbox{ if }   & p\in C\setminus
{(0,0,0)}.
\end{array}\right.$$
Its  maximum value indicates the first center that has to be blown up:
$(0,0,0)$. 

\

After a blow-up at the origin, the maximum of the w-order in the
second level, $\mbox{max w-ord}^{(2)}_1$, has dropped, so  the
function $n^{(2)}_1$ plays a role in counting old exceptional divisors.
 The sequence defined by the algorithm takes the following form:

\

\mbox{{\textsf{Starting point: 3-dimensional basic object}}}\hrule

$$\begin{array}{lccccc}
\mbox{{\textsf{Couples:}}} & ({\mathbb A}^3_k, \{\emptyset\}) &
\leftarrow & (V_1^{(3)}, \{H_1\}) & \leftarrow & (V_2^{(3)},
  \{H_1, H_2\})\\
 & & & & & \\
\mbox{{\textsf{Pairs:}}} &  (\langle z^2+(x^2-y^3)^2\rangle,2) &
\leftarrow & (\langle z_1^2+I(H_1)^2(x^2_1-y_1)^2 \rangle,2) &
\leftarrow & (\langle z_2^2+I(H_1)^2I(H_2)^2(1-y_2)^2\rangle,2)
\end{array}$$

\

\mbox{{\textsf{Restricting: 2-dimensional basic object}}}\hrule

$$
\begin{array}{lccccc}
\mbox{{\textsf{Couples:}}}  & ({\mathbb A}^2_k, \{\emptyset\}) &
\leftarrow & (V_1^{(2)}, \{\overline{H}_1\})& \leftarrow &
(V_2^{(2)},
 \{\overline{H}_1, \overline{H}_2\})\\
 & & & & & \\
\mbox{{\textsf{Pairs:}}} & (\langle (x^2-y^3)^2\rangle, 2) &
\leftarrow & (\langle I(\overline{H}_1)^2(x^2_1-y_1)^2\rangle,2) & \leftarrow &
(\langle I(\overline{H}_1)^2I(\overline{H}_2)^2(1-y_2)^2\rangle,2).
\end{array}$$

After two permissible transformations the Step A of resolution has
been accomplished since we have reached the monomial case in
dimension 2. Notice  that  a resolution of $({\mathbb
A}_k^2,(x^2-y^3)^2, E^{(2)}=\{\emptyset\})$ induces a resolution
of the original 3-dimensional basic object. }

\end{Example}

\begin{Paragraph} {\bf About this paper.}
{\rm The purpose of this paper is to show that Step A of the
resolution process is characteristic free, i.e., this part of the algorithmic resolution can be performed in any characteristic if we replace the restriction to  hypersurfaces of maximal contact with a different form of induction that requires {\em projections}. This
new approach  is formulated in terms of Rees algebras  instead of pairs.
There is a dictionary that translates between pairs and Rees algebras (see
\ref{paraleloragadprs}, \ref{ragadprs3} and \ref{rgn}). Their role in resolution
problems will be explained in  Section \ref{sectionorder}. When translating the
algorithmic resolution to this new setting we encounter the corresponding technical problems as described in \ref{technical}:

$\bullet$ \underline{Projections to smooth schemes.} For a
fixed simple Rees algebra there may be numerous suitable  projections. In Main Theorem \ref{ordertau} we
show that the w-ord-functions defined after projecting are
independent of the choice of the projections. This leads to
Definition  \ref{deford}, and, consequently,
satellite functions can also be defined in our context.

$\bullet$ \underline{Commutative diagrams of projections and
permissible transformations.} This is addressed  in Section
\ref{eliminationmonoidal}.

$\bullet$ \underline{Association of simple Rees algebras  to
non-simple Rees algebras.} For a fixed Rees algebra there may be
different choices of simple Rees algebras that can be associated
to it. In Section \ref{nonsimple} we show that they are canonical choices, and in 
particular they automatically globalize.

\

Using these results, it follows   that the so called ``reduction to
the monomial case" is possible in positive characteristic, meaning that  the
monomial case arises in some lower dimension via
induction (cf. \cite{positive}).}
\end{Paragraph}



\part{Rees algebras}\label{partRees}

\section{Rees algebras}\label{Rees}
We begin by introducing Rees algebras, an essential tool for the 
outcome of this paper. Special attention should be
paid to  Example  \ref{ExampleA}  where Rees algebras are studied in the typical situations that we are interested in.

\begin{Definition} \label{Reesalg}{\rm Let $B$ be a Noetherian ring, and let
$\{I_n\}_{n\geq 0}$ be a sequence of ideals in $B$ satisfying the
following conditions:
\begin{enumerate}
\item[i.] $I_0=B$;
\item[ii.] $I_k\cdot I_l\subset I_{k+l}$.
\end{enumerate}
Then the graded subring ${\mathcal G}=\oplus_{n\geq 0}I_nW^n$ of
the polynomial ring $B[W]$ is said to be a {\em Rees algebra} if
it is a finitely generated $B$-algebra. }
\end{Definition}
\begin{Remark}{\rm A Rees algebra can be described by
giving a finite set of   generators
$$\{f_{n_1}W^{n_1},\ldots,f_{n_s}W^{n_s}\}$$ with $f_{n_i}\in B$ for $i=1\ldots,s$.  An element
$g\in I_n$ will be of the form $g=F_n(f_{n_1},\ldots,f_{n_s})$ for
some weighted homogeneous polynomial in $s$-variables
$F_n(Y_1,\ldots,Y_s)$ where $Y_i$ has weight $n_i$ for
$i=1,\ldots,s$.}
\end{Remark}

\begin{Example} {\rm A canonical example of a Rees algebra is the  Rees
ring of an ideal: fix an ideal $J\subset B$, and let ${\mathcal
G}=\oplus_nJ^nW^n$. In fact, a Rees algebra is  not very far away
from being the Rees ring of an ideal  in a sense that we make precise
in the following lines. Let ${\mathcal G}=\bigoplus_{n\geq
0}I_nW^{n}\subset B[W]$ be the Rees algebra generated by
$\{f_{n_1}W^{n_1},\ldots,f_{n_s}W^s\}$ with $f_i\in B$, and let
$N$ be a common multiple of all integers $n_i$, $i=1,\ldots,s$.
Then
$$\bigoplus_{k\geq 0}I_N^kW^{kN}\subset \bigoplus_{n\geq 0}I_nW^{n}$$
is a finite extension of Rees algebras (cf.
\cite[2.3]{integraldifferential}). So, up to integral closure, a
Rees algebra can be thought of as the Rees ring of a suitable
ideal (see \cite{hpositive}). }
\end{Example}

\begin{Remark} {\rm Given a Rees algebra ${\mathcal G}=\oplus_{n\geq
0}I_nW^n$ another Rees algebra can be defined by setting
$$I_n^{\prime}=\sum_{r\geq n}I_r,$$
and letting ${\mathcal L}=\oplus_{n\geq 0}I_n^{\prime}W^n$. Then
${\mathcal L}$ is contained in the integral closure of ${\mathcal
G}$ (cf. \cite[Remark 2.2 (2)]{hpositive}), and has the additional
property that $I^{\prime}_k\supset I^{\prime}_s$ if $s\geq k$. So,
up to integral closure it can always be assumed that a Rees
algebra ${\mathcal G}=\oplus_{n\geq
0}I_nW^n$ fulfills the additional condition that $I_n\supset 
I_m$ if $m\geq n$.}
\end{Remark}

\begin{Paragraph}{\bf Rees algebras on schemes.} {\rm Let $V$ be a scheme and let $\{I_n\}_{n\geq 0}$ be a
sequence of sheaves of ideals in ${\mathcal O}_V$ with
$I_0={\mathcal O}_V$ and such that $I_k\cdot I_l\subset I_{k+l}$
for all non-negative integers  $k,l$. The graded subsheaf of
algebras ${\mathcal G}=\oplus_{n\geq 0}I_nW^n$ of ${\mathcal O}_V[W]$
is said to be a {\em sheaf of Rees algebras}, or simply a {\em 
Rees algebra  on 
$V$} for short, if there is an affine
open cover $\{U_i\}$  of $V$,  such that ${\mathcal G}(U_i)\subset
{\mathcal O}_V(U_i)[W]$ is a Rees  ${\mathcal O}_V(U_i)$-algebra
by Definition \ref{Reesalg}.}
\end{Paragraph}

\begin{Paragraph} \label{singularlocus}{\bf The singular locus of a Rees algebra.} {\rm Let $V$ be a non-singular scheme   and let
${\mathcal G}=\oplus_nI_nW^n$ be a sheaf of Rees algebras. Let
$\nu_x(J)$ denote the order of an ideal $J$ in the regular local
ring ${\mathcal O}_{V,x}$. The {\em singular locus of ${\mathcal
G}$}, denoted by $\mbox{Sing }{\mathcal G}$, is the closed set
 of all points $x\in V$ such that $\nu_x(I_n)\geq n$ for
all non-negative integers $n$, i.e.,
$$\mbox{Sing }{\mathcal G}=\bigcap_n\{x\in V: \nu_x(I_n)\geq n, \mbox{ for all } n\in {\mathbb Z}_{\geq 0}\}.$$}
\end{Paragraph}

\begin{Example} \label{ExampleA} {\rm Let $\langle f\rangle\subset {\mathcal O}_V$ be the
ideal of an affine  hypersurface $H$ in an affine smooth scheme
$V$. Also let $b$ be a non-negative integer, and let ${\mathcal G}$ be
the Rees algebra generated by $f$ in degree $b$, i.e., ${\mathcal G}={\mathcal O}_V[fW^b]$. 
Then  $\mbox{Sing
}{\mathcal G}$ is the closed set of points of multiplicity at
least $b$ of $H$ (this may be empty). The same holds  if
$J\subset {\mathcal O}_V$ is a sheaf of ideals, $b$ is a
non-negative integer and ${\mathcal G}$ is the Rees algebra
generated by $J$ in degree $b$, namely ${\mathcal G}={\mathcal O}_V[JW^b]$.  Then
  the singular locus of
  ${\mathcal G}$ consists of
the points of $V$ where the order of $J$ is at least $b$ (which may
be empty).}
\end{Example}

\begin{Paragraph}{\bf Singular locus  and integral closure.}
{\rm  Two Rees algebras with the same integral closure have the same singular locus. In other 
words, 
if  ${\mathcal G}_1,{\mathcal G}_2\subset
{\mathcal O}_V[W]$ have the same integral closure in ${\mathcal O}_V[W]$, then $\mbox{Sing
}{\mathcal G}_1=\mbox{Sing }{\mathcal G}_2$ (see \cite[Proposition
4.4 (1)]{integraldifferential}. Hence the singular locus of a Rees algebra  is defined up 
to integral closure. }
\end{Paragraph}

\begin{Paragraph} \label{orderRees}{\bf The order of a Rees algebra at a point.}
  {\rm Let $x\in \mbox{Sing }{\mathcal
G}=\bigoplus_{n\geq 0}I_nW^n$, and let $fW^n\in I_nW^n$. Then set
$$\mbox{ord}_x(f)=\frac{\nu_x(f)}{n}\in {\mathbb Q},$$
where $\nu_x(f)$ denotes the order of $f$ in the regular local ring
${\mathcal O}_{V,x}$. Notice that $\mbox{ord}_x(f)\geq 1$ since
$x\in \mbox{Sing }{\mathcal G}$. Now define
$$\mbox{ord}_x{\mathcal
G}=\mbox{inf}\{\mbox{ord}_x(f):fW^n\in I_nW^n, n\geq 1\}.$$ If
${\mathcal G}$ is generated by
$\{f_{n_1}W^{n_1},\ldots,f_{n_m}W^{n_m}\}$ then
$$\mbox{ord}_x{\mathcal G}=\mbox{min}\{\mbox{ord}_x(f_{n_i}):
i=1,\ldots,m\},$$ and therefore, if $x\in \mbox{Sing }{\mathcal
G}$ then $\mbox{ord}_x{\mathcal G}$ is a rational number that is greater
or equal to one. Furthermore if $N$  is a common multiple of all
$n_i$, then
$$\mbox{ord}_x{\mathcal G}=\frac{\nu_x(I_N)}{N}.$$
If ${\mathcal G}_1,{\mathcal G}_2\subset {\mathcal
O}_V[W]$ have the same integral closure, then
$\mbox{ord}_x{\mathcal G}_1=\mbox{ord}_x{\mathcal G}_2$ at any
point  $x\in \mbox{Sing }{\mathcal G}_1=\mbox{Sing }{\mathcal
G}_2$ (cf. \cite[Proposition 6.4]{EV}).}

\end{Paragraph}

\begin{Paragraph}\label{paraleloragadprs} {\bf Rees algebras vs. pairs.}
 {\rm The information encoded in a Rees algebra is  
essentially the same as the one encoded by  Hironaka's notion of {\em pair} (see \cite{Hironaka77}).
  We assign to a pair $(J,b)$ over a smooth scheme $V$ the Rees
algebra:
\begin{equation}\label{eqq34}
\mathcal{G}_{(J,b)}=\calo_V[J^bW^b],
\end{equation}
which is a graded subalgebra in $\calo_V[W]$. It turns out that
every Rees algebra over $V$ is a finite extension of
$\mathcal{G}_{(J,b)}$ for a suitable pair $(J,b)$
(see \cite[Proposition 2.9]{positive} for details).

\

Observe  that for $\mathcal{G}_{(J,b)}=\calo_V[J^bW^b]$ there is
an equality of closed sets
$$\Sing(\mathcal{G}_{(J,b)})=\Sing(J,b),$$ and also of functions
$$\mbox{ord}^{}_{\mathcal{G}_{(J,b)}}=\mbox{ord}^{}_{(J,b)},$$
where the left-hand side is that defined in \ref{orderRees}. 

\

Hence,  up to integral closure, any Rees algebra is equivalent to a pair, and a resolution of the latter is equivalent to a resolution of the former (see \ref{ragadprs3} and \ref{rgn},  and \cite{mariluz}). }
\end{Paragraph}

\section{Differential Rees algebras}\label{DiffRees}

As was indicated in the previous section (see Example \ref{ExampleA})
we are particularly interested in the  multiplicity of
embedded hypersurfaces. For this purpose, we will use
  a class of Rees algebras  that are, in a sense, compatible
with differential operators. This point will be clarified in
\ref{difsingular}.

\

Let $V$ be a smooth scheme  over a field $k$. Then, for any
non-negative integer $s$, the sheaf of $k$-differential operators
of order $s$, $Diff^s_k$, is a coherent sheaf locally free over
$V$. If $s=0$, the sheaf $Diff^0_k$ can be naturally identified
with ${\mathcal O}_V$ and for each $s\geq 0$ there are natural
inclusions $Diff_k^s\subset Diff_k^{s+1}$.

\begin{Definition}\label{DiffAlg}{\rm A   Rees algebra ${\mathcal G}=\oplus_nI_nW^n$
is said to be a {\em differential Rees algebra},
 a {\em differential Rees algebra relative to $k$} or an 
{\em absolute differential Rees algebra} if the
following conditions hold:
\begin{enumerate}
\item[i.] For all non-negative integers $n$ there is an inclusion
$I_n\supset I_{n+1}$.
\item[ii.] There is an affine open covering of $V$, $\{U_i\}$, such
that for any $D\in Diff^r_k(U_i)$ and any $h\in I_n(U_i)$ we have
that $D(h)\in I_{n-r}(U_i)$ provided that $n\geq r$.
\end{enumerate} }
\end{Definition}

\begin{Paragraph}\label{extensiondif} {\bf The differential Rees 
algebra generated by a Rees algebra.} {\rm Let ${\mathcal G}$ be a
Rees algebra on a smooth scheme $V$  over a field $k$. There is a
natural way to construct a   differential Rees  algebra
 containing
 ${\mathcal G}$ with the property of being the smallest differential Rees  algebra
 containing it. This Rees algebra will be denoted
 by $\mathbb{D}\mbox{iff}({\mathcal G})$. If ${\mathcal G}$ is locally generated
 on an affine open set $U$ by $\{f_{n_1}W^{n_1},\ldots,f_{n_s}W^s\}$,
 then it can be shown that
 ${\mathbb D}\mbox{iff}({\mathcal G}(U))$ is generated by
 $$\{D(f_{n_i})W^{n_i^{\prime}-r}: D\in Diff^r_k, 
 0\leq r<n_i^{\prime}\leq n_i, i=1,\ldots,s\},$$
(see  \cite[Theorem 3.4]{integraldifferential}).}
\end{Paragraph}

\begin{Paragraph}\label{difsingular}{\bf Differential Rees algebras and singular locus.}
{\rm On a smooth scheme  $V$, of finite type over a perfect field $k$,
the sheaves of differentials $Diff^r_k$ for different values of
$r$  allow us  to study the order of a sheaf of
ideals. Similarly, differential Rees algebras are the right
structures for studying the singular locus of a Rees algebra. More
precisely, given a Rees algebra ${\mathcal
G}=\bigoplus_nI_nW^n$ on $V$,
$$\mbox{Sing }{\mathcal G}=\cap_{r\geq 0}V(Diff^{r-1}_k(I_r)),$$
(see \cite[Definition 4.2]{integraldifferential}). This definition
coincides with the one given  in Definition \ref{singularlocus} (see
\cite[Proposition 4.4]{integraldifferential}). In fact if ${\mathbb
D}\mbox{iff}({\mathcal G})$ is the differential Rees algebra generated by
a Rees algebra ${\mathcal G}$ then
$$\mbox{Sing }{\mathcal
G}=\mbox{Sing }{\mathbb D}\mbox{iff}({\mathcal G});$$ also if $x\in
\mbox{Sing }{\mathcal G}=\mbox{Sing }{\mathbb
D}\mbox{iff}({\mathcal G})$ then
$$\mbox{ord}_x{\mathcal
G}=\mbox{ord}_x{\mathbb D}\mbox{iff}({\mathcal G})$$ (cf.
\cite[Proposition 6.4]{EV}). Furthermore, if  ${\mathcal G}$ is a
differential Rees algebra, then $\mbox{Sing }{\mathcal G}=V(I_r)$ for any
positive integer $r$ (see \cite[Proposition
4.4]{integraldifferential}).}
\end{Paragraph}

\begin{Paragraph} \label{difintegral}{\bf Differential Rees algebras and integral
closure.} {\rm In many problems concerning resolution of
singularities it is natural to consider ideals up to integral
closure. For instance two ideals with the same integral closure
have the same embedded principalizations (Log-resolutions). In the
use of differential Rees  algebras as a tool to
understand singularities, we need to consider algebras up to
integral closure, so we need to understand how integral closure relates to
differential Rees algebras. This issue is treated in \cite[Section
6]{integraldifferential} where it is proven that if ${\mathcal
G}_1\subset {\mathcal G}_2$ is a finite extension of differential
algebras on a smooth scheme $V$  over a field $k$, then ${\mathbb
D}\mbox{iff}({\mathcal G}_1)\subset {\mathbb
D}\mbox{iff}({\mathcal G}_2)$ is also a finite extension. In other
words, if ${\mathcal G}_1$ is equal to ${\mathcal G}_2$ up to
integral closure, then so are ${\mathbb D}\mbox{iff}({\mathcal
G}_1)$ and ${\mathbb D}\mbox{iff}({\mathcal G}_2)$.}
\end{Paragraph}

\noindent {\bf Relative Differential Rees Algebras}

\begin{Paragraph} \label{difrel}
{ \em Let $\phi: V^{(d)} \to V^{(e)}$  be a smooth morphism of
smooth schemes of dimensions $d$ and $e$ respectively. Then, for any
non-negative integer $s$, the sheaf of relative differential
operators of order $s$,  $Diff^s( V^{(d)} /V^{(e)})$, is   locally
free over $V^{(d)} $. }
\end{Paragraph}

\begin{Definition}\label{def82}{\rm Let $\phi: V^{(d)} \to V^{(e)}$  be a smooth morphism of
smooth schemes of dimensions $d$ and $e$ respectively. A   Rees
algebra ${\mathcal G}=\oplus_nI_nW^n\subset {\mathcal
O}_{V^{(d)}}[W]$ is said to be a {\em $\phi$-relative differential
algebra}
 or simply a  {\em
$\phi$-differential Rees algebra } if:
\begin{enumerate}
\item[i.] For all non-negative integers $n$ there is an inclusion
$I_n\supset I_{n+1}$.
\item[ii.] There is an affine open covering $\{U_i\}$ of $V^{(d)} $   such
that for any $D\in Diff^s( V^{(d)} /V^{(e)})(U_i)$ and any $h\in
I_n(U_i)$, $D(h)\in I_{n-s}(U_i)$ provided that $n\geq
s$. \end{enumerate} }
\end{Definition}

Relative differential Rees algebras will play a central role in our
arguments  due to their relation to  a form of elimination that we
shall discuss in the next sections.
The case of relative dimension one, $V^{(d)} \to V^{(d-1)}$, is of particular interest.

\section{Rees algebras, permissible
transformations and weak equivalence}\label{transformationsweak}
In the previous section we attached to a Rees algebra a closed
set, its singular locus, and a function along this closed set: the order of a 
Rees algebra at a point. 
 The purpose of this section is to introduce the concept of {\em
weak equivalence}. Two Rees algebras with the same integral closure
will be weakly equivalent; a Rees algebra and the differential 
Rees algebra expanded by it will be weakly equivalent too. To this end, we consider three kinds of
transformations of Rees algebras: monoidal transformations,
restrictions to open sets, and products of smooth schemes with affine
spaces. These will be used to define the equivalence relation. If
two Rees algebras are equivalent according to this relation, then
they will have the same {\em resolution} (this concept  will be defined in
the next sections). In fact, within an equivalence class of Rees
algebras there is a natural procedure to choose one up to integral
closure (see
 Theorem \ref{examplesweak}). This equivalence relation will play a role in
Section \ref{nonsimple}.

\begin{Paragraph} \label{monoidal} {\bf Monoidal transformations.}
{\rm Let ${\mathcal G}=\oplus_{n}J_nW^n\subset{\mathcal O}_{V}[W]$
be a Rees algebra. A monoidal transformation with center $Y\subset
V$,  $V\leftarrow V^{\prime}$,  is said to be   {\em permissible} if
$Y\subset \mbox{Sing }{\mathcal G}$ is a smooth closed subscheme. If
$H\subset V^{\prime}$ is the exceptional divisor, then for each
$n\in {\mathbb N}$,
$$J_n{\mathcal O}_{V^{\prime}}={\mathcal I}(H)^nJ_n^{\prime}$$
for some  sheaf of ideals $J_n^{\prime}\subset {\mathcal
O}_{V^{\prime}}$.  Then the {\em weighted transform of ${\mathcal
G}$} is defined as
$${\mathcal G}^{\prime}:=\oplus_{n}J_n^{\prime}W^n.$$
}
\end{Paragraph}

The next proposition gives a local description of the weak
transform of a Rees algebra ${\mathcal G}$ after a permissible
monoidal transformation.
\begin{Proposition}\label{localt}\cite[Proposition 1.6]{EV} Let ${\mathcal G}=\oplus_nJ_nW^n$ be a
Rees algebra on a smooth scheme $V$   over a field $k$, and let
$V\leftarrow V^{\prime}$ be a permissible transformation. If
${\mathcal G}$ is generated by
$\{g_{n_1}W^{n_1},\ldots,g_{n_s}W^{n_s}\}$ then ${\mathcal
G}^{\prime}$ is generated by
$\{g_{n_1}^{\prime}W^{n_1},\ldots,g_{n_s}^{\prime}W^{n_s}\}$,
where   $g_{n_i}^{\prime}$ denotes the weighted transform of $g_{n_i}$
for  $i=1,\ldots,s$.
\end{Proposition}

\begin{Paragraph}\label{propertiest} {\rm {\bf Integral closure,
 differential operators and  weighted  transforms. \cite[4.1]{EV}} {\bf Giraud's Lemma.} Let ${\mathcal
G}_1 \subset {\mathcal G}_2\subset {\mathcal G}_3$ be an inclusion
of Rees algebras, such that ${\mathcal G}_3$ is the differential
algebra spanned by ${\mathcal G}_1$, and let $V\leftarrow
V^{\prime}$ be a permissible monoidal transformation with center
$Y\subset \mbox{Sing } {\mathcal G}_1$. Then:
\begin{enumerate}
\item[(i)] There is an inclusion of weighted transforms
$${\mathcal G}_1^{\prime} \subset {\mathcal G}_2^{\prime}\subset
{\mathcal G}_3^{\prime}.$$
\item[(ii)] The three algebras ${\mathcal G}_1^{\prime} \subset {\mathcal G}_2^{\prime}\subset
{\mathcal G}_3^{\prime}$ span the same differential Rees algebra.
\item[(iii)] If ${\mathcal G}_1 \subset
{\mathcal G}_2$ is a finite extension,  then ${\mathcal
G}_1^{\prime} \subset {\mathcal G}_2^{\prime}$ is a finite extension as well.
\end{enumerate}
}
\end{Paragraph}

\

\noindent{\bf A notion of equivalence for  Rees algebras}

\begin{Paragraph}\label{specialsmooth}
{\rm  If ${\mathcal G}^{}=\bigoplus I_kW^k$ is a differential
${\mathcal O}_V$-algebra and $V^{\prime\prime}\to V$ is a smooth
morphism, then the natural extension ${\mathcal
G}^{\prime\prime}=\bigoplus I_k{\mathcal O}_{V^{\prime\prime}}W^k$
is also a differential Rees algebra (cf. \cite[Proposition
5.1]{integraldifferential}). Moreover, if $\phi: T\to V$ is a
morphism of smooth schemes then $\phi^{*}({\mathcal G})$ is a
differential Rees algebra on $T$ and $\mbox{Sing }\phi^{*}({\mathcal
G})=\phi^{-1}(\mbox{Sing }{\mathcal G})$ (cf. \cite[Theorem
5.4]{integraldifferential}).

\

There are two types of smooth morphisms that we are specially
interested in:
\begin{itemize}
\item[i.] If $U\subset V$ is an open
subset (in Zariski's or \'etale topology), then the restriction of ${\mathcal G}$ to $U$ is a Rees
algebra, and if ${\mathcal G}$ is a differential Rees algebra, so is
its restriction.
\item[ii.] If $\phi: T=V\times {\mathbb A}^k\to V$ is the projection, then the pull back $\phi^*{\mathcal G}$
 is a Rees algebra. Moreover, if ${\mathcal G}$
is a differential Rees algebra, then so is $\phi^*{\mathcal G}$.
\end{itemize}
}

\end{Paragraph}

\begin{Definition} \label{permissibletransformation} {\rm Let
${\mathcal G}$ be a Rees algebra. A morphism  $V^{\prime}\to V$ is 
{\em permissible } if it is either a permissible
monoidal transformation as in Definition \ref{monoidal}, or a smooth
morphism as described in \ref {specialsmooth} (i) or (ii).}
\end{Definition}

We shall consider a smooth scheme $V$ together with a set $E$ of
smooth hypersurfaces having normal crossings, so we  present our
data as  
\begin{equation}
(V,\mathcal{G},E), 
\end{equation}
 which we call a {\em basic object} 
paralleling the notation used with pairs. 

\begin{Definition}
{\rm  A {\em local sequence of     basic objects}
takes the following form:
\begin{equation}\label{ABCranso}
(V,\mathcal{G},E) \longleftarrow
(V'_1,\mathcal{G}_1,E_1)\longleftarrow \cdots \longleftarrow
(V'_s,\mathcal{G}_s,E_s),
\end{equation}
where $(V,\mathcal{G},E) \longleftarrow (V'_1,\mathcal{G}_1,E_1)$, and each
$(V'_i,\mathcal{G}_i,E_i) \longleftarrow (V'_{i+1},\mathcal{G}_{i+1},E_{i+1})$,
 is a pull-back, or a pull-back followed by a permissible monoidal transformation
 defined  with a center $Y_i \subset \mbox{Sing }(\mathcal{G}_i)$ having normal
 crossings with the hypersurfaces in $E_i$. In this last case $E_{i+1}$ consists
 of the strict transforms of the hypersurfaces in $E_i$ together with the exceptional
 hypersurface introduced by the monoidal transformation.}
\end{Definition}

\begin{Definition} \label{Defweak} {\rm Two Rees algebras $\mathcal{G}_i$ , $i=1,2$, or two basic
objects $(V,  \mathcal{G}_i,E)$ , $i=1,2$, are said to be {\em
weakly equivalent} if: $\mbox{Sing }(\mathcal{G}_1)=\mbox{Sing
}(\mathcal{G}_2),$ and if any local sequence  of one of them, say,
\begin{equation*}
(V',\mathcal{G}_i,E') \longleftarrow
(V'_1,\mathcal{G}_{i,1},E'_1)\longleftarrow \cdots \longleftarrow
(V'_s,\mathcal{G}_{i,s},E_s),
\end{equation*}
defines a local sequence of transformation of the other, and
$\mbox{Sing }(\mathcal{G}_{1,j})=\mbox{Sing }(\mathcal{G}_{2,j})$
for $0\leq j \leq s$.}
\end{Definition}

The following Theorem is derived from the cited result of Hironaka and 
Giraud's Lemma (see \ref{propertiest}). This fact, and many applications of it, are studied in \cite{mariluz}.

\begin{Theorem} \label{examplesweak} \cite[p. 119]{Hironaka05}, \cite{Hironaka03}  
The following hold for Rees algebras defined on a smooth scheme $V$ over a perfect field $k$: 
\begin{enumerate}
\item ${\mathcal G}$ and ${\mathbb D}\mbox{iff}({\mathcal G})$ are weakly equivalent. 
\item If $\overline{{\mathcal G}}_1=\overline{{\mathcal G}}_2$ 
then  ${\mathcal G}_1$ and
${\mathcal G}_2$ are weakly equivalent. 
\item Local-Global Principle. Two Rees algebras ${\mathcal G}_1$ and
${\mathcal G}_2$ are weakly equivalent if and only if $\overline{{\mathbb D}\mbox{iff}({\mathcal G}_1)}=\overline{{\mathbb D}\mbox{iff}({\mathcal G}_2)}$. 
\end{enumerate}
\end{Theorem}

\begin{Theorem} \label{weakorder}  \cite[p. 101]{Hironaka05}, \cite{Hironaka03}
If ${\mathcal G}_1$ and ${\mathcal G}_2$ are weakly equivalent, then
 $\mbox{ord}_x{\mathcal G}_1=\mbox{ord}_x{\mathcal
G}_2$ for each  $x\in \mbox{Sing }{\mathcal G}_1=\mbox{Sing
}{\mathcal G}_2$.

\end{Theorem}

\section{Simple points and tangent cones}\label{tangentcones}
Let ${\mathcal G}=\bigoplus_{n\geq 0}I_nW^n$ be a Rees algebra on
a $d$-dimensional smooth scheme $V$ over a field $k$.  We present
here the notion of {\em  $\tau$-invariant} at a singular point $x\in
\mbox{Sing }{\mathcal G}$,  say $\tau_{{\mathcal G},x}$. This is defined as the codimension 
of a subspace in the tangent space of a point $x\in \mbox{Sing }{\mathcal G}$.   
It provides
local information on the singularity since it   is a bound on the
local codimension of the singular locus (see Theorem
\ref{teodlcod}). From the point of view of resolution, it indicates the number of variables can  be ``eliminated", via elimination algebras, as we shall see in
a coming section (see \ref{elagatai}).

\begin{Definition} \label{simplepoint} {\rm A point $x\in
\mbox{Sing }{\mathcal G}$ is {\em simple} if for some index  $k\geq 1$
the order of $I_k$ in $x$, $\nu_x (I_k)$, is $k$, (i.e., if
$\mbox{ord}_x {\mathcal G}=1$).}
\end{Definition}

\begin{Paragraph} \label{tangentcone} {\bf The tangent cone.}  \cite[4.2]{hpositive}  {\rm Let
$x\in \mbox{Sing }{\mathcal G}$ be a closed point. Consider
the graded algebra associated to the closed point's  maximal ideal $m_x$,
$\mbox{Gr}_{m_x}({\mathcal O}_{V,x})$ (which is isomorphic to a
polynomial ring in $d$-variables, say $k^{\prime}[Z_1,\ldots,Z_d]$). This is the coordinate ring
associated to the tangent space of $V$ at $x$, namely
$\mbox{Spec}(\mbox{Gr}_{m_x}({\mathcal O}_{V,x}))={\mathbb
T}_{V,x}$. 

\

The {\em initial ideal} or {\em tangent ideal} of
${\mathcal G}$ at $x$, ${\mbox{In}}_{x}({\mathcal G})$, is the
ideal of $\mbox{Gr}_{m_x}({\mathcal O}_{V,x})$ generated by the
elements $\mbox{In}_x(I_n)$ for all $n\geq 1$.  The zero set of the tangent ideal in $\mbox{Spec
}(\mbox{Gr}_{m_x}({\mathcal O}_{V,x}))$ is the {\em tangent cone} of
${\mathcal G}$ at $x$, ${\mathcal C}_{{\mathcal G},x}$.

\

>From the algebraic point of view,   $\tau_{{\mathcal G},x}$
is defined as  the minimum number of variables needed to describe
${\mbox{In}}_{x}{\mathcal G}$. From the geometric point of view,
$\tau_{{\mathcal G},x}$ is the codimension of the largest linear
subspace  ${\mathcal L}_{{\mathcal G},x}\subset {\mathcal
C}_{{\mathcal G},x}$ such that $u+v\in {\mathcal C}_{{\mathcal
G},x}$ for all $u\in {\mathcal C}_{{\mathcal G},x}$ and all $v\in
{\mathcal L}_{{\mathcal G},x}$.

\

Some facts about tangent ideals:\begin{enumerate}\item[(i)] The tangent ideal
${\mbox{In}}_{x}{\mathcal G}$ is zero unless $x\in \mbox{Sing
}{\mathcal G}$ is a simple point.
\item[(ii)] If  ${\mathcal G}$ is a differential
algebra and if  $k^{\prime}$ (the residue field at $x$) is a field of
characteristic zero, then ${\mbox{In}}_{x}{\mathcal G}$ is
generated by linear forms. If $k^{\prime}$ is a field of positive
characteristic $p$, then there is a sequence  $ e_0 < e_1 < \cdots
< e_r$ in $\mathbb{Z}_{\geq 0}$, such that ${\mbox{In}}_{x}{\mathcal G}$ is
generated by elements of the form
\begin{equation}\label{lis}
l_1,\ldots,l_{s_0},l_{s_0+1},\ldots,l_{s_1},\ldots,l_{s_r-1},\ldots,l_{s_r}. 
\end{equation}
Here    $l_1,\ldots,l_{s_0}$ are linear  combinations of powers
$Z_i^{p^{e_0}}$, and  if $t\geq 0$,
$$l_{s_t+1},\ldots,l_{s_{t+1}}$$ are linear combinations of powers
$Z_i^{p^{e_t}}$. Furthermore the  $s_r$ homogeneous elements in (\ref{lis}) form a regular sequence at
$\mbox{Gr}_{m_x}({\mathcal O}_{V,x})$.

So
$\langle l_1,\ldots,l_{s_r}\rangle$ define a subscheme of codimension
$s_r$ in ${\mathbb
T}_{V,x}$. If $k'$ is a perfect field the radical of this ideal is spanned by linear forms, defining a subspace of codimension $s_r$ in ${\mathbb
T}_{V,x}$.

The integer $s_r$ is  the $\tau$-invariant of the
singularity at $x$ which we have denoted   by $\tau_{{\mathcal G},x}$. If
$p^{e_0}$ is the smallest power of $p$ in (\ref{lis}) then the order
of $I_n$ in ${\mathcal O}_{V,x}$ is $n$ if and only if $n$ is a
multiple of $p^{e_0}$.

\item[(iii)] If ${\mathcal G}$ is a differential Rees algebra then: $${\mathcal L}_{{\mathcal G},x}= {\mathcal
C}_{{\mathcal G},x}.$$
\item[(iv)] For any Rees algebra ${\mathcal G}$, the inclusion ${\mathcal G}\subset {\mathbb D}\mbox{iff}({\mathcal G})$ defines an inclusion ${\mathcal
C}_{{\mathbb D}\mbox{iff}({\mathcal G}),x} \subset {\mathcal
C}_{{\mathcal G},x}$, and:
$$({\mathcal
L}_{{\mathbb D}\mbox{iff}({\mathcal G}),x}=){\mathcal
C}_{{\mathbb D}\mbox{iff}({\mathcal G}),x}={\mathcal L}_{{\mathcal
G},x}.$$
\item[(v)]
If $Y\subset \mbox{Sing }{\mathcal G}$ is a permissible center, then
$ {\mathbb T}_{Y,x}\subset  {\mathbb T}_{V,x}$, is a linear
subspace, and forthermore ${\mathbb T}_{Y,x}\subset {\mathcal
L}_{{\mathcal G},x}$ for all $x\in Y\subset \mbox{Sing }{\mathcal
G}$. In particular $\tau_{{\mathcal G},x}$ bounds the local
codimension of the regular scheme $Y$ in $V$, i.e., $\mbox{co-dim}_x
Y\geq \tau_{{\mathcal G},x}$.

\end{enumerate}}

\end{Paragraph}

The following Theorem is due to Hironaka:
\begin{Theorem} \label{tauweak}  If ${\mathcal G}_1$ and ${\mathcal
G}_2$ are weakly equivalent, then for each $x\in \mbox{Sing
}{\mathcal G}_1=\mbox{Sing }{\mathcal G}_2$ there is an equality
between their $\tau$-invariants, i.e., $\tau_{{\mathcal
G}_1,x}=\tau_{{\mathcal G}_2,x}$.
\end{Theorem}

\begin{Definition}\label{cdtype}
{\rm A Rees algebra ${\mathcal G}$ is said to be of {\em
codimensional type $\geq e$} if $\tau_{{\mathcal G},x}\geq e$ for
all $x\in \mbox{Sing }{\mathcal G}$.}
\end{Definition}



\begin{Theorem}\label{teodlcod}
Let $x\in \mbox{Sing }{\mathcal G}$. Then:
\begin{enumerate}
\item $\mbox{Co-dim}_x\mbox{Sing }{\mathcal G}\geq \tau_{{\mathcal
G},x}$;
\item If  the equality $\mbox{co-dim}_x\mbox{Sing }{\mathcal G}=\tau_{{\mathcal
G},x}$  holds, then $\mbox{Sing }{\mathcal
G}$ is smooth locally at $x$.
\end{enumerate}

\end{Theorem}
We shall prove that locally at $x$, $\mbox{Sing }{\mathcal G}$ is
included in a complete intersection  scheme of codimension
$\tau_{{\mathcal G},x}$ (see Corollary \ref{PA}), which proves the
first assertion. The second claim  will be addressed in Lemma
\ref{codimensionone} and Remark \ref{rkcdime}. 

\

Let ${\mathcal G}$ be a Rees algebra on a $d$-dimensional 
smooth scheme $V$. If $\mbox{co-dim}_x(\mbox{Sing }{\mathcal G})>  \tau_{{\mathcal
G},x}$, then we will associate to ${\mathcal G}$  an
 {\em elimination algebra} in dimension $d-\tau_{{\mathcal
G},x}$, locally, in a neighborhood of $x$. Hence we will  eliminating $d-\tau_{{\mathcal
G},x}$ variables in our original problem posed in a $d$-dimensional space.

\part{Elimination}\label{partElimination}

\section{Elimination via universal invariants}\label{evui}
As indicated in the Introduction and specially in Part
\ref{Char0},  the resolution of singularities of a hypersurface over a field of
characteristic zero can be achieved in two Steps A and B.
In Step A, a suitable stratification of the locus of  maximum multiplicity
 is constructed using an inductive argument;  we
briefly explained how the notion of maximal contact
plays a  role in this stratification.

\

We will introduce a new approach, using projections and  {\em universal elimination
algebras}  in our inductive arguments. When the characteristic of the base field is zero, projections are natural substitutes of restrictions to smooth 
hypersurfaces of maximal contact. In this case, our approach leads 
to the same resolution invariants already obtained in constructive resolution 
of singularities,   where maximal contact was used. 
 
\

This section is organized as follows: in \ref{Motivation} we
discuss the motivation  for  using elimination (see Example
\ref{Example2} for more details);  universal elimination algebras
are defined in \ref{par61}, and their relation to differential
operators is described in \ref{pdiffop}. Finally, in Theorem
\ref{elimRamMult} we explain how these universal invariants
specialize to provide information on Rees algebras.

\begin{Paragraph} \label{Motivation} {\bf The motivation.} {\rm Assume that $S$ is a regular
ring containing a field $k$. Let $$f(Z)=Z^n+a_1Z^{n-1}+\ldots+a_n\in
S[Z]$$ and denote by $\Upsilon_n$ the set of points in $\{f(Z)=0\}$
with multiplicity $n$. The natural inclusion $S\subset S[Z]$ induces
a smooth morphism $\beta$, and a  finite restriction
$\overline{\beta}$,
\begin{equation}\label{casoesp}\begin{array}{rcc} \mbox{Spec}
\left(S[Z]/\langle f(Z)\rangle\right) &
 \hookrightarrow & \mbox{Spec} (S[Z])\\
 \overline{\beta}\searrow & & \beta \downarrow \\
  & & \mbox{Spec} (S).
  \end{array}
  \end{equation}
Our goal is to find functions on  the coefficients of $f(Z)$ that
describe the image in $\mbox{Spec} (S)$ of $\Upsilon_n$. The {\em
elimination algebra of the polynomial $f(Z)$} will be  the $k$-subalgebra of $S$
generated by these functions.

\

Notice that  $B= S[Z]/\langle f(Z)\rangle$ is a free $S$-module of
rank $n$.  Let $Q$ be a prime ideal in $B$ dominating $S$ at a
prime $P$. Under these conditions, Zariski's projection formula
for multiplicities   ensures that the multiplicity of $B_Q$ is at
most $n$. Moreover, if   this multiplicity is exactly  $n$, then $Q$ is the unique
prime in $B$ daminating  $P$, and $B_Q$ and $A_P$ have the
same residue fields (see \cite[Corollary 1, p. 299]{ZS}). The
morphism $\overline{\beta}$ is said to be purely ramified over a
point $x\in  \mbox{Spec} (S)$ if the {\em geometric fiber} over
$x$ is a unique point.

\begin{quote}
{\em So the multiplicity formula says that $\Upsilon_n$ is
contained in the set of points where $\overline{\beta}$ is purely
ramified.}
\end{quote}

\begin{Example} \label{Example2} {\rm Suppose  $n=2$,  
 let $f(Z)=Z^2+a_1Z+a_2$ with $a_1, a_2\in S$, and set $\overline{\beta}$ as above. In
this case    the discriminant-namely $a_1^2-4a_2\in S$-describes
the image under $\overline{\beta}$ of the purely ramified locus in
$\mbox{Spec} (S)$. Notice that $a_1^2-4a_2\in S$ is a weighted
homogeneous polynomial of degree two,  provided that we assign weight
one to $a_1$ and weight two to  $a_2$. It is not hard to check
that if the characteristic of $S$ is not 2, then the closed subset  in  $\mbox{Spec}(S)$ where the discriminant has
order at least two, is exactly the image via $\overline{\beta}$ of the two
 fold points.

\

Using the language of Rees algebras, our datum is the $S[Z]$-algebra generated by $f(Z)$ in degree two, say ${\mathcal
G}=S[Z][fW^2] (\subset S[Z][W])$. Therefore its singular locus is
the set of $2$-fold points of $\{f(Z)=0\}$. As we shall see, in
this case the elimination algebra associated to ${\mathcal G}$
is the Rees algebra over $S$ generated by $a_1^2-4a_2$ in
degree 2, say ${\mathcal R}_{\mathcal G}=S[(a_1^2-4a_2)W^2]$. If the the characteristic 
is different from 2, then $\mbox{Sing } {\mathcal R}_{\mathcal G}$   is the image of the two-fold points of
$\{f(Z)=0\}$.

 \

However, this argument  breaks down  if the characteristic is two, a situation that
requires some attention.  This problem can be remedied  by
extending ${\mathcal G}$ to a differential Rees algebra. This is done 
by considering a Rees algebra with more generators by 
 applying all differential
operators to $f(Z)$. This already forces us to extend the notion of
elimination algebra to the case of several polynomials, since
typically a differential Rees algebra will have more than one
generator. When a Rees algebra is
differential, then its singular locus can be identified with the
singular locus of its elimination algebra. These ideas will be explored further, specially in  Section
\ref{localprojection} (see \ref{eliminationalg}). }
\end{Example}
}
\end{Paragraph}


\begin{Paragraph} \label{par61} 
{\bf  The universal elimination algebra.} {\rm  Let $k$ be a
field. Consider the polynomial ring in $n$ variables
$k[Y_1,\ldots,Y_n]$, and  the {\em universal polynomial} of degree $n$,
$$F_n(Z)=(Z-Y_1)\cdots(Z-Y_n)=Z^n-s_{n,1}Z^{n-1}+\ldots+(-1)^ns_{n,n}\in
k[Y_1,\ldots,Y_n,Z],$$ where for $i=1,\ldots,n$, $s_{n,i}\in
k[Y_1,\ldots,Y_n,Z]$ denotes the $i$-th symmetric polynomial in
$n$ variables.

\

Observe that the diagram
\begin{equation}\label{universalespecial}\begin{array}{rcc} \mbox{Spec
}\left(k[s_{n,1}, \ldots,s_{n,n}][Z]/\langle F_n(Z)\rangle \right)
 &  \hookrightarrow &  \mbox{Spec} (k[s_{n,1}, \ldots,s_{n,n}][Z])  \\
\overline{\alpha}\searrow &   & \alpha\downarrow   \\ & &
\mbox{Spec} (k[s_{n,1}, \ldots,s_{n,n}])
\end{array}
\end{equation}
specialices to 
(\ref{casoesp}) via
\begin{equation}
\label{especial}
\begin{array}{rrcl}
\Theta: & k[s_{n,1},\ldots,s_{n,n}] & \longrightarrow & S \\
 & (-1)^is_{n,i} & \to & a_i.
 \end{array}
 \end{equation}

In the following lines we  consider the  universal case in 
(\ref{universalespecial}). Our motivation 
is to find equations on the coefficients of the polynomial $F_n$
that describe the image of the $n$-fold points of $F_n=0$. We 
begin  by looking for functions on the coefficients that
describe the purely ramified locus of the morphism; we
reproduce arguments from \cite[Section 1]{hpositive}.

\

First notice that the group of permutations of $n$ elements, $S_n$,
acts linearly on $k[Y_1,\ldots,Y_n]$ and that the subring of
invariants is
$$k[Y_1,\ldots,Y_n]^{S_n}=k[s_{n,1},\ldots,s_{n,n}].$$
 Set $\mathbb{T}= k[s_{n,1},\ldots,s_{n,n}]$ and observe that
$\mathbb{T}\subset k[Y_1,\ldots,Y_n]$ is an inclusion of graded
rings  since the action of $S_n$ in $k[Y_1,\ldots,Y_n]$ is linear
(i.e., it preserves the grading provided that we assign  weight $i$ 
 to $s_{n,i}$
).

\

In the setting of (\ref{casoesp}) the purely ramified locus 
is independent on the choice of $Z$, at least 
if we stick to changes of the form 
\begin{equation} \label{changes} uZ+\alpha,
\end{equation}
 with $\alpha, u\in S$ and $u$
 invertible. So, in finding equations in $S$ describing
$\overline{\beta}(\Upsilon_n)$  we have to look for functions  on
the coefficients of $f$ that are invariant under changes as in
(\ref{changes}).

\

We consider first changes of the form $Z+\alpha$. In the universal case, these changes of variable
 can be expressed as
$$F_n(Z+T)=(Z-(Y_1-T))\cdots(Z-(Y_n-T)) \in
k[Y_1-T,\ldots,Y_n-T]^{S_n}[Z].$$

\

The group $S_n$ also acts linearly on $k[Y_i-Y_j]_{1\leq
i,j\leq n}$ defining a graded subring
$$k[Y_i-Y_j]_{1\leq i,j\leq n}^{S_n}\subset
\mathbb{T}=k[Y_1,\ldots,Y_n]^{S_n}=k[s_{n,1},\ldots,s_{n,n}].$$
The elements on the left hand side   
 are functions on the coefficients of the universal polynomial  which are clearly invariant by any change of the form $Z+T$. Through $\Theta$ (see (\ref{especial})) they define functions on the coeficients of $f(Z)$ which are invariant under changes of the form $Z+\alpha$.

\

Let ${\mathcal U}$ be the  $k$-subalgebra of $\mathbb{T}[Z]$
generated by $F_n(Z)$ in degree $n$, say $k[F_n(Z)]$. Note that $F_n(Z)$ is weighted 
homogeneous in $k[s_{n,1},\ldots,s_{n,n}][Z]$. Define the 
{\em universal elimination algebra of 
${\mathcal U}$}, as  
\begin{equation} \label{universalgenerators}
{\mathcal R}_{\mathcal U}:=k[Y_i-Y_j]_{1\leq i,j\leq
n}^{S_n}.
\end{equation}
${\mathcal R}_{\mathcal U}$ can be generated by homogeneous elements 
$H_{m_1},\ldots,H_{m_r}$ of degrees $m_1,\ldots,m_r$; i.e.,
$${\mathcal R}_{\mathcal U}=[H_{m_1},\ldots,H_{m_r}]. $$
Here each $H_{m_i}$ is a homogeneous
polynomial in degree $m_i$ in $Y_1,\ldots,Y_n$, 
 and it is also a weighted homogeneous
polynomial in $s_{n,1},\ldots,s_{n,n}$.  For instance, in Example \ref{Example2} the
elimination algebra is generated by the discriminant in degree
two.

\

Express  ${\mathcal R}_{\mathcal U}=\oplus_NI_NW^N$.   For each positive index $N$, the homogeneous polynomials in $I_N$  form a finite dimensional vector space over $k$. Polynomials in this vector space are weighted homogeneous on the coefficients  $s_{n,1},\ldots,s_{n,n}$. 

\
Via  (\ref{especial}) we can define two Rees algebras, one over $S[Z]$ and 
another over $S$. For any specific monic polynomial
of degree $n$, say $f(Z)=Z^n+a_1Z^{n-1}+\ldots+a_n\in
S[Z]$ over a $k$-algebra $S$, via (\ref{especial}) the ideals $I_N$ lead to ideals     $J_N$  spanned by  weighted homogeneous functions on the coefficients $a_i$. These ideals $J_N (\subset S)$  will be invariant under any change of variables as in (\ref{changes}). Note also that for any two positive integers $N$, $M$:  $J_N\cdot J_M\subset J_{N+M}$. Now  $\oplus_NJ_NW^N\subset S[W]$ 
is a Rees algebra free of the variable $Z$ which we refer to as  an {\em elimination algebra of}  
$S[Z][f(Z)W^n]$.

}
\end{Paragraph}


\

\noindent {\bf Differential operators}

\

\begin{Paragraph} \label{pdiffop} {\rm
Now our purpose is  to get a better understanding of the
information encoded in the universal elimination algebra. Given a
diagram as in (\ref{universalespecial}), we want  to study how the
universal elimination algebra can be used to describe the image
under $\alpha$ of the purely ramified locus, and in turn, of the
set of $n$-fold points of $\{F_n=0\}\subset \mbox{Spec
}(K[s_{n,1},\ldots,s_{n,n}][Z])$.

\

 To understand how to reach this goal, we cite the following  Lemma from
\cite{hpositive},  which  relates the multiple roots of a polynomial
to the vanishing of its derivatives:
\begin{Lemma}\label{lemakey}\cite[Lemma 1.3]{hpositive} Let $K$ be an algebraically closed field and
let $f(Z)\in K[Z]$ be a polynomial of degree $n$. Then the
following are equivalent:
\begin{enumerate}
\item[i.] $\Delta^kf(z)$ is nilpotent in $K[Z]/\langle
f(z)\rangle$ for $0\leq k< n$.
\item[ii.] $f(Z)=(Z-\alpha)^n$ for some $\alpha\in K$.
\end{enumerate}
\end{Lemma}

\

We start by introducing  differential operators  in the
universal case:  let $T,Z$ be variables, and let
$k[Y_1,\ldots,Y_n,Z,T]$ be the polynomial ring in $n+2$ variables.
Consider the $k[Y_1\ldots,Y_n]$-morphism:
$$\begin{array}{rrcl}
Tay: &  k[Y_1,\ldots,Y_n][Z] & \longrightarrow &
k[Y_1,\ldots,Y_n][Z,T]\\
 & Z & \to & Z+T.
 \end{array}$$
For each polynomial $G(Z)\in k[Y_1,\ldots,Y_n][Z]$,
$$Tay(G(Z))=\sum_{k\geq 0}G_k(Z)T^k,$$
and for each index $k$, we can define the operators:
\begin{equation}\label{p26}
\begin{array}{rrcl}
\Delta^k: &  k[Y_1,\ldots,Y_n][Z] & \longrightarrow &
k[Y_1,\ldots,Y_n][Z]\\
 & G(Z) & \to & \Delta^k(G(Z)):=G_k(Z).
 \end{array}
 \end{equation}
For $k\geq 0$ the $\Delta^k$ are particular  {\em differential
operators of degree $k$, relative to the inclusion}
$k[Y_1,\ldots,Y_n]\subset k[Y_1,\ldots,Y_n][Z]$. }
\end{Paragraph}

Now consider the universal monic polynomial of degree $n$:
$$F_n(Z)=(Z-Y_1)\cdots(Z-Y_n)=Z^n-s_{n,1}Z^{n-1}+\ldots+(-1)^ns_{n,n}\in
k[Y_1,\ldots,Y_n,Z].$$

Observe that
$$Tay(F_n(Z))=F_n(Z+T)=(Z+T-Y_1)\cdots(Z+T-Y_n)=(T-(-Z+Y_1))\cdots(T-(-Z+Y_n)),$$
and that the coefficients of this polynomial in the variable $T$
are precisely the symmetric polynomials in the variables
$Z-Y_1,\ldots,Z-Y_n$, i.e.,
\begin{equation}
\label{differentialsim}
\Delta^k(F_n(Z))=(-1)^{n-k}s_{n,n-k}(-Z+Y_1,\ldots,-Z+Y_n),
\end{equation}
for $k=1,\ldots,n-1$. In this setting,  the action of $S_n$ in
$k[Y_1,\ldots,Y_,Z]$ can be considered as a permutation of
$Y_1,\ldots,Y_n$ that fixes $Z$. Hence,
\begin{equation}
\label{differentialrelative}
k[Z-Y_1,\ldots,Z-Y_n]^{S_n}=k[F_n(Z),\{\Delta^k(F_n(Z))\}_{k=1,\ldots,n-1}].
\end{equation}
Let us stress here that  $\Delta^k(F_n(Z))$ is homogeneous of degree
$n-k$ for $k=1,\ldots,n-1$, so $k[F_n(Z),\{\Delta^k(F_n(Z))\}_{k=1,\ldots,n-1}]$ is a graded
subring in $k[Y_1,\ldots,Y_,Z]$.

\

Now, since $Y_i-Y_j=(Z-Y_j)-(Z-Y_i)$, we have that
$$k[Y_i-Y_j]_{1\leq i,j \leq n}\subset k[Z-Y_1,\ldots,Z-Y_n].$$
Hence there is an inclusion of graded algebras
\begin{equation}\label{eqr}
\begin{array}{c}
k[H_{n_1},\ldots,H_{n_r}]=k[Y_i-Y_j]_{1\leq i,j \leq
n}^{S_n}\subset\\
\subset
k[Z-Y_1,\ldots,Z-Y_n]^{S_n}=k[F_n(Z),\{\Delta^k(F_n(Z))\}_{k=1,\ldots,n-1}],
\end{array}
\end{equation}
and therefore for $i=1,\ldots,r$  each $H_{n_i}$ is also weighted
homogeneous in the elements 
\begin{equation}
\label{differentialrelativebis}
\{F_n(Z),\{\Delta^k(F_n(Z))\}_{k=1,\ldots,n-1}\}.
\end{equation}

\

To conclude, it can be shown that $k[s_{n,1},
\ldots,s_{n,n}][Z]/\langle F_n(Z)\rangle\simeq k[s_{n,1},
\ldots,s_{n,n}][Y_1]$, and setting $Z=Y_1$ in $k[F_n(Z),\{\Delta^k(F_n(Z))\}_{k=1,\ldots,n-1}]$:
\begin{equation}\label{eqrtay}
\begin{array}{c}
k[H_{n_1},\ldots,H_{n_r}]=k[Y_i-Y_j]_{1\leq i,j \leq
n}^{S_n}\subset\\
\subset
k[F_n(Y_1),\{\Delta^k(F_n(Y_1))\}_{k=1,\ldots,n-1}](\subset
k[s_{n,1}, \ldots,s_{n,n}][Y_1]),
\end{array}
\end{equation}
is a finite extension of graded rings. This result,  in combination with Lemma
\ref{lemakey}, gives the following theorem.

\begin{Theorem}\cite[Theorem 1.16]{hpositive}
\label{elimRamMult} Let $S$ be a $k$-algebra, let
$f(Z)=Z^n+a_1Z^{n-1}+\ldots+a_{n-1}Z+a_n\in S[Z]$ and consider a
commutative diagram
\begin{equation}\label{casoesp1}\begin{array}{rcc} \mbox{Spec}
\left(S[Z]/\langle f(Z)\rangle\right) &
 \hookrightarrow & \mbox{Spec} (S[Z])\\
 \overline{\beta}\searrow & & \beta \downarrow \\
  & & \mbox{Spec} (S).
  \end{array}
  \end{equation}
as described  in (\ref{casoesp}). Let ${\mathcal G}$ be the Rees algebra
generated  by $f(Z)$ in degree $n$, say ${\mathcal G}=S[Z][f(Z)W^n]$, and let $\Upsilon_n$ denote
the set of $n$-fold points of $\{f(Z)=0\}\subset
\mbox{Spec}(S[Z])$, i.e.,
 $\Upsilon_n=\mbox{Sing }{\mathcal G}$. Consider the specialization morphism,
$$\begin{array}{rcl}
\mathbb{T}=k[s_{n,1},\ldots,s_{n,n}] & \longrightarrow & S\\
s_{n,i} & \to & (-1)^ia_i
\end{array}$$
which gives rise to the elimination algebra associated to
${\mathcal G}$,
$${\mathcal R}_{\mathcal G}=S[H_{m_j}(a_1,\ldots,a_n)W^{m_j}, j=1,\ldots,r]\subset S[W],$$
where $m_j$ denotes the degree of the weighted homogeneous polynomial $H_{m_j}(s_{n,1},\ldots,s_{n,n})$ as mentioned above.

Then:
\begin{enumerate}
\item[i)] The closed set
$V(H_{m_j}(a_1,\ldots,a_n);j=1,\ldots,r)\subset \mbox{Spec} (S)$
is the image of the set of points where $\overline{\beta}$ is
purely ramified (see \ref{Motivation}).
\item[ii)] If $S$ is regular, then
\begin{equation}\label{inclusionSing}
\beta(\Upsilon_n)=\beta(\mbox{Sing }{\mathcal G})\subset
\mbox{Sing }{\mathcal R}_{\mathcal G}. \end{equation} If in
addition, the characteristic of $S$ is zero, then the inclusion
in (\ref{inclusionSing}) is an equality.
\end{enumerate}

\end{Theorem}

\begin{Paragraph}\label{comeafter} {\bf Elimination algebras in the
general case.} {\rm So far, we have defined elimination algebras 
for Rees algebras generated by one element, say ${\mathcal G}=S[fW^n]$. 
 Elimination algebras can also be defined for
Rees algebras with more than one generator. 
A case of  particular interest  is that of the differential Rees algebra generated by ${\mathcal G}$,
namely ${\mathbb D}\mbox{iff}({\mathcal G})$  (see \ref{extensiondif}). 
In this setting   Theorem \ref{elimRamMult} can be qualitatively improved,
since:
\begin{equation}\label{inclusionDiffSing}
\beta(\Upsilon_n)=\beta(\mbox{Sing }{\mathbb D}\mbox{iff}(
{\mathcal G}))= \mbox{Sing }{\mathcal R}_{{\mathbb D}\mbox{iff}(
{\mathcal G})}
\end{equation}
in any characteristic (see \cite[Corollary 4.12]{hpositive}, or
\ref{eliminationalg} below). We refer to \cite[1.23-1.40]{hpositive}
for more  details on the construction of the elimination universal
algebra associated to a Rees algebra with more than one generator. Also,  
we refer to
\ref{explicitdescription} and  \ref{eliminationalg}
where we indicate  how elimination algebras can be computed. In
Example \ref{canguro}  we provide a concrete example.}
\end{Paragraph}

\section{A smooth  local projection and the elimination algebra} \label{localprojection}
Once the universal case has been treated in the previous section, we
are ready to study the case to be considered here: how to define
elimination in terms of algebras on smooth schemes.  Let $V=V^{(d)}$
be a $d$-dimensional smooth scheme of finite type over a perfect  field $k$.
Let ${\mathcal G}=\oplus_{n\in {\mathbb N}}I_nW^n$ be a sheaf of
Rees algebras and let $x\in \mbox{Sing }{\mathcal G}$ be a simple
point not contained in any component of codimension one  of
$\mbox{Sing }{\mathcal G}$ (see Definition \ref{simplepoint}). In
the following we describe how to construct:
\begin{itemize}
\item A suitable smooth 
local projection (in an  \'etale neighborhood of $x$), or say, a smooth morphism 
(in an \'etale neighborhood of $x$), 
$$\beta_{d,d-1}: V^{(d)}\to V^{(d-1)},$$ with $\beta_{d,d-1}(x)=x_1$. In  this
case $\beta^*_{d,d-1}:{\mathcal O}_{V^{(d-1)},{x_1}}\to {\mathcal
O}_{V^{(d)},{x}}$ is an  inclusion.
\item An elimination algebra $${\mathcal R}_{{\mathcal
G},\beta_{d,d-1}}\subset {\mathcal O}_{V^{(d-1)}}[W]$$
locally at  $x_1$.
\end{itemize}
So we start with an algebra $\mathcal G\subset
{\mathcal O}_{V^{(d)}}[W]$, and define ${\mathcal R}_{{\mathcal
G},\beta_{d,d-1}}\subset {\mathcal O}_{V^{(d-1)}}[W]$.
Although ${\mathcal R}_{{\mathcal G},\beta_{d,d-1}}$    depends
on the projection $\beta_{d,d-1}$, it will satisfy some nice
properties, and our main   invariants will  derive from them (see
\ref{eliminationalg}).

\begin{Definition}\label{def83}
{\rm Let ${\mathcal G}$ be a Rees algebra on a smooth
$d$-dimensional scheme $V^{(d)}$ over a perfect field $k$, and let $x\in
\mbox{Sing }{\mathcal G}$ be a simple closed point. We   say that a local
projection to a smooth $(d-1)$-dimensional scheme, $V^{(d-1)}$,
$$\begin{array}{rrcl}
\beta_{d,d-1}: & V^{(d)}  & \to & V^{(d-1)}\\
 & x & \to & x_1
\end{array}$$
is {\em  ${\mathcal G}$-admissible locally at $x$} if the
following conditions hold:
\begin{enumerate}
\item[(i)] The closed point $x$ is not contained in any
component of codimension one of    $\mbox{Sing }{\mathcal G}$.
\item[(ii)] The Rees algebra ${\mathcal G}$ is a $\beta_{d,d-1}$-relative differential
algebra (see Definition \ref{def82}).
\item[(iii)] {\em Transversality:} $\ker d\beta_{d,d-1}\cap   {\mathcal C}_{{\mathcal
G},x}=\{\overrightarrow{0}\}\subset {\mathbb T}_{{ V},x}$ .
\end{enumerate}}
\end{Definition}

Essentialy, the idea is that if $\beta_{d,d-1}: V^{(d)}\to V^{(d-1)}$ is 
${\mathcal G}$-admissible locally at $x$, then an 
${\mathcal O}_{V^{(d-1)}}$-Rees algebra can be assigned to 
${\mathcal G}$. If condition (i) in Definition \ref{def83} is not fullfilled, 
the assignment is not needed (see Lemma \ref{codimensionone}).  Now we will explain the role of conditions
(ii) and (iii) in Definition \ref{def83}.

\begin{Paragraph}\label{p81} {\bf Condition (ii)  in Definition
\ref{def83}: relative differential Rees algebras.} {\rm  Notice that if
${\mathcal G}=\oplus_{n\in {\mathbb N}}I_nW^n$ is an absolute
differential Rees algebra, then it is  also a relative differential
algebra for any smooth morphism $\beta_{d,d-1}: V^{(d)}\to
V^{(d-1)}$ of schemes over a field $k$ defined in a neighborhood
of $x\in \mbox{Sing }{\mathcal G}$. Thus this is a generic condition.  
Our starting point will always  be a differential Rees algebra ${\mathcal G}$. 
If a closed point $x$  is within condition (i), then it is simple to construct 
a smooth scheme $V^{(d-1)}$ and a ${\mathcal G}$-admissible morphism 
as avobe. Since ${\mathcal G}$ is an absolute differential Rees algebra, it is will 
a $\beta$-relative differential Rees algebra. This is as much as we will use to produce 
an elimination algebra over $V^{(d-1)}$.  

 A key point in our development
is the study of properties of {\em relative} differential
algebras, and their stability by monoidal transformations (see
Section \ref{eliminationmonoidal}); whereas transforms of absolute differential Rees algebras are not absolute differential (see Remark \ref{rk85} for further details). }
\end{Paragraph}

\begin{Paragraph}\label{localproj}{\rm {\bf Condition (iii) in Definition \ref{def83}: local smooth projections
and transversality. }  Almost any local smooth projection, or equivalently, almost any smooth morphism defined  locally,  in a
neigborhood of a simple point in the singular locus of a Rees
algebra,   will fulfill   condition (iii) in Definition \ref{def83}.   In \ref{opentran} we show that this condition is open: it holds at any singular point in a neighborhood of $x$. 
First  we will explain the meaning of this condition  in terms of local rings, and then we describe a procedure to construct a
smooth morphism satisfying this geometric condition at $x$. 

Suppose that a local smooth projection  to a
$(d-1)$-dimensional regular scheme is defined, say 
\begin{equation}
\label{firstprojection}
\begin{array}{rrcl} \beta_{d,d-1}:& V^{(d)}& \longrightarrow & V^{(d-1)}\\
& x & \to & x_1.
\end{array}
\end{equation}
A regular system of
parameters $\{y_1,\ldots,y_{d-1}\}\subset {\mathcal
O}_{V^{(d-1)},x_1}$  extends to parameters 
$$\{y_1,\ldots,y_{d-1},y_d\}\subset {\mathcal
O}_{V^{(d)},x}.$$ Then condition (iii) in Definition \ref{def83} holds if and only if
$\{\mbox{In}_xy_1=0,\ldots,\mbox{In}_xy_{d-1}=0\}\subset {\mathbb
T}_{V,x}$ is not contained in the tangent cone of ${\mathcal G}$
at $x$, ${\mathcal C}_{{\mathcal G},x}$.

\

This also shows how to produce local smooth projections that fulfill Condition (iii): let ${\mathcal G}=\oplus_nI_nW^n$ be a Rees algebra, and let $x\in
\mbox{Sing }{\mathcal G}$ be a simple closed  point. The
graded ideal $\mbox{In}_x\mathcal{G}$ defines the subscheme
${\mathcal C}_{{\mathcal G},x}$ of $ {\mathbb T}_{{ V},x}$ (in
fact, recall that if ${\mathcal G}$ is a differential Rees algebra, then
${\mathcal C}_{{\mathcal G},x}={\mathcal L}_{{\mathcal G},x}$, see
\ref{tangentcone},(iii)). Now select a regular system of
parameters $\{y_1,\ldots,y_{d-1},y_d\}\subset {\mathcal
O}_{V^{(d)},x}$  such that
$\{\mbox{In}_xy_1=0,\ldots,\mbox{In}_xy_{d-1}=0\}\subset {\mathbb
T}_{V,x}$ is not contained in ${\mathcal C}_{{\mathcal G},x}$.
  Note that  there is a natural injective map
from the ring of polynomials in $(d-1)$-variables with
coefficients in $k$ into ${\mathcal O}_{V^{(d)},x}$, and
localizing we get an inclusion    of regular local rings,
$$
\begin{array}{rcl}\label{local}
k[Y_1,\ldots,Y_{d-1}]_{\langle Y_1,\ldots,Y_{d-1}\rangle}&  \longrightarrow &  {\mathcal O}_{V^{(d)},x}\\
   Y_i & \to & y_i
 \end{array}
$$
This is one way to produce a local smooth projection as (\ref{firstprojection}), to a
$(d-1)$-dimensional regular scheme ($V^{(d-1)}={\mathbb A}^{d-1}_k$), satisfying condition (iii) in Definition \ref{def83}.}
\end{Paragraph}

\begin{Paragraph}\label{Zariski} {\bf Transversality and Zariski's multiplicity
formula for projections.} {\rm With the same notation as in
\ref{localproj}, fix a local smooth projection as  in
(\ref{firstprojection}).  Now our goal  is to study the image of
$\mbox{Sing }{\mathcal G}$ under the morphism $\beta_{d,d-1}$ in a
neighborhood of $x$.  We will show that if $Y$ is a smooth center in
$\mbox{Sing }{\mathcal G}$ containing $x$, then $Y$ and $
\beta_{d,d-1}(Y)$ are isomorphic; in particular both are smooth.

\

Since $x\in \mbox{Sing }{\mathcal G}$ is a
simple point, there  an index  $n\in {\mathbb Z}_{>0}$ and an element
$f\in I_n$ of order exactly $n$ at $x$. Therefore,
\begin{equation}
\label{singf} \mbox{Sing }{\mathcal G}\subset
\{n-\mbox{fold-points of }f=0\}\subset V(\langle f\rangle)
\end{equation}
and moreover $ {\mathcal C}_{{\mathcal G},x}\subset V(\mbox{In}_xf)$. Let
$\{y_1,\ldots, y_{d-1}\}$  be a regular system of parameters in
${\mathcal O}_{V^{(d-1)},x_1}$. Since $\beta_{d,d-1}: V^{(d)}\to
V^{(d-1)}$ is smooth, $\{y_1,\ldots, y_{d-1}\}$ can be extended to a
regular system of parameters $\{y_1,\ldots, y_{d-1}, Z\}$ in
${\mathcal O}_{V^{(d)},x}$.

\

The condition of transversality imposed in Definition \ref{def83}
(iii) ensures that $f\in I_n$ can be chosen so that
$V(\mbox{In}_xf)$ and
$\{\mbox{In}_xy_1=0,\ldots,\mbox{In}_xy_{d-1}=0\}$ intersect only at
the origin of the vector space ${\mathbb T}_{V,x}$. This last
condition can be reformulated by saying that $\mbox{In}_xf\in
\mbox{Gr}_{m_x}({\mathcal O}_{V^{(d)},x})$ is a homogeneous polynomial
of degree $n$ in the variables
$\{\mbox{In}_xy_1,\ldots,\mbox{In}_xy_{d-1},\mbox{In}_xZ \}$, in
which the monomial $(\mbox{In}_xZ)^n$ appears with non-zero
coefficient.

\



Since Weierstrass Preparation Theorem  holds in an \'etale
neighborhood of ${\mathcal O}_{V^{(d-1)},x}$, we may replace
 ${\mathcal O}_{V^{(d-1)},x_{1}}$  and
 ${\mathcal O}_{V^{(d)},x}$ in (\ref{local}) by suitable \'etale neighborhoods if needed,
   and  thus   assume that there is a  regular system of parameters,
$\{y_1,\ldots,y_{d-1}\}\in {\mathcal O}_{V^{(d-1)},x_{1}}$ that
extends to a regular system of parameters in   ${\mathcal
O}_{V^{(d)},x}$, $\{y_1,\ldots,y_{d-1},z\}$ so that
$f=z^n+a_1z^{n-1}+\ldots+a_{n-1}z+a_n$ with $a_i\in \langle
y_1,\ldots,y_{d-1}\rangle^i$. Setting $z$ as $Z$,
\begin{equation}
\label{wfinite} f(Z)=Z^n+a_1Z^{n-1}+\ldots+a_{n-1}Z+a_n\in
{\mathcal O}_{V^{(d-1)},x_1}[Z].
\end{equation}
The map
\begin{equation}\label{finite} {\mathcal O}_{V^{(d-1)},x_1}\to
{\mathcal O}_{V^{(d-1)},x_1}[Z]/\langle f(Z)\rangle
 \end{equation} is  a
finite morphism of local rings  which induces a finite projection
$$\beta: V(f)\to V^{(d-1)},$$
mapping $x\in V(f)$ to,  say, $x_1$. Also, notice that the extension of the
maximal ideal $m_{x_1}\subset {\mathcal O}_{V^{(d-1)},x_1}$ to
${\mathcal O}_{V^{(d-1)},x_1}[Z]/\langle f(Z)\rangle$ is a reduction
of the maximal ideal $M=\langle
\overline{y}_1,\ldots,\overline{y}_{d-1},\overline{Z}\rangle\subset
{\mathcal O}_{V^{(d-1)},x_1}[Z]/\langle f(Z)\rangle$.

\

Let $V_n(f)$ be the closed set of $n$-fold points of  the
hypersurface $V(f)$ in $V^{(d)}$. The map $\beta: V(f)\to V^{(d-1)}$
is  defined  in a neighborhood of $x$ and it    is the restriction
of the smooth morphism $\beta_{d,d-1}: V^{(d)} \to V^{(d-1)}$.  Then
by Zariski's multiplicity formula for projections (see
\ref{Motivation}):

(A) The projection $\beta_{d,d-1}$ induces a bijection between
$V_n(f)$ and its image $\beta(V_n(f))$.

(B) For any irreducible scheme $Y\subset V_n(f)$, the finite map
\begin{equation}\label{eqzki}
\beta: Y \to \beta(Y)
\end{equation}
is  birational.

\

Therefore, by (\ref{singf}) and (A), there is a bijection between
 $\mbox{Sing
}{\mathcal G}$ and $\beta_{d,d-1}(\mbox{Sing }{\mathcal G})$.  From (B) it follows that if $Y\subset
\mbox{Sing} {\mathcal G}$ is an irreducible closed subscheme, then $\beta:
Y\to \beta(Y)$ is a finite birational morphism. Moreover, (A) also
ensures that $\beta: Y\to \beta(Y)$ defines a bijection of the
underlying topological spaces.

\

Assume, in addition, that $x\in Y\subset  \mbox{Sing }{\mathcal G}$
and that $Y$ is a regular center. Then $x$ is the unique point of
$Y$ mapping to $\beta(x)=x_1\in \beta(Y)$, and we claim now that
$(Y,x)$ is \'etale over $(\beta(Y), x_1)$. This together with the
previous properties would show that $\beta(Y)$ is regular at $x_1$,
and that the finite birational map $\beta: Y\to \beta(Y)$ is in fact
an isomorphism in an open neighborhood of $x$.

\

We will argue geometrically  to prove that $\beta: Y\to \beta(Y)$ is
\'etale at $x$. The smooth morphism $\beta_{d,d-1}: V^{(d)}\to
V^{(d-1)}$ induces a linear map of tangent spaces, $d\beta_{d,d-1}:
{\mathbb T}_{V^{(d)},x} \to {\mathbb T}_{V^{(d-1)},x_1}$. The claim
is that $\mbox{ker }d\beta_{d,d-1}\cap {\mathbb
T}_{{Y},x}=\{\overrightarrow{0}\}\subset {\mathbb T}_{V^{(d)},x}$. This follows from
our choice of $f\in I_n$ and the transversality condition in
\ref{def83} (iii). In fact ${\mathbb T}_{{Y},x}\subset
V(\mbox{In}_xf)$ and $V(\mbox{In}_xf)$ and $\mbox{ker
}d\beta_{d,d-1} =\{\mbox{In}_xy_1,\ldots,\mbox{In}_xy_{d-1}\}$
intersect a the origin of the vector space ${\mathbb
T}_{V^{(d)},x}$. This proves that $Y$ and $\beta(Y)$ are isomorphic
in a suitable neighborhood of $x$.

\



The previous discussion also shows that    there is a change of
variable of the form $Z'= Z-a$ in ${\mathcal O}_{V^{(d-1)},x_1}[Z]$,
for a suitable $a\in {\mathcal O}_{V^{(d-1)},x_1}$, such  that
$I(Y)_x=\langle Z', v_1,\dots ,v_s \rangle$, where $\{ v_1,\dots
,v_s \}$ is part of a regular system of parameters at $ {\mathcal
O}_{V^{(d-1)},x_1}$ (see also the proof of Theorem  \ref{1theo}). In
fact, if $\overline{Z}$ denotes the restriction of $Z$ to $Y$, then
there is an element $a\in {\mathcal O}_{V^{(d-1)},x_1}$ which
restricts to the same function on ${\mathcal
O}_{\beta(Y),x_1}={\mathcal O}_{Y,x}$. So $Z-a$ will vanish along
$Y$, and the claim follows from this fact.

}
\end{Paragraph}

\begin{Remark}\label{opentran} {\rm The local smooth projection constructed 
in \ref{localproj} and the
arguments and results described in \ref{Zariski} are also valid in
an open neighborhood of $x$ in $\mbox{Sing }{\mathcal G}$. To show
this it is enough to prove that the transversality condition on
$f$ holds   in an open neighborhood of $x$. To this end, we first note that the transversality 
can be expressed in terms of differential operators taht are relative 
to  $\beta_{d,d-1}$.  In fact  the operators
$\Delta^k$  defined in (\ref{p26})  are also defined as ${\mathcal O}_{V^{(d-1)},x_1}$-relative
differential operators on ${\mathcal O}_{V^{(d)},x}$, say:
\begin{equation}\label{eqion1}
\Delta^k: \mathcal{O}_{V^{(d)},x}\to \mathcal{O}_{V^{(d)},x}.
\end{equation}
Here the inclusion $k[Y_1,\dots,Y_n] \subset k[Y_1,\dots,Y_n]
[Z]$ is replaced  by the inclusion of regular rings
$\mathcal{O}_{V^{(d-1)},x_1}\subset \mathcal{O}_{V^{(d)},x}$, and
  $Z$ is as in (\ref{Zariski}).

\

One can check that the condition of transversality imposed on $f\in
I_n$ at the closed point $x\in \mbox{Sing }{\mathcal G}$ can also be
formulated by requiring  that $\Delta^n(f)$ be a unit at the regular
ring ${\mathcal O}_{V^{(n)},x}$, or, formally,  that:

 \begin{equation} \label{qeab}
 \Delta^n(f)(x)\neq 0,
 \end{equation}
 which also shows that if the geometric condition in Definition \ref{def83} (iii) holds at $x$,
 then it  holds for all
 singular points in an open neighborhood
 of $x$.
}
\end{Remark}

\begin{Remark}\label{rk85} {\rm Fix a polynomial ring $S[Z]$. A morphism,
$$Tay:    S[Z]   \longrightarrow  S[Z,T]; \mbox{ }
   Z   \to   Z+T, $$
and   operators $\Delta^k: S[Z]\to S[Z]$ are defined by setting
$$Tay(G(Z))=\sum_{k\geq 0}\Delta^k(G(Z))T^k.$$
Each $\Delta^k$ is a differential operator of order $k$ over the
ring $S$, and furthermore, for each positive integer $N$,
$\{\Delta^k,  k=0,1, \dots N\}$ is a basis of $Diff^N(S[Z]/S)$,
the free $S$-module of $S$-differential operators of order $N$.

\

Consider  a finite number of monic polynomials, say
$$f_i(Z)=Z^{n_i}+a^i_1Z^{n_i-1}+\ldots+a^i_{n_i}, \\\ i=1, \dots, r,$$
and define a subalgebra  of $S[Z][W]$ of the form
$$S[Z][\{f_i(Z)W^{n_i}, {i=1, \dots ,r}\}].$$
In general this Rees algebra will not be compatible with
$S$-differential operators in the sense of Definition \ref{def82},
(ii). However in \cite[Theorem 2.9]{integraldifferential}  it is
shown that there is a smallest extension of this algebra to one
having this property, and such extension is
\begin{equation}\label{eq85b}
S[Z][\{f_i(Z)W^{n_i},\{\Delta^k(f_i(Z))W^{n_i-k}\}_{k=1,\ldots,n_i-1}\}_{i=1, \dots ,r}].
\end{equation}
But  this extension is, in turn, the pull-back of an algebra in the universal setting in (\ref{eqrtay}),  by a suitable morphism on $S[W]$ (see Theorem
\ref{elimRamMult} and  \ref{comeafter}). A setting like this can  be  always assumed to hold by 
Weirstrass Preparation Theorem. }
\end{Remark}

\begin{Paragraph} \label{eliminationalg}{\rm {\bf The elimination algebra ${\mathcal R}_{{\mathcal
G}_{\beta_{d,d-1}}}$.} Let $V^{(d)}$ be a $d$-dimensional smooth
scheme over a perfect field $k$, let ${\mathcal G}=I_nW^n\subset {\mathcal
O}_{V^{(d)}}$ be a Rees algebra, and assume that $x\in \mbox{Sing
}{\mathcal G}$ is a simple closed point not contained in any component
 of codimension one of $\mbox{Sing }{\mathcal G}$.
Construct
  a smooth morphism $\beta_{d,d-1}: V^{(d)}\to V^{(d-1)}$
to some smooth $(d-1)$-dimensional scheme transversal to
${\mathcal G}$ in a neighborhood of $x$ (this can be done, for
instance, following the arguments given in \ref{localproj}). If in
addition ${\mathcal G}$ is a $\beta_{d,d-1}$-relative differential
algebra, then   $\beta_{d,d-1}: V^{(d)}\to V^{(d-1)}$ is locally
${\mathcal G}$-admissible at $x$. In this case   an {\em elimination
algebra}
$${\mathcal R}_{{\mathcal
G}_{\beta_{d,d-1}}}\subset {\mathcal O}_{V^{(d-1)},x_1}[W]$$ can be
defined (see \cite[1.25, Definitions 1.42 and 4.10]{hpositive}).  To
do so,  first   choose a positive integer $n$, and an element $f\in
I_n$ of order $n$ at ${\mathcal O}_{V^{(d)},x}$, and then produce a
monic polynomial $f(Z) \in I_n$ as in (\ref{wfinite}) in a suitable
\'etale neighborhood of $x$. Then, it can be checked that, up to
integral closure, we may assume that ${\mathcal G}$ is as in
(\ref{eq85b}), for $S={\mathcal O}_{V^{(d-1)},x_1}$, and suitable
monic polynomials $f_i(Z)$, $i=1, \dots ,r$. In particular,
${\mathcal G}$ is locally (and up to integral closure) the pull-back
of the universal algebra  so we define  ${\mathcal R}_{{\mathcal
G}_{\beta_{d,d-1}}}\subset {\mathcal O}_{V^{(d-1)},x_1}[W]$
following the procedures indicated in Theorem  \ref{elimRamMult} and
\ref{comeafter}.

\

This elimination algebra depends on the projection $\beta_{d,d-1}$
but by construction it does not depend on the choice of $f$ once
the projection is fixed,  and it satisfies the  following
conditions:

\begin{enumerate}
\item[i.] The inclusion
$\beta^*_{d,d-1}:{\mathcal O}_{V^{(d-1)},{x_1}}\to {\mathcal
O}_{V^{(d)},{x}}$ induces  an inclusion of Rees algebras
${\mathcal R}_{{\mathcal G},\beta_{d,d-1}}\subset {\mathcal G}$
(this follows now from (\ref{eqrtay}); see  also \cite[Theorem
4.13]{hpositive}).
\item[ii.] If ${\mathcal G}$ is a differential Rees algebra, then   ${\mathcal R}_{{\mathcal
G},\beta_{d,d-1}}$ is  a differential Rees algebra.
\item[iii.] There is an inclusion of  closed subsets
$$\beta_{d,d-1}(\mbox{Sing
}{\mathcal G})\subset \mbox{Sing }{\mathcal R}_{{\mathcal
G},\beta_{d,d-1}}$$ and equality holds if ${\mathcal G}$ is a
differential Rees algebra, or if $\mbox{char }k$ is coprime 
with the degree of $f(Z)$ (cf. \cite[Corollary 4.12]{hpositive}).
\item[iv.] The order of ${\mathcal R}_{{\mathcal
G},\beta_{d,d-1}}$ at $x_1$ does not depend on the projection, in
other words, $\mbox{ord}_{x_1}{\mathcal R}_{{\mathcal
G},\beta_{d,d-1}}$  is independent of $\beta_{d,d-1}$ (see
\cite[Theorem 5.5]{hpositive}).
\end{enumerate}}
\end{Paragraph}

\begin{Paragraph}\label{explicitdescription} {\bf Another description of ${\mathcal R}_{{\mathcal
G},\beta_{d,d-1}}$.} {\rm Assume that  locally, in a neighborhood
of a simple point $x$, a Rees algebra ${\mathcal G}=\oplus_nI_nW^n$ is
generated by
$$\{f_{n_1}W^{n_1},\ldots,f_{n_s}W^{n_s}\}.$$ With  the same notation as in
\ref{Zariski}, we can   choose an element $fW^n$, assuming  that
$f$ has order $n$ in ${\mathcal O}_{V^{(d)},x}$. We can  also
assume that, after multiplying by a unit,
$f=F(Z)=Z^n+a_1Z^{n-1}+\ldots+a_{n-1}Z+a_n\in {\mathcal
O}_{V^{(d-1)},x_1}[Z]$.

\

Note that multiplying  by an element $f_{n_i}W^{n_i}$ induces an
endomorphism
$$L_{f_{n_i}}: \left({\mathcal
O}_{V^{(d-1)},x_1}[Z]/\langle F(Z)\rangle \right)[W]\to
\left({\mathcal O}_{V^{(d-1)},x_1}[Z]/\langle
F(Z)\rangle\right)[W].$$ Since ${\mathcal
O}_{V^{(d-1)},x_1}[Z]/\langle F(Z)\rangle)[W]$ is a free
${\mathcal O}_{V^{(d-1)},x_1}[W]$-module of rank $n$,  each
endomorphism $L_{f_{n_i}}$ has a characteristic polynomial of
degree $n$, say, 
\begin{equation}
\label{coefficients} T^n+g_{1,n_i}T^{n-1}+\ldots+g_{n,n_i}. \
\end{equation} 
It can be proved that for $i=1,\ldots,n$, each 
 $g_{j,n_i}\in {\mathcal O}_{V^{(d-1)},x_1}[W]$ is homogeneous   and that the elimination algebra ${\mathcal
R}_{{\mathcal G},\beta_{d,d-1}}$ is generated by these coefficients
up to integral closure (see \cite[Corollary 4.12]{hpositive} and
Example \ref{canguro} for a   computation in a concrete example).
Furthermore, in a suitable neighborhood of $x$, and up to integral
closure:
\begin{equation}\label{eq876}
{\mathcal G}={\mathcal O}_{V}[\Delta^e(F(Z))W^{n-e},0\leq e\leq
n-1]\odot \beta_{d,d-1}^*({\mathcal R}_{{\mathcal
G},\beta_{d,d-1}}),
\end{equation}
where the right hand side is the smallest subalgebra in ${\mathcal
O}_{V^{(d)}}[W]$ containing both algebras (cf.\cite{B}, or \cite{BVV}). Here the 
$\Delta^e$ are the  differential operators introduced in (\ref{p26}) or in Remark \ref{rk85}. }

\end{Paragraph}

\begin{Paragraph}
{\rm  {\bf Elimination algebras and integral closure.} An
important property of this form of elimination is its link with
integral closure of graded algebras. Using  the same notation as in
\ref{eliminationalg}, consider the following diagram:
$$\begin{array}{rcl}
{\mathcal O}_{V^{(d)},x}[W] & \stackrel{\gamma^*}{\longrightarrow}
& {\mathcal O}_{V^{(d)},x}/\langle
f_n\rangle [W]\simeq {\mathcal O}_{V^{(d-1)},x_1}[Z]/\langle F(Z)  \rangle  [W]\\{\large }
\beta_{d,d-1}^*\uparrow & \nearrow & \\
{\mathcal O}_{V^{(d-1)},x_1}[W]& & \end{array}$$ where $\gamma^*$ denotes the natural restriction. Then the image of
${\mathcal R}_{{\mathcal G},\beta_{d,d-1}}$ in  ${\mathcal
O}_{V^{(d-1)},x_1}[Z]/\langle F(Z)\rangle [W]$ is contained in
$\gamma^*({\mathcal G})$, and this defines a finite extention of graded rings
   (see \cite[Theorem 4.11]{hpositive}). This Theorem also asserts that if an
inclusion of Rees algebras ${\mathcal G}\subset {\mathcal
G}^{\prime}$ is finite, then ${\mathcal R}_{{\mathcal
G},\beta_{d,d-1}}\subset {\mathcal R}_{{\mathcal
G}^{\prime},\beta_{d,d-1}}$ is finite.}
\end{Paragraph}

\begin{Paragraph} {\bf Notation.} {\rm In what
follows, given a
     Rees algebra ${\mathcal G}={\mathcal G}^{(d)} $ on a $d$-dimensional smooth scheme $V^{(d)}$
     of finite type over a field $k$,
   we will refer to an elimination algebra as  ${\mathcal R}_{{\mathcal G},\beta_{d,d-1}}$ if
   we need to emphasize the projection, or just as ${\mathcal G}^{(d-1)}\subset {\mathcal O}_{V^{(d-1)}}[W]$
   if the choice of the projection is not relevant in the discussion.}
   \end{Paragraph}

\begin{Paragraph}
{\rm  {\bf Elimination algebras and the $\tau$-invariant.}
\label{elagatai} The equality   in (\ref{eq876}) is to be considered 
 up to integral closure (both algebras have the same
integral closure). A consequence of Theorem \ref{tauweak} is that
the $\tau$-invariant at  a point is well defined up to integral
closure.

\

In \cite{B} it is proven that if  ${\mathcal G}$ is a differential
algebra, then $\tau_{{{\mathcal R}_{{\mathcal
G},\beta_{d,d-1}}},x_1}=\tau_{{\mathcal G},x}-1$. In summary, this proof shows    that (\ref{eq876})  holds for any $n$ and
any $f_n=F(Z)=Z^n+a_1Z^{n-1}+\ldots+a_{n-1}Z+a_n\in {\mathcal
O}_{V^{(d-1),x_1}}[Z]$, with the only conditions that $f_n\in I_n$
have order $n$ at the local ring ${\mathcal O}_{V^{(d)},x}$ and be transversal to 
$\beta_{d,d-1}$. 
Finally, for the case of differential Rees algebras we may choose $f_n$
so that the initial form defines a linear subspace of (codimension
one) in  $\mathbb{T}_{V^{(d)},x}$.

\

So, in general, if ${\mathcal G}$ is of codimensional type $\geq
e\geq 1$ in a neighborhood of $x$ (i.e., if  $\tau_{{\mathcal
G},x}\geq e$ in $U\subset \mbox{Sing }{\mathcal G}$) then we can
expect to iterate the arguments in \ref{localproj} $e$-times, and
 a sequence of local  projections can be defined:
$$\begin{array}{ccccccc}
V^{(d)} & \stackrel{\beta_{d,d-1}}{\longrightarrow} &  V^{(d-1)} & \to & \ldots &
\stackrel{\beta_{d-(e-1),d-e}}{\longrightarrow} & V^{(d-e)}\\
x=x_0 & \to & x_1 & \to & \ldots & \to & x_e,
\end{array}$$
which by composition induces a local smooth projection from $V^{(d)}$ to
some $(d-e)$-dimensional smooth space  $V^{(d-e)}$. In this way,
by iteration,  we can define elimination algebras $${\mathcal
G}^{(d-1)}\subset {\mathcal O}_{V^{(d-1)}}[W],\ldots, {\mathcal
G}^{(d-e)}\subset {\mathcal O}_{V^{(d-e)}}[W]$$ if for each
$i=1,\ldots,e$, the projection
$$\beta_{d-(i-1),d-i}: V^{(d-(i-1))} \to   V^{(d-i)}$$ is ${\mathcal G}^{(d-(i-1))}$-admissible locally at
$x_{i-1}$.  By  \cite[Corollary 4.12]{hpositive}, there is an
inclusion of closed subsets
$$\beta_{d-(i-1),d-i}(\mbox{Sing }{\mathcal G}^{(d-(i-1))})\subset \mbox{Sing }{\mathcal G}^{(d-i)},$$ which is an equality when
${\mathcal G}^{(d-(i-1))}$ is a differential Rees algebra for
$i=1,\ldots,e$.}
 \end{Paragraph}

\begin{Definition}
\label{def838} {\rm Let ${\mathcal G}$ be a Rees algebra on a
smooth $d$-dimensional scheme $V^{(d)}$ over a perfect field $k$ and  let
$x\in \mbox{Sing }{\mathcal G}$ be a simple point with
$\tau_{{\mathcal G}^{(d)},x}\geq e$. We will say that a local
projection to a smooth $(d-e)$-dimensional scheme over $k$
$$\begin{array}{rrcl}
\beta_{d,d-e}: & V^{(d)}  & \to & V^{(d-e)}\\
 & x & \to & x_e
\end{array}$$ is {\em locally ${\mathcal G}$-admissible at $x$} if
it factorizes as a sequence of  local ${\mathcal
G}^{(d-i)}$-admissible projections as in Definition \ref{def83},
$$\begin{array}{ccccccc}
 (V^{(d)},x_0=x)  & \stackrel{\beta_{d,d-1}}{\longrightarrow} & (V^{(d-1)},x_1) & \to & \ldots &
 \stackrel{\beta_{d-(e-1),d-e}}{\longrightarrow} &  (V^{(d-e)},x_e) \\
  {\mathcal G}^{(d)} &   & {\mathcal G}^{(d-1)} &  & \ldots &   & {\mathcal
 G}^{(d-e)},
\end{array}$$
where for $i=1,\ldots,e$, each ${\mathcal G}^{(d-i)}\subset
{\mathcal O}_{V^{(d-i)},x_i}[W]$ is the elimination algebra of
${\mathcal G}^{(d-(i-1))}\subset {\mathcal
O}_{V^{(d-(i-1))},x_{i-1}}[W]$, and
$\beta_{d-(i-1),d-i}(x_{i-1})=x_i$.}

\end{Definition}

\section{Elimination algebras and permissible monoidal
transformations}\label{eliminationmonoidal} The purpose of this
section is to study the behavior of admissible projections under
permissible monoidal transformations (see Definition \ref{def83}),
a result that will  play a key role in the inductive construction
of the function in Theorem \ref{ordertau}.

\begin{Theorem}\label{1theo}
Let ${\mathcal G}^{(d)}$ be a Rees algebra on a smooth
$d$-dimensional scheme  $V^{(d)}$ over a perfect field $k$,  and let $x\in \mbox{Sing
}{\mathcal G}^{(d)}$ be a simple point (i.e., $\tau_{{\mathcal
G}^{(d)},x}\geq 1$). Suppose that a local ${\mathcal
G}^{(d)}$-admissible projection  is given, defining  an
elimination algebra:
$$
\begin{array}{rrcl} \beta_{d,d-1}:& (V^{(d)},x)& \longrightarrow & (V^{(d-1)},x_1)\\
& {\mathcal G}^{(d)} &   & {\mathcal G}^{(d-1)}.
\end{array}$$
Let $Y\subset \mbox{Sing }{\mathcal G}^{(d)}$ be a permissible
center. Then, locally in a neighborhood of  $x$:
\begin{enumerate}
\item[(i)] The closed set
$\beta_{d,d-1}(Y)\subset \mbox{Sing }{\mathcal G}^{(d-1)}\subset
V^{(d-1)}$ is   a permissible center for ${\mathcal G}^{(d-1)}$.
\item[(ii)] Consider  the monoidal transformations on
$V^{(d)}$ and $V^{(d-1)}$ with centers $Y$ and $\beta_{d,d-1}(Y)$
respectively. Then there is a projection $\beta_{d,d-1}^{\prime}$
defined in  a suitable open set, and a commutative diagram of
projections, 
$$
\xymatrix@R=0pc@C=0pc{ (V^{(d)},x)  & & & & & \quad U\subset
V^{(d)^{\prime}}\ \ar[lllll]_{\pi^{(d)}}\\
{\mathcal G}^{(d)}\ar[dddd] ^{\beta_{d,d-1}}   & &  &  & & {\mathcal G}^{{(d)}^{\prime}}\ar[dddd]^{\beta_{d,d-1}^{\prime}} \\
\\
& & & \circlearrowleft \\
\\
(V^{(d-1)},x_1) & &  &  & & V^{(d-1)^{\prime}} \ar[lllll]_{\pi^{(d-1)}}\\
{\mathcal G}^{(d-1)}  & &  &  & & {\mathcal G}^{(d-1)^{\prime}},
}$$
\end{enumerate}
where ${\mathcal G}^{{(d)}^{\prime}}$ denotes the weighted transform of 
${\mathcal G}^{{(d)}}$ in $V^{(d)^{\prime}}$, and 
${\mathcal G}^{{(d-1)}^{\prime}}$ denotes the weighted transform of 
${\mathcal G}^{{(d-1)}}$ in $V^{(d-1)^{\prime}}$. Furthermore, if $x'\in \mbox{Sing }{\mathcal G}^{(d)^\prime}\neq
\emptyset$ maps to $x$, then,
\begin{enumerate} \item[(a)]  The
projection $V^{(d)^{\prime}}\to V^{(d-1)^{\prime}}$ is ${\mathcal
G}^{(d)^\prime}$-admissible locally at $x^{\prime}$ (see
Definition \ref{def83}). In particular ${\mathcal G}^{(d)^\prime}$
is a $\beta_{d,d-1}^{\prime}$-relative-differential Rees algebra,
defining an elimination algebra  ${\mathcal G}^{\prime(d-1)}$.
\item[(b)] Let $x^{\prime}_1=\beta_{d,d-1}^{\prime}(x^{\prime})$.
 Locally in an open neighborhood of $x^{\prime}_1$, there is a natural inclusion ${\mathcal
G}^{(d-1)^{\prime}}\subset {\mathcal G}^{\prime(d-1)}$   which is a 
 finite extension.

\end{enumerate}
\end{Theorem}

To prove the Theorem we need  to study the behavior  of elimination
algebras in the universal case, which is the purpose of the next
lines. The proof of the Theorem is given in \ref{proofblowing}.

\begin{Paragraph} \label{elimono} {\bf Monoidal transformations and weak
transforms.} {\rm Our goal is to understand:
\begin{enumerate}
\item[i.] How   elimination algebras behave under permissible monoidal
transformations.
\item[ii.] The behavior of differential Rees algebras
under permissible monoidal transformations.
\end{enumerate}
Both points will be treated in the universal case.

\

Let $F_n(Z)=(Z-Y_1)\cdot (Z-Y_2)\cdots (Z-Y_n)$ be the universal
monic polynomial of degree $n$, and let $V$ be a new variable.  We
would like to make sense of  the expression
$$\Big(\frac{1}{V}\Big)^nF_n(Z)\in k[Y_1,\dots,Y_n][V,V^{-1}][Z]$$
which  corresponds  to the notion of weighted transform of a degree-$n$
element in a Rees algebra.

\

Notice that
$$\left(\frac{1}{V}\right)^nF_n(Z)=\left(\frac{Z}{V}-\frac{Y_1}{V}\right)\cdot \left(\frac{Z}{V}-\frac{Y_2}{V}\right)
\cdots \left(\frac{Z}{V}-\frac{Y_n}{V}\right)\in
k[Y_1,\dots,Y_n][V,V^{-1}][Z].$$

\

Set
$$F'_n\left(\frac{Z}{V}\right):=\left(\frac{Z}{V}-\frac{Y_1}{V}\right)\cdot
\left(\frac{Z}{V}-\frac{Y_2}{V}\right) \cdots
\left(\frac{Z}{V}-\frac{Y_n}{V}\right).$$ Then
$F'_n\left(\frac{Z}{V}\right)=\frac{1}{V^n}F_n(Z)$ is  a monic polynomial in
the ring
 $k\left[\frac{Y_1}{V}, \frac{Y_2}{V},\dots ,\frac{Y_n}{V}\right]\left[\frac{Z}{V}\right]$.

\

Let  $\Delta^{{k}}_1$ be a differential operator on
$k[\frac{Y_1}{V}, \frac{Y_2}{V},\dots ,\frac{Y_n}{V}][\frac{Z}{V}]$
 relative to   $k[\frac{Y_1}{V}, \frac{Y_2}{V},\dots
,\frac{Y_n}{V}]$ for some $k<n$. Then by (\ref{differentialsim}),
\begin{equation}
\label{comportamientodiff}
\Delta^{{k}}_1\left(F'_n\left(\frac{Z}{V}\right)\right)=
\frac{1}{V^{n-k}}\cdot\Delta^{k}(F_n(Z)).
\end{equation}

\

Now let $S_n$ act on $k[Y_1,\dots,Y_n][V,V^{-1}][Z]$ by permuting
 the variables $Y_i$ and fixing both $V$ and $Z$. In this way $S_n$ acts
also by permutation on the variables $\frac{Y_i}{V}$ for
$i=1,\ldots,n$  and   fixes  $\frac{Z}{V}$. Therefore
$$k\left[\frac{Y_i}{V}-\frac{Y_j}{V}\right]^{S_n}=k[H'_{m_1},\dots,H'_{m_r}], $$
where
\begin{equation}\label{weaktransfelimination}
H'_{m_i}=\frac{1}{V^{m_i}}\cdot H_{m_i}\end{equation} for $H_{m_i}$ and $m_i$ 
as in (\ref{universalgenerators}).}

\end{Paragraph}

\begin{Paragraph} \label{proofblowing} {Proof of Theorem \ref{1theo}.} {\rm (i) There is an inclusion
$\beta_{d,d-1}(\mbox{Sing }{\mathcal G})\subset \mbox{Sing
}{\mathcal G}^{(d-1)}$ (see \ref{eliminationalg} (iii)), mapping
$Y$ to $\beta(Y)$, which is an isomorphism (see \ref{Zariski}).
Since
 $Y$ and $ \beta(Y)$ are isomorphic, then $\beta(Y)$ is a regular
permissible center for ${\mathcal G}^{(d-1)}$.

\

(ii) We claim that the  monoidal transformations of $V^{(d)}$ and
$V^{(d-1)}$ with centers $Y$ and $\beta_{d,d-1}(Y)$ respectively,
produce, in a suitable open neighborhood of $\mbox{Sing }{\mathcal G}^{(d)^{\prime}}$, say 
 $U\subset  V^{(d)^{\prime}}$, a
commutative diagram of projections:
$$
\xymatrix@R=0pc@C=0pc{ (V^{(d)},x) \ar[dddd] ^{\beta_{d,d-1}}  & &
& & & \quad (U\subset
V^{(d)^{\prime}},x^{\prime}) \ \ar[lllll]_{\pi^{(d)}} \ar[dddd]^{\beta_{d,d-1}^{\prime}} \\
\\
& & & \circlearrowleft \\
\\
(V^{(d-1)},x_1) & &  &  & & (V^{(d-1)^{\prime}}, x_1^{\prime})
\ar[lllll]_{\pi^{(d-1)}}. }$$

 Set ${\mathcal G}=\oplus I_nW^n\subset
{\mathcal O}_{V^{(d)}}[W]$ and ${\mathcal G}^{\prime}=\oplus
I_n^{\prime}W^n\subset {\mathcal O}_{V^{(d)^{\prime}}}[W]$. To
prove (ii) we argue as in \ref{Zariski}. First, choose an integer $n$ and
select an element $f\in I_n$ of order $n$ transversal to
$\beta_{d,d-1}$ at $x$, and consider its  weighted transform,
$f^{\prime}\in I_n^{\prime}$. The hypothesis of transversality on
$f$ can be reformulated by saying  that the relative differential
operator of order $n$, $\Delta^n\in Diff^n{\mathcal
O}_{V^{(d)}/V^{(d-1)}}$, is such that $\Delta^nf$ is a unit in a
neighborhood of $x$ (see \ref{opentran}, specially (\ref{qeab})).

\

We claim that the law of transformation of relative differentials
in (\ref{comportamientodiff}) specializes to show   that   there
is a relative differential operator $\Delta^{\prime^n}\in
Diff^n{\mathcal O}_{V^{(d)^{\prime}}/V^{(d-1)^{\prime}}}$ such
that $\Delta^{\prime^n}(f^{\prime})$ is a unit in a neighborhood
of $x^\prime$. Therefore  $x^{\prime}$ is a simple point (i.e.,
$\tau_{{\mathcal G}^{(d)^{\prime}},x^{\prime}}\geq 1$). We also
claim that the law of transformation in (\ref{comportamientodiff})
already shows that ${\mathcal G}^{(d)^{\prime}}$ is a
$\beta^{\prime}_{d,d-1}$-relative differential Rees algebra locally at
$x'$.

\

To clarify these points note first that we may assume that, after
multiplying  by a unit, $f$ is monic of degree $n$, i.e.,
$f=Z^n+a_1Z^{n-1}+\ldots+a_{n-1}Z+a_n\in {\mathcal
O}_{V^{(d-1),x_1}}[Z]$ in an \'etale neighborhood of the point
$x$.  According to   (i),   the center $Y$ maps isomorphically to $\beta(Y)$; in
particular the class (restriction) on $Y$ of any element of
${\mathcal O}_{V^{(d)},x}$ is also the class of an element of
${\mathcal O}_{V^{(d-1)},x_1}$. Thus, after a suitable change of
variable of the form $Z-\alpha$, $\alpha \in {\mathcal
O}_{V^{(d-1)},x_1}$, we may assume that $Z$ vanishes identically
along $Y$, and that $I(Y)=\langle Z, y_1, \dots , y_s\rangle $,
where $\{ y_1,\dots , y_s\}$ is part of a regular system of
parameters  in ${\mathcal O}_{V^{(d-1)},x_1}$, and each
coefficient $a_i$ has order $\geq i$ along the regular center $Y$.
The closed set $\mbox{Sing }{\mathcal G}^{(d)}$ is included in the
closed set of $n$-fold points of the hypersurface $V(\langle f
\rangle)$, and $\mbox{Sing }{\mathcal G}'^{(d)}$ is included in
the closed set of  $n$-fold points of $V(\langle f' \rangle)$.

\

Consider the open set $U\subset V^{(d)^{\prime}}$ which is the
union of the charts $\mbox{Spec}\left({\mathcal
O}_{V^{(d)}}\left[\frac{Z}{y_j}, \frac{y_1}{y_j},\dots ,
\frac{y_s}{y_j}\right]\right)$,  for $j=1,\dots , s$. Note that $V(\langle f^{\prime}\rangle)\subset 
V^{(d)^{\prime}}$ is the strict transform of $V(\langle f\rangle)\subset V^{(d)}$, and 
that $V(\langle f^{\prime}\rangle)\subset U$.       The
inclusions ${\mathcal O}_{V^{(d-1)}}\left[\frac{y_1}{y_j},\dots ,
\frac{y_s}{y_j}\right]\subset {\mathcal
O}_{V^{(d)}}\left[\frac{Z}{y_j}, \frac{y_1}{y_j},\dots ,
\frac{y_s}{y_j}\right]$ define $$\beta_{d,d-1}^{\prime}: U\to
V^{(d-1)^{\prime}},$$ as above. The point $x'\in \mbox{Sing
}{\mathcal G}^{(d)^\prime}$ is included the $n$-fold points of
$f'$, and
$$f'=\frac{Z^n}{y_j^n}+\frac{a_1}{y_j}\frac{Z^{n-1}}{y_j^{n-1}}+
\ldots+\frac{a_{n-1}}{y^{n-1}_j}\frac{Z}{y_j}+\frac{a_{n}}{y^{n}_j}\in
{\mathcal
O}_{V^{(d-1)^{\prime}},x_1^{\prime}}\left[\frac{Z}{y_j}\right].$$

\

Moreover, as $x^{\prime}$ is a point of multiplicity $n$ of $f^{\prime}=0$,
the residue fields of ${\mathcal O}_{V^{(d)^{\prime}},x^{\prime}}$
and of ${\mathcal O}_{V^{(d-1)^{\prime}},x_1^{\prime}}$ are the
same. Thus there is an element $a\in {\mathcal
O}_{V^{(d-1)^{\prime}},x_1^{\prime}}$ such that
$Z_1=\frac{Z}{y_j}-a$ vanishes at $x^{\prime}$, and $f^{\prime}$ is
a monic polynomial of degree $n$ in ${\mathcal
O}_{V^{(d-1)^{\prime}},x_1^{\prime}}[Z_1]$. So,  if
$\{z_1,\ldots, z_{d-1}\}$ is a regular system of parameters at
${\mathcal O}_{V^{(d-1)^{\prime}},x_1^{\prime}}$, then
we can choose $\{z_1,\ldots, z_{d-1}, Z_1\}$ so that  it be  a regular system of parameters at
${\mathcal O}_{V^{(d)^{\prime}},x^{\prime}}$. Now the local
arguments in  \ref{localproj} and \ref{eliminationalg} can be
repeated for
$$
\begin{array}{rrcl} \beta_{d,d-1}^{\prime}:& (V^{(d)^{\prime}},x^{\prime})& \longrightarrow &
(V^{(d-1)^{\prime}},x_1^{\prime})\\
& {\mathcal G}^{(d)^{\prime}} & \to & {\mathcal G}^{(d-1)^{\prime}}
\end{array}$$
to prove the statement in (b). More precisely notice that:

\begin{itemize}
\item Up to integral closure ${\mathcal G}^{(d)}$ is generated by monic polynomials in the setting of (\ref{eq85b}), as indicated in \ref{eliminationalg}.
\item  The weighted transforms of the local generators of
${\mathcal G}^{(d)}$ generate ${\mathcal G}^{{(d)}^{\prime}}$ (see
Proposition \ref{localt}).

\item  The closed point  $x^{\prime}$ is contained in
$V(\langle f'\rangle )\subset V^{(d)^{\prime}}$.

\item We claim now that  the weighted transform of
$f$ in ${\mathcal G}^{{(d)}^{\prime}}$, $f^{\prime}$, can be used
to define the elimination algebra in a neighborhood of
$x^{\prime}$ as in \ref{localproj}. In Theorem  \ref{elimRamMult}
and   \ref{comeafter}  it is given  an explicit description of the
elimination algebra as a specialization of the universal algebra
elimination algebra. It follows from (\ref{comportamientodiff})
that   ${\mathcal G}^{\prime}$ is a relative differential Rees algebra.
It also follows   from (\ref{weaktransfelimination}) that up to integral
closure, ${\mathcal G}^{(d-1)^{\prime}}={\mathcal
G}^{{\prime}^{(d-1)}}$.\qed
\end{itemize}
}
\end{Paragraph}

\part{Main Theorem and inductive invariants}\label{partMain}

\section{Main Theorem \ref{ordertau}}\label{sectionorder}

In this section we discuss resolutions of Rees algebras (see
\ref{ragadprs3}),  where the main invariant is  the function
$\mbox{ord }_{}{\mathcal G}$, defined by
$\mbox{ord}_{x_{}}{\mathcal G}$ in \ref{orderRees} for each $x\in
\mbox{Sing }{\mathcal G}$. When the characteristic is zero and
$\mbox{ord}_{x_{}}{\mathcal G}=1$ there is a smooth hypersurface
of maximal contact at $x$, and a new Rees algebra
$\overline{{\mathcal G}}$  is defined along this smooth
hypersurface. In particular a new value
$\mbox{ord}_{x_{}}\overline{{\mathcal G}}$ can be defined. It is
then shown that this value is an invariant; in other words, independent of
the choice of the hypersurface of maximal contact. This result is
the main outcome of the so called Hironaka trick; an alternative
and enlightening proof of this result is due to J.  W\l odarczyk (see
\cite{WLL}).

\

In this work  hypersurfaces of maximal contact are replaced by
suitable projections, and we begin this section by formulating
this result in this setting, in Main Theorem \ref{ordertau}. We
then recall briefly how resolution of Rees algebras can be
achieved by induction (see \ref{ragadprs3}   and \ref{rgn})
parallelling the ideas given in Part \ref{Char0}.

\

 Let
${\mathcal G}^{(d)}$ be a Rees algebra on a smooth $d$-dimensional
scheme $V^{(d)}$  over a perfect field $k$, and let $x\in
\mbox{Sing }{\mathcal G}^{(d)}$ be a closed point with
$\tau_{{\mathcal G}^{(d)},x}\geq e$. Assume that  $x$   is not contained
in any  component of codimension $e$ of $\mbox{Sing }{\mathcal G}$, and
that there are two different admissible projections to a
$(d-e)$-dimensional smooth space (see (\ref{def838})),
$$\begin{array}{rrclcrrcl}
\beta_{1_{d,d-e}}: & V^{(d)}& \longrightarrow &  V^{(d-e)}_1& \ \
\ &
\beta_{2_{d,d-e}}: & V^{(d)}& \longrightarrow &  V^{(d-e)}_2\\
& x & \to & x_{e,1} & & & x & \to & x_{e,2}. \end{array}$$  Then the
question is to compare $\mbox{ord}_{x_{e,1}}{\mathcal G}_1^{(d-e)}$
and $\mbox{ord}_{x_{e,2}}{\mathcal G}_2^{(d-e)}$.

\

In \cite[Theorem 5.5]{hpositive} it is shown that if  $e\geq 1$ then
$$\mbox{ord}_{x_{1,1}}{\mathcal
G}_1^{(d-1)}=\mbox{ord}_{x_{1,2}}{\mathcal G}_2^{(d-1)}.$$ In this
section we generalize this result, which leads to Definition
\ref{deford}. More precisely we prove the following theorem:

\begin{Theorem}[{\bf Main Theorem}]\label{ordertau}
Let $V^{(d)}$ be a $d$-dimensional scheme smooth  over a perfect field
$k$, let ${\mathcal G}^{(d)}\subset {\mathcal O}_{V^{(d)}}[W]$ be
a differential Rees algebra, let  $x\in \mbox{Sing }{\mathcal G}^{(d)}$
be a simple closed point, and let $m\leq \tau_{{\mathcal G},x}$.
Consider two different ${\mathcal G}^{(d)}$-admissible local
projections to some $(d-m)$-dimensional smooth schemes with their
corresponding elimination algebras:
\begin{equation}
\label{transorder}
\begin{array}{rrclcrrcl} \beta_{1_{d,d-m}}: &
(V^{(d)},x) & \longrightarrow &  (V^{(d-m)}_1, x_{m,1})& \ &
\beta_{2_{d,d-m}}: & (V^{(d)},x) & \longrightarrow &  (V^{(d-m)}_2,x_{m,2})\\
& {\mathcal G}^{(d)} & \to & {\mathcal G}^{(d-m)}_1 & & &
{\mathcal G}^{(d)} & \to & {\mathcal G}^{(d-m)}_2.
\end{array}
\end{equation}
Then:
$$\mbox{ord}_{x_{m,1}}{\mathcal
G}_1^{(d-m)}=\mbox{ord}_{x_{m,2}}{\mathcal G}_2^{(d-m)}.$$

\noindent Moreover, if  $V^{(d)} \leftarrow V^{(d)^{\prime}}$ is a
composition of permissible monoidal transformations, $x^{\prime}\in
\mbox{Sing }{\mathcal G}^{(d)^{\prime}}$   a closed point dominating
$x$,  and
$$
\xymatrix@R=0pc@C=0pc{ (V^{(d)},x)  & & & & & \quad (U\subset
V^{(d)^{\prime}},x^{\prime})\ \ar[lllll]\\
{\mathcal G}^{(d)}\ar[dddd]    & &  &  & & {\mathcal G}^{{(d)}^{\prime}}\ar[dddd] \\
\\
& & & \circlearrowleft \\
\\
(V^{(d-m)}_j,x_{m,j}) & &  &  & & (V^{(d-m)^{\prime}}_j,x_{m,j}^{\prime}) \ar[lllll]\\
{\mathcal G}^{(d-m)}_j  & &  &  & & {\mathcal
G}^{(d-m)^{\prime}}_j }$$ is the corresponding  commutative
diagram of elimination algebras and admissible projections for
$j=1,2$,  then
$$\mbox{ord}_{x_{m,1}^{\prime}}{\mathcal
G}_1^{(d-m)^{\prime}}=\mbox{ord}_{x_{m,2}^{\prime}}{\mathcal
G}_2^{(d-m)^{\prime}}.$$

\end{Theorem}

The Theorem provides the following upper semi-continuous functions:
\begin{Definition}\label{deford}
{\rm (i) Let $V^{(d)}$ be a $d$-dimensional scheme smooth  over a
field $k$,  let ${\mathcal G}^{(d)}\subset {\mathcal
O}_{V^{(d)}}[W]$ be a differential Rees algebra, let  $x\in \mbox{Sing
}{\mathcal G}^{(d)}$ be a simple closed point, and let $m\leq
\tau_{{\mathcal G},x}$. Then, in a neighborhood of $x$,  we define
the function
$$\begin{array}{rrcl}\mbox{ord}^{(d-m)}_{{\mathcal G}^{(d)}}:& \mbox{Sing }\mathcal{G}^{(d)}&\to& \mathbb{Q}\\
 & z& \to & \mbox{ord}_{z_m}{{\mathcal G}^{(d-m)}}\end{array}$$
where ${\mathcal G}^{(d-m)}$ is an elimination algebra defined by an
arbitrary $\mathcal{G}^{(d)}$-admissible local smooth projection to some
$(d-m)$-dimensional smooth scheme, $\beta_{d,d-m}:V^{(d)}\to
V^{(d-m)}$, and $z_m=\beta_{d,d-m}(z)$ (notice that the function is
well defined since it does not depend on the projection by Theorem
\ref{ordertau}).

(ii) Let ${\mathcal G}^{(d)}\subset {\mathcal O}_{V^{(d)}}[W]$ be a
differential Rees algebra as in (i), let  $\beta_{d,d-m}: V^{(d)}\to V^{(d-m)}$
be any ${\mathcal G}^{(d)}$-admissible local smooth projection in a
neighborhood of a simple point $x\in \mbox{Sing }{\mathcal
G}^{(d)}$, and let  $V^{(d)} \leftarrow V^{(d)^{\prime}}$ be a
composition of permissible monoidal transformations. Let
$x^{\prime}\in \mbox{Sing }{\mathcal G}^{(d)^{\prime}}$  be a closed
point dominating $x$, and consider the corresponding  commutative
diagram of elimination algebras and admissible projections as in
Theorem \ref{1theo},
$$
\xymatrix@R=0pc@C=0pc{ (V^{(d)},x)  & & & & & \quad U\subset
V^{(d)^{\prime}}\ \ar[lllll]\\
{\mathcal G}^{(d)}\ar[dddd] ^{\beta_{d,d-m}}   & &  &  & & {\mathcal G}^{{(d)}^{\prime}}\ar[dddd]^{\beta_{d,d-m}^{\prime}} \\
\\
& & & \circlearrowleft \\
\\
(V^{(d-m)},x_m) & &  &  & & V^{(d-m)^{\prime}} \ar[lllll]\\
{\mathcal G}^{(d-m)}  & &  &  & & {\mathcal G}^{(d-m)^{\prime}}.
}$$
 Then in a neighborhood of $x^{\prime}$  the function
$$\begin{array}{rrcl}\mbox{ord}^{(d-m)^{\prime}}_{{\mathcal G}^{(d)^{\prime}}}:& \mbox{Sing }\mathcal{G}^{(d)^{\prime}}&\to& \mathbb{Q}\\
 & z^{\prime}& \to & \mbox{ord}_{z^{\prime}_m}{{\mathcal G}^{(d-m)^{\prime}}}\end{array}$$
 with $z_m^{\prime}=\beta_{d,d-m}^{\prime}(z^{\prime})$ is well define since by Theorem \ref{ordertau} it is independent of the projection.
}

\end{Definition}

\

The proof of Theorem \ref{ordertau} will be postponed to the next
section. In the rest of this section we analyze some properties 
  of the functions in \ref{deford}. The study of these 
functions will lead us  to the so-called 
{\em reduction to the monomial case} to be  treated in the coming sections (see Part 4).

\

Theorem \ref{ordertau} is stated for a simple Rees algebra on a
$d$-dimensional smooth scheme $V^{(d)}$. We will indicate
why  simple Rees algebras arise in resolution problems. 

\

There is a dictionary between  Rees algebras and pairs as
indicated in \ref{paraleloragadprs}. Resolution is expresssed in terms of 
pairs, where   
simple pairs (see  \ref{racsim}) play a central role.  The notion of simple 
pair is analogous to the notion of simple Rees algebra.

\begin{Paragraph}{\bf Resolution of Rees algebras. \cite[5.10]{positive}} \label{ragadprs3}{\rm
As   pointed out  in \ref{paraleloragadprs},  there is a strong
link between $(J,b)$ and the Rees algebra
$\mathcal{G}=\mathcal{G}_{(J,b)}$ (see (\ref{eqq34})). So a
sequence of transformations of pairs and basic objects as in
(\ref{Atransfuno}) defines a sequence of transformations or Rees
algebras:
\begin{equation}\label{Atransfuno1}
(V,\mathcal{G},E) \longleftarrow (V_1,\mathcal{G}_1,E_1)\longleftarrow \cdots
\longleftarrow (V_s,\mathcal{G}_s,E_s).
\end{equation}
It follows from our notion of transformation of Rees algebras that each $\mathcal{G}_i=\mathcal{G}_{(J_i,b)}$, so
$$\Sing(\mathcal{G}_i)=\Sing(J_i,b).$$
Furthermore, if $d$ denotes the dimension of $V$, then  the  functions
\begin{equation}\label{eqstdeg1}
\vword^{(d)}_{{\mathcal G}_i} :\Sing(\mathcal{G}_i) \to
\mathbb{Q}
\end{equation}
  are defined with the same properties as in the case of pairs (see \ref{lasat1}).

\

We say that a sequence of transformations,
\begin{equation}\label{Atransfuno2}
(V,\mathcal{G},E) \longleftarrow (V_1,\mathcal{G}_1,E_1)\longleftarrow \cdots
\longleftarrow (V_s,\mathcal{G}_s,E_s).
\end{equation}
is a {\em resolution} of $(V,\mathcal{G},E)$ (or simply a {\em resolution}
of $\mathcal{G}$ if $E$ is empty), if  $
\Sing(\mathcal{G}_s)=\emptyset $. }
\end{Paragraph}

\begin{Paragraph}\label{rgn}{\bf The monomial case.  \cite[6.11]{positive}} {\rm
Let $V^{(d)}$ be a $d$-dimensional scheme smooth  over a perfect field
$k$, let ${\mathcal G}^{}\subset {\mathcal O}_{V^{(d)}}[W]$ be a
differential Rees algebra of codimensional type $\geq m$ (\ref{cdtype}). In
this case the  function
$$\mbox{ord}^{(d-m)}_{{\mathcal G}^{}}:\mbox{Sing }\mathcal{G}^{}\to \mathbb{Q}$$
is described   in Definition \ref{deford}. The  discussion on
basic objects and its resolution, which was  presented in Section
\ref{algoritmicpairs}, also extends to this context and satellite
functions $\vword^{(d-m)}_{{\mathcal G}}$ are defined, with the
property that $\vword^{(d-m)}_{{\mathcal
G}}=\mbox{ord}^{(d-m)}_{{\mathcal G}^{}},$ and  if
\begin{equation}\label{Accro2}
(V^{(d)},\mathcal{G}^{},E) \longleftarrow (V^{(d)}_1,\mathcal{G}^{}_1,E_1)\longleftarrow \cdots
\longleftarrow (V^{(d)}_s,\mathcal{G}^{}_s,E_s).
\end{equation}
is a sequence of monoidal transformations with center $Y_i \subset \vMax  \vword^{(d-m)}_{{\mathcal G}}$, then
$$ \max   \vword^{(d-m)}_{{\mathcal G}}  \geq \max  \vword^{(d-m)}_{{\mathcal G}_1} \geq \dots \geq \max
 \vword^{(d-m)}_{{\mathcal G}_s} .$$
  When this holds we say that $(V_s,\mathcal{G}_s,E_s)$ is in the {\em monomial case} if
$ \max \vword^{(d-m)}_{{\mathcal G}_s}=0$ (here $m$ could be zero). This amounts to saying that 
${\mathcal G}$ is, up to integral closure, of the form ${\mathcal G}_{(J,b)}$ for some monomial 
ideal $J$.

\begin{Remark} \label{char0} {\rm As indicated in the introduction, when the characteristic is zero  the upper semi-continuos functions $\mbox{ord}^{(d-m)}_{{\mathcal G}^{}}$ coincide with the classical $\mbox{ord}^{(d-m)}$ defined for the corresponding pairs. Elimination algebras and coefficient ideals produce the same invariants via the fact that locally admissible projections are nothing but  
restrictions to hypersurfaces of maximal contact in the characteristic zero case 
(this is discuss in full datail  in \cite{mariluz}). }
\end{Remark}}
\end{Paragraph}

\section{Proof of Theorem \ref{ordertau}}\label{sectionorder1}

\noindent {\bf The strategy.} In this section we address the proof
of Theorem \ref{ordertau}. Recall our starting point: we assume the
existence of two locally ${\mathcal G}^{(d)}$-admissible projections
to  $(d-m)$-smooth dimensional schemes (\ref{def838}),
$$\begin{array}{ccccc}
 & & (V^{(d)},x) & & \\
 & & {\mathcal G}^{(d)} & & \\
&\swarrow &  & \searrow & \\
(V^{(d-m)}_1,x_{m,1}) & & & & (V^{(d-m)}_2,x_{m,2})\\
{\mathcal G}_1^{(d-m)} & & & & {\mathcal G}_2^{(d-m)}
\end{array}$$
and using the hypothesis of Theorem \ref{ordertau} we want to show
that  $\mbox{ord}_{x_{m,1}} {\mathcal
G}_1^{(d-m)}=\mbox{ord}_{x_{m,2}} {\mathcal G}_2^{(d-m)}$. To this
end, as indicated  in the following proposition, it will be enough
to find a suitable local ring $(B,m)$ and suitable maps
$${\mathcal O}_{V^{(d-m)}_i,x_{m,i}}\to B$$ so that the images of
${\mathcal G}^{(d-m)}_i$ in $B[W]$ under these maps have the same
integral closure for $i=1,2$. More precisely:

\begin{Proposition}\cite[Lemma 5.7 and Corollary 5.8]{hpositive} \label{preorder} Let
$(B,m)$ be a local ring, let $$(S_1,m_1), (S_2,m_2)\subset (B,m)$$
be two regular local rings, let ${\mathcal H}=\oplus I_kW^k\subset
B[W]$ be a Rees algebra  and let     $${\mathcal H}_1=\oplus
J_{1,k}W^k\subset S_1[W] \ \ \mbox{  and  } \ \ {\mathcal
H}_2=\oplus J_{2,k}W^k\subset S_2[W]$$ be  Rees algebras with
inclusions
$${\mathcal H}_1, {\mathcal H}_2\subset {\mathcal H}.$$
Assume that for $i=1,2$:
\begin{enumerate}
\item[(i)]  The inclusions $S_i\subset B$ are finite and flat extensions
of local rings.
\item[(ii)] The ideals $m_iB\subset B$ are reductions of $m$.
\item[(iii)] The inclusions ${\mathcal H}_i=\oplus
J_{i,k}W^k \subset {\mathcal H}=\oplus I_kW^k$ are  finite for $i=1,2$.
\end{enumerate}
Then $$\mbox{ord}_{S_1}({\mathcal H}_1)=\mbox{ord}_{S_2}({\mathcal
H}_2).$$
\end{Proposition}

The basic idea for the proof of Theorem \ref{ordertau} is that
under its assumptions  we can find a suitable sequence of elements
\begin{equation} \label{suitable} f_1W^{n_1},\ldots,f_mW^{n_m}\in {\mathcal
G}^{(d)}
\end{equation} so that the hypotheses of Proposition \ref{preorder} hold
for:
\begin{itemize}
\item $B={\mathcal O}_{V^{(d)},x}/\langle
f_1,\ldots,f_m\rangle$;
\item ${\mathcal H}$  the image of ${\mathcal G}^{(d)}$ in $B[W]$ under
the natural quotient map ${\mathcal O}_{V^{(d)},x}\to {\mathcal O}_{V^{(d)},x}/\langle f_1,\ldots,f_m\rangle$; \item $S_1={\mathcal
O}_{V^{(d-m)}_1,x_{m,1}}$, $S_2={\mathcal
O}_{V^{(d-m)}_2,x_{m,2}}$;
\item ${\mathcal H}_1={\mathcal G}^{(d-m)}_1$,  ${\mathcal
H}_2={\mathcal G}^{(d-m)}_2$.
\end{itemize}

\

\begin{Paragraph}
\label{ideadeprueba}{\bf Idea of the proof  of Theorem
\ref{ordertau}. }{\rm  Observe that  there are two statements in  Theorem \ref{ordertau}: the first is a result about differential Rees algebras,  while  the second part  is the corresponding  statement for the  weighted transform of a differential Rees algebra after a finite sequence of monoidal transformations.

\

The first part of the Theorem will be proven in two steps: 1 and 2.  In step 1 we will show that differential
algebras contain
  sequences of elements with special properties. This will be  used in step 2 to accomplish  the first part of Theorem \ref{ordertau}.

\

Similarly, the  proof of the second part of the Theorem will be shown  in
two steps:  $1^{\prime}$,  and $2^{\prime}$. In step 1$^{\prime}$ we  will prove that,  after a finite sequence of monoidal  permissible transformations, the  weighted transform of 
a differential Rees algebra contains a
 sequence of elements with special properties. This will be  used in step
2$^{\prime}$, where the second part of the  Theorem  \ref{ordertau}  will be
given.

\

\noindent{\bf Idea of the proof of the first part of Theorem \ref{ordertau}}

\

\noindent {\bf Step 1.}   Assume that    $V^{(d)}$ is  a $d$-dimensional scheme
smooth over a perfect field $k$, and that ${\mathcal G}^{(d)}\subset {\mathcal
O}_{V^{(d)}}[W]$ is a differential Rees algebra. Let  $x\in \mbox{Sing
}{\mathcal G}^{(d)}$ be a simple closed point, and let $m\leq
\tau_{{\mathcal G},x}$. Suppose that
\begin{equation}
\label{projtest}
(V^{(d)},x)\to
(V^{(d-m)},x_{m})
\end{equation} is a ${\mathcal G}^{(d)}$-admissible projection
locally at $x$ (see Definition \ref{def838}). Under these assumptions   there there is factorization of  (\ref{projtest}) into local
admissible projections, together with elimination algebras, and we will show that 
there are elements
\begin{equation}
\label{flag}
\begin{array}{ccccccccc} (V^{(d)},x) & \to &
(V^{(d-1)},x_1)& \to & \ldots & \to &
(V^{(d-(m-1))},x_{m-1}) & \to & (V^{(d-m)},x_{m})\\
f_1\in {\mathcal G}^{(d)} & &  f_2\in {\mathcal G}^{(d-1)} & &
\ldots & & f_{m}\in {\mathcal G}^{(d-(m-1))} & & {\mathcal
G}^{(d-m)},
\end{array}
\end{equation}
where $f_1,\ldots,f_m\in {\mathcal G}^{(d)}$ via the inclusions
$${\mathcal G}^{(d-(m-1))}\subset\ldots\subset {\mathcal G}^{(d-1)}
\subset {\mathcal G}^{(d)},$$ and where each $f_i$ is  transversal
to
$$(V^{(d-(i-1))},x_{i-1}) \to (V^{(d-i)},x_i),$$
for $i=1,\ldots,m$ (here we take $x_0=x$). Set $B={\mathcal
O}_{V^{(d)},x}/\langle f_1,\ldots,f_m\rangle$ and let $m_B$ be its maximal ideal.

\

\noindent {\bf Step 2.}  Under the assumptions of step 1, suppose  that,  in some  neighborhood of $x$,  {\bf an  arbitrary}   ${\mathcal G}^{(d)}$-admissible
local smooth projection to some $(d-m)$-dimensional smooth scheme, 
and an elimination algebra  are  given, say  
\begin{equation}
\label{transorder1}
\begin{array}{rrcl} \beta_{1_{d,d-m}}: &
(V^{(d)},x) & \longrightarrow &  (V^{(d-m)}_1, x_{m,1})\\
& {\mathcal G}^{(d)} &  & {\mathcal G}^{(d-m)}_1.
\end{array}
\end{equation}  Consider the local ring $( {\mathcal O}_{V^{(d-m)}_1, x_{m,1}}, m_{x_{m,1}})$.  
We   show that the natural map is an inclusion of 
  of local rings, say
\begin{equation}
\label{inclusionkey}
{\mathcal O}_{V^{(d-m)}_1, x_{m,1}} \subset  B={\mathcal O}_{V^{(d)},x}/\langle
f_1,\ldots,f_m\rangle,
\end{equation}
with the following properties:
\begin{enumerate}
\item[(i)]  The inclusion is finite and flat;
\item[(ii)] The ideal $m_{x_1}B$ is a reduction of $m_B$;
\item[(iii)] The  inclusion of Rees algebras
 $${\mathcal G}^{(d-m)}_1 \subset \overline{{\mathcal G}^{(d)}}(\subset B[W])$$
is an integral extension in $B[W]$ (here
$\overline{{\mathcal G}^{(d)}}$ denotes the image of ${\mathcal
G}^{(d)}$ in $B[W]$).
\end{enumerate}
Since (\ref{transorder1}) is an arbitrary admissible projection and $B$ has been fixed in step 1,
the first part of Theorem \ref{ordertau} will follow   from Proposition \ref{preorder}.

\

\noindent{\bf Idea of the proof of the second  part of Theorem \ref{ordertau}}

\

\noindent {\bf Step 1$^{\prime}$. }   Fix    $V^{(d)}$,  ${\mathcal G}^{(d)}\subset {\mathcal
O}_{V^{(d)}}[W]$,  $x\in \mbox{Sing
}{\mathcal G}^{(d)}$, $m\leq
\tau_{{\mathcal G},x}$ and an admissible projection as in step 1:
\begin{equation}
\label{projtesttrans}
(V^{(d)},x)\to
(V^{(d-m)},x_{m})
\end{equation}  together with  the factorization given  in (\ref{flag}), and the elements
$f_1,\ldots,f_m\in {\mathcal G}^{(d)}$  with the properties stated in step 1.

\

Let
$V^{(d)}\leftarrow V^{(d)^{\prime}}$ be a composition of permissible
monoidal transformations mapping $x^{\prime}$ to $x$. Then by
Theorem \ref{1theo} sequence (\ref{flag}) can be lifted to a
sequence of local admissible projections for the weighted transform of
${\mathcal G}^{(d)}$, ${\mathcal G}^{(d)^{\prime}}$, inducing a
commutative diagram of permissible transformations, local admissible
projections and elimination algebras,
\begin{equation}
\label{flagbu}
\begin{array}{ccccccccc}
(V^{(d)^{\prime}},x^{\prime}=x^{\prime}_0) & \to &
(V^{(d-1)^{\prime}},x_1^{\prime})& \to & \ldots & \to &
(V^{(d-(m-1)^{\prime})},x_{m-1}^{\prime}) & \to &
(V^{(d-m)^{\prime}},x_{m}^{\prime}) \\
{\mathcal G}^{(d)^{\prime}} & &  {\mathcal G}^{(d-1)^{\prime}} & &
\ldots & & {\mathcal G}^{(d-(m-1))^{\prime}} & & {\mathcal
G}^{(d-m)^{\prime}}
\\
\downarrow & & \downarrow & & \ldots & & \downarrow & & \downarrow \\
(V^{(d)},x=x_0) & \to & (V^{(d-1)},x_1)& \to & \ldots & \to &
(V^{(d-(m-1))},x_{m-1}) & \to & (V^{(d-m)},x_{m})\\
 {\mathcal G}^{(d)} & &   {\mathcal G}^{(d-1)}
& & \ldots & &
 {\mathcal G}^{(d-(m-1))} & & {\mathcal G}^{(d-m)}
\end{array}
\end{equation} 
where  $x_i^{\prime}$ maps to $x_i$ for $i=0,\ldots,m$.  Observe that the first 
vertical map is the blow up, while the other vertical arrows denote the induce 
blow ups in lower dimensions.  Notice that then    the
strict transforms of $f_1,\ldots,f_m$ in ${\mathcal
O}_{V^{(d)^{\prime}}}$, $f_1^{\prime},\ldots,f_m^{\prime}$, are in
${\mathcal G}^{(d)^{\prime}}$,  and  that moreover,
\begin{equation}
\label{lflag}
\begin{array}{ccccccc}
f_1^{\prime}\in {\mathcal G}^{(d)^{\prime}} & &  f_2^{\prime}\in
{\mathcal G}^{(d-1)^{\prime}} & & \ldots & & f_{m-1}^{\prime}\in
{\mathcal G}^{(d-(m-1))^{\prime}} .
\end{array}
\end{equation}
We will show that  each $f_i^{\prime}$ is  transversal
to
$$(V^{(d-(i-1))^{\prime}},x_{i-1}^{\prime}) \to (V^{(d-i)^{\prime}},x_i^{\prime}),$$
for $i=1,\ldots,m$ (here we take $x_0^{\prime}=x^{\prime}$).
 Set $B^{\prime}={\mathcal O}_{V^{(d)^{\prime}},x^{\prime}}/\langle
f_1^{\prime},\ldots,f_m^{\prime}\rangle$ and let $m_{B^{\prime}}$ be its maximal ideal.

\

\noindent {\bf Step 2$^{\prime}$.}   Under the assumptions of step 1$^{\prime}$, assume  that in some  neighborhood of $x$ {\bf an arbitrary}  ${\mathcal G}^{(d)}$-admissible
local smooth projection to some $(d-m)$-dimensional smooth scheme  is given (as in step 2):
\begin{equation}
\label{transorder1prime}
\begin{array}{rrcl} \beta_{1_{d,d-m}}: &
(V^{(d)},x) & \longrightarrow &  (V^{(d-m)}_1, x_{m,1})\\
& {\mathcal G}^{(d)} &   & {\mathcal G}^{(d-m)}_1,
\end{array}
\end{equation}
which by definition is constructed with a factorization as in 
Definition   \ref{def838}. Here (\ref{transorder1prime}) is probably 
different 
from  (\ref{projtesttrans}), so we use the subindex ''1´´ to distinguish (\ref{transorder1prime}) 
from (\ref{projtesttrans}).  Then the  composition of monoidal
permissible transformations $V^{(d)}\leftarrow V^{(d)^{\prime}}$ from step 1$^{\prime}$
induces a composition of permissible transformations and elimination
algebras  in (\ref{transorder1prime}),
$$
\xymatrix@R=0pc@C=0pc{ (V^{(d)},x)  & & & & & \quad (U\subset
V^{(d)^{\prime}},x^{\prime})\ \ar[lllll]\\
{\mathcal G}^{(d)}\ar[dddd]    & &  &  & & {\mathcal G}^{{(d)}^{\prime}}\ar[dddd] \\
\\
& & & \circlearrowleft \\
\\
(V^{(d-m)}_1,x_{m,1}) & &  &  & & (V^{(d-m)^{\prime}}_1,x_{m,1}^{\prime}) \ar[lllll]\\
{\mathcal G}^{(d-m)}_1  & &  &  & & {\mathcal G}^{(d-m)^{\prime}}_1.
}$$

Consider the local ring $( {\mathcal O}_{V^{(d-m)^{\prime}}_1, x_{m,1}^{\prime}}, m_{x_{m,1}^{\prime}})$.  Then we will show that the natural 
map  is an   inclusion of local rings, say 
$${\mathcal O}_{V^{(d-m)^{\prime}}_1, x_{m,1}^{\prime}} \subset  B^{\prime}={\mathcal O}_{V^{(d)^{\prime}},x^{\prime}}/\langle
f_1^{\prime},\ldots,f_m^{\prime}\rangle.$$
In addition we will see that this map is finite and flat, that the ideal $ m_{x_{m,1}^{\prime}} B$ is a reduction of $m_{B ^{\prime}}$ and that
 $${\mathcal G}^{(d-m)^{\prime}}_1 \subset 
\overline{{\mathcal G}^{(d)^{\prime}}}(\subset B^{\prime}[W])$$
is an integral extension in $B^{\prime}[W]$ (here
$\overline{{\mathcal G}^{(d)^{\prime}}}$ denotes the image of ${\mathcal
G}^{(d)^{\prime}}$ in $B^{\prime}[W]$).  Since  $B^{\prime}$ has been
 fixed in step  $1^{\prime} $, and $\beta_{1_{d,d-m}}$ is arbitrary, 
the second part of Theorem \ref{ordertau} will follow   from Proposition \ref{preorder}.

\

\noindent{\bf About steps  1 and 1$^{\prime}$}

\

The main difficulty in the proof of Theorem \ref{ordertau} is the accomplishment of steps 1 and $1^{\prime}$. More
precisely, and with the same notation as above, given a differential Rees algebra,  a suitable ring $B$ is constructed and fixed in
step 1. Then, in step 2, we have to show that for  {\bf any }   admissible projection there is an inclusion as in (\ref{inclusionkey})  that satisfies properties (i), (ii) and (iii).    Moreover, it is not immediate, either, that this situation  can be carried out after a finite sequence of permissible transformations. Thus,  the key for the proof is to find a suitable sequence of elements as in step 1 and step 1$^{\prime}$. Most part of this section
will be devoted to proving the existence of these particular sequences of elements for  a differential Rees  algebra.  First,
we have to introduce some definitions and  prove
auxiliary results:

\

- Given a Rees-algebra ${\mathcal G}$ and a simple closed point
$x\in \mbox{Sing }{\mathcal G}$, we introduce the notion of {\em
$\tau_{{\mathcal G},x}$-sequence} (see Definition \ref{tauseq}).
We will see that the existence of such sequences is guaranteed
when ${\mathcal G}$ is a differential Rees algebra.

\

- However it is not clear that $\tau_{{\mathcal G},x}$-sequences
behave well under permissible monoidal transformations, so they
are not suitable for proving   Theorem \ref{ordertau}. This
problem is overcome by introducing  {\em ${\mathcal G}$-nested
sequences} (see Definition \ref{nested}). Nested
sequences have some interesting  properties as
listed in \ref{PropNestedSeq}.  In particular, they behave well under  permissible
transformations.  It is worth pointing out that, while the notion of $\tau$-sequence is  intrinsic to    ${\mathcal G}$, the concept of nested sequence is relative to a particular smooth local smooth projection and a suitable factorization of it, as in  (\ref{flag}).

\

- The existence of nested sequences is not obvious:  in
Proposition \ref{614} and in Corollary \ref{PA} we show how to
construct ${\mathcal G}$-nested sequences starting from a $\tau_{{\mathcal
G},x}$-sequence.

\

Once the existence of nested sequences is established, the
proof Theorem \ref{ordertau} will  follow from their
properties. The proof of Theorem \ref{ordertau} wiil be addressed  in
\ref{prepreorder}.}
\end{Paragraph}

\

\noindent{\bf About $\tau$-sequences}

\

\begin{Definition}\label{tauseq}
{\rm Let ${\mathcal G}=\oplus_nI_nW^n$ be a Rees algebra in a
$d$-dimensional smooth scheme $V$ over a field $k$, let $x\in
\mbox{Sing }{\mathcal G}$ be a simple closed point, and let $k^{\prime}$
be the residue field at $x$. We will say that a set of homogeneous
elements $f_1W^{n_1},\ldots,f_sW^{n_s}\in {\mathcal G}$ is a {\em
$\tau_{{\mathcal G},x}$-sequence of length $s$} if
  for $j=1,\ldots,s$:
  \begin{itemize}
  \item[i.]   $n_j=p^{e_j}$;
  \item[ii.]
$\mbox{In}_xf_j, \in Gr_{{\mathcal O}_{V,x}}\simeq
k^{\prime}[Z_1,\ldots,Z_d]$
 is a $k^{\prime}$-linear combination of
 $Z_1^{p^{e_j}},\ldots,Z_d^{p^{e_j}}$ for some $e_j\in {\mathbb
 N}$;
 \item[iii.]  The class of $\mbox{In}_xf_j$
 is a regular element at the graded ring $Gr_{{\mathcal O}_{V,x}}/\langle \mbox{In}_xf_i: i\neq j\rangle$.
 \end{itemize}
 By definition, if  $f_1W^{n_1},\ldots,f_sW^{n_s}\in {\mathcal G}$ is a
$\tau_{{\mathcal G},x}$-sequence of length $s$,  then $s\leq
\tau_{{\mathcal G},x}$.}
 {\rm A $\tau_{{\mathcal G},x}$-sequence $f_1W^{n_1},\ldots,f_sW^{n_s}\in
{\mathcal G}$ is said to be a {\em
 maximal-$\tau_{{\mathcal G},x}$-sequence} if $s=\tau_{{\mathcal
 G},x}$. }
\end{Definition}

\begin{Paragraph} {\bf On the conditions of Definition
\ref{tauseq}.}\label{tausseq}  {\rm Let $f_1W^{n_1},\ldots,f_sW^{n_s} \in
{\mathcal G}$ be a $\tau_{{\mathcal G},x}$-sequence.   If char
$k=0$ then condition (ii)  says that $\mbox{In}_xf_1,\ldots,
\mbox{In}_xf_s\in Gr_{{\mathcal O}_{V,x}}$ are linear forms, while
 condition (iii) means that they are linearly independent. If  char $k=p>0$, then,
 up to a change of the base field, it can be assumed that $\mbox{In}_xf_j \in Gr_{{\mathcal O}_{V,x}}$
 is some $p^{e_j}$-th power of a linear form for $j=1,\ldots,s$. Condition (iii)
 indicates that these linear forms are independent (see \ref{tangentcone}). Notice that if $f_1W^{n_1},f_2W^{n_2}\ldots,f_sW^{n_s}$ is a $\tau$-sequence, then so is $(f_1)^p(W^{n_1})^p,f_2W^{n_2}\ldots,f_sW^{n_s}$. In particular it can always be assumed that
 $n_1=\ldots=n_s$. }
\end{Paragraph}

\begin{Remark}
{\rm When ${\mathcal G}^{(d)}$ is a differential Rees algebra, then there is
a maximal $\tau_{{\mathcal G}^{(d)},x}$-sequence at any  simple point
$x\in \mbox{Sing }{\mathcal G}^{(d)}$ (see \ref{tangentcone}). However
if
$$\begin{array}{ccc}
(V^{(d)},x) & \leftarrow & (V^{(d)^{\prime}},x^{\prime})\\
{\mathcal G}^{(d)} & & {\mathcal G}^{(d)^{\prime}} \end{array}$$ is a
permissible monoidal transformation  it is, in general, not true that the
strict transforms of a $\tau_{{\mathcal G}^{(d)},x}$-sequence form
a $\tau_{{\mathcal G}^{(d)^{\prime}},x^{\prime}}$-sequence. }
\end{Remark}

\

\noindent{\bf About nested sequences}

\

\begin{Definition}\label{nested}
{\rm Let  ${\mathcal G}^{(d)}$ be a Rees algebra,  and let $x\in \mbox{Sing
}{\mathcal G}^{(d)}$ be a simple point with $\tau_{{\mathcal G}^{(d)},x}\geq s$.
Suppose that  there is a ${\mathcal G}^{(d)}$-admissible projection
to some  $(d-s)$-dimensional smooth scheme  in a neighborhood of $x$,
$$(V^{(d)},x) \to  (V^{(d-s)},x_{s}) ,  $$
and a factorization into admissible projections
\begin{equation}
\label{nasociada}
\begin{array}{ccccccc}
(V^{(d)},x) & \stackrel{\beta_{d,d-1}}{\longrightarrow} &
   \ldots & \to &
(V^{(d-(s-1))},x_{s-1}) & \stackrel{\beta_{d-(s-1),d-s}}{\longrightarrow} & (V^{(d-s)},x_{s})\\
{\mathcal G}^{(d)} &     & \ldots &   &
 {\mathcal G}^{(d-(s-1))} & & {\mathcal
G}^{(d-s)}. 
\end{array}
\end{equation}
A set of homogeneous elements $f_1^{(d)}W^{n_1},
f_2^{(d-1)}W^{n_2},\ldots,f_s^{(d-(s-1))}W^{n_s}\in {\mathcal
G}^{(d)}$ is said to be a {\em ${\mathcal G}^{(d)}$-nested
sequence relative to sequence   (\ref{nasociada}}) if
$$
\begin{array}{ccccccc}
(V^{(d)},x=x_0) & \stackrel{\beta_{d,d-1}}{\longrightarrow} &
   \ldots & \to &
(V^{(d-(s-1))},x_{s-1}) & \stackrel{\beta_{d-(s-1),d-s}}{\longrightarrow} & (V^{(d-s)},x_{s})\\
f_1^{(d)}W^{n_1}\in {\mathcal G}^{(d)} &     & \ldots &   &
f_s^{(d-(s-1))}W^{n_s}\in {\mathcal G}^{(d-(s-1))} & & {\mathcal
G}^{(d-s)},
\end{array}
$$
and   $f_i^{(d-(i-1))}$ is transversal to $\beta_{d-(i-1),d-i}$
for $i=1,\ldots,s$ (see \ref{localproj} and \ref{Zariski} for the notion of transversality and its role in constructing admissible  local smooth projections). }
\end{Definition}

\begin{Paragraph} \label{PropNestedSeq}{\bf Some facts about nested
sequences.} {\rm Assume that $$f_1^{(d)}W^{n_1},
f_2^{(d-1)}W^{n_2}, \ldots,f_s^{(d-(s-1))}W^{n_s}\in {\mathcal
G}^{(d)}$$ is a ${\mathcal G}^{(d)}$-nested-sequence in a
neighborhood of $x$ as in Definition \ref{nested}, relative to a
sequence as in (\ref{nasociada}). Then:

\

\underline{1. Nested sequences define complete intersections.} In
other words, the quotient
$${\mathcal O}_{V^{(d)},x}/\langle
f_1^{(d)}, f_2^{(d-1)}, \ldots,f_s^{(d-(s-1))}\rangle$$ is
  a complete intersection.

\

To see this, notice that since $f_i^{(d-(i-1))}\in {\mathcal
G}^{(d-(i-1))}$ is transversal to $\beta_{d-(i-1),d-i}:
V^{(d-(i-1))}\to V^{(d-i)}$,
 for each $i=1,\ldots,s$  the local ring homomorphism
$${\mathcal
O}_{V^{(d-i)},x_{i}} \to {\mathcal
O}_{V^{(d-(i-1))},x_{i-1}}/\langle f_{i}^{(d-(i-1))}\rangle$$ is
finite and flat (up to an \'etale change of base, see
\ref{Zariski}). As a consequence,
$${\mathcal O}_{V^{(d)},x}/\langle f_1^{(d)}, f_2^{(d-1)},
\ldots,f_s^{(d-(s-1))}\rangle$$ is a finite free ${\mathcal
O}_{V^{(d-s)},x_{s}}$-module. Hence, the quotient is
Cohen-Macaulay, and moreover, a complete intersection.

\

\underline{2. Nested sequences and reductions.} If $m_{x_{s}}$
denotes the maximal ideal in ${\mathcal O}_{V^{(d-s)},x_{s}}$ then
$$m_{x_{s}}{\mathcal O}_{V^{(d)},x}/\langle f_1^{(d)},
f_2^{(d-1)}, \ldots,f_s^{(d-(s-1))}\rangle$$ is a reduction of the
maximal ideal in ${\mathcal O}_{V^{(d)},x}/\langle f_1^{(d)},
f_2^{(d-1)}, \ldots,f_s^{(d-(s-1))}\rangle$ (see \ref{Zariski},
specially the arguments involving formula (\ref{wfinite})).

\

\underline{3. Nested sequences lift to nested sequences after
  permissible monoidal transformations.} Let $V^{(d)}\leftarrow
V^{{(d)}^{\prime}}$ be a permissible monoidal transformation, let
${\mathcal G}^{{(d)}^{\prime}}$ be the  weighted transform of ${\mathcal
G}^{(d)}$ in $V^{{(d)}^{\prime}}$, and let $x^{\prime}_0\in
\mbox{Sing }{\mathcal G}^{{(d)}^{\prime}}$ be a closed point
dominating $x_0$. Then the strict transforms of $f_1^{(d)},
f_2^{{(d-1)}}, \ldots,f_s^{{(d-(s-1))}}$ in $V^{{(d)}^{\prime}}$,
which we denote by $f_1^{{(d)}^{\prime}}, f_2^{{(d-1)}^{\prime}},
\ldots,f_s^{{(d-(s-1))}^{\prime}}$,  form a ${\mathcal
G}^{{(d)}^{\prime}}$-nested sequence relative to the transform of
sequence (\ref{nasociada})   (see Theorem \ref{1theo} and its proof). Therefore the quotient
$${\mathcal O}_{V^{{(d)}^{\prime}},x_0^{\prime}}/\langle
f_1^{{(d)}^{\prime}}, f_2^{{(d-1)}^{\prime}},
\ldots,f_s^{{(d-(s-1))}^{\prime}}\rangle
$$
 defines a complete intersection, and hence it is flat over ${\mathcal O}_{V^{{(d-s)}^{\prime}},x_{0,s}^{\prime}}$.
 If ${\mathcal I}(E)\subset {\mathcal O}_{V^{{(d)}^{\prime}}}$ denotes
 the ideal sheaf of the exceptional divisor, then the strict transform of the ideal
$$\langle f_1^{(d)}, f_2^{(d-1)}, \ldots,f_s^{(d-(s-1))}\rangle \subset 
{\mathcal O}_{V^{{(d)}},x}$$ in ${\mathcal
O}_{V^{{(d)}^{\prime}},x_0^{\prime}}$ is given by the increasing union of 
colon ideals 
$$\cup_{n\geq 0}\left(\langle f_1^{(d)}, f_2^{(d-1)}, 
\ldots,f_s^{(d-(s-1))}\rangle
{\mathcal O}_{V^{{(d)}^{\prime}},x^{\prime}_0}: {\mathcal I}(E)^n
\right)$$
which in this case equals to
$$\langle f_1^{{(d)}^{\prime}}, f_2^{{(d-1)}^{\prime}},
\ldots,f_s^{{(d-(s-1))}^{\prime}}\rangle \subset {\mathcal
O}_{V^{{(d)}^{\prime}},x^{\prime}_0}$$ since $${\mathcal
O}_{V^{{(d)}^{\prime}},x_0^{\prime}}/\langle f_1^{{(d)}^{\prime}},
f_2^{{(d-1)}^{\prime}}, \ldots,f_s^{{(d-(s-1))}^{\prime}}\rangle
$$ is flat over ${\mathcal
O}_{V^{{(d-m)}^{\prime}},x_{0,m}^{\prime}}$. Here ${\mathcal I}(E)\subset 
{\mathcal O}_{V^{(d)^{\prime}}}$ is obtained by lifting 
the ideal of the exceptional divisor from $V^{(d-m)}\leftarrow V^{(d-m)^{\prime}}$. 

\

\underline{4. Nested sequences and finite extensions.} There is a
diagram
$$
\xymatrix@R=0pc@C=0pc{    & &  &  &  & {\mathcal
G}^{(d)} \subset   {\mathcal O}_{V^{(d)},x}[W] \ar[dddd]^{\gamma^*_0}\\
\\
& & &   \\
\\
{\mathcal G}^{(d-s)}\subset {\mathcal O}_{V^{(d-s)} ,x_s}[W]
\ar[rrrrr]^{\gamma^*_s}  & & &  & &   {\mathcal
O}_{V^{(d)},x}/\langle f_1^{(d)},\ldots,f_s^{(d-(s-1))}\rangle[W]
}$$
 where
$\gamma^*_{s}$ is a finite map, $\gamma^*_0$ is the natural
surjection, and
$$\gamma^*_{s}({\mathcal
G}^{(d-s)})\subset \gamma^*_0({\mathcal G}^{(d)})$$ is a finite
extension of graded algebras in ${\mathcal O}_{V^{(d)},x}/\langle
f_1^{(d)},\ldots,f_s^{(d-(s-1))}\rangle[W]$.

\

To see this, notice that by  property (1), for each
$i=1,\ldots,s$, the map
\begin{equation}
\label{compo} {\mathcal O}_{V^{(d-i)},x_i} {\longrightarrow}
{\mathcal O}_{V^{(d)},x}/\langle f_1^{(d)},
\ldots,f_i^{(d-(i-1))}\rangle
\end{equation} is  finite, and       factorizes as
$${\mathcal O}_{V^{(d-i)},x_i}\stackrel{\eta^*_{d-i,d-(i-1)}}{\longrightarrow} {\mathcal
O}_{V^{(d-(i-1))},x_{i-1}}/\langle f_i^{(d-(i-1))}\rangle
\stackrel{\delta^*_{d-(i-1)}}{\longrightarrow} {\mathcal
O}_{V^{(d)},x}/\langle f_1^{(d)}, \ldots,f_i^{(d-(i-1))}\rangle.$$

\

For $i=1,\ldots,s$ consider the diagram

$$
\xymatrix@R=0pc@C=0pc{   & & & & & & &  & & & {\mathcal
G}^{(d-(i-1))} \subset  {\mathcal O}_{V^{(d-(i-1))},x_{i-1}}[W] \ar[dddd]^{\alpha_{d-(i-1)}^*}\\
 & & & & & & & & & &  \\
 & & & & & & & & & &  \\
  & & & & & & & & & &  \\
{\mathcal G}^{(d-i)}\subset {\mathcal O}_{V^{(d-i)} ,x_i}[W]
\ar[rrrrrrrrrr]^{\eta^*_{d-i,d-(i-1)}}  & & &  & & & & & & &
{\mathcal O}_{V^{(d-(i-1))},x_{i-1}}/\langle
f_i^{(d-(i-1))}\rangle[W] \ar[dddd]^{\delta^*_{d-(i-1)}}\\
 & & & & & & & & & &  \\
 & & & & & & & & & &  \\
  & & & & & & & & & &  \\
    & & & & & & & & & &  \\
& &  &  & & & & & & &  {\mathcal O}_{V^{(d)},x}/\langle f_1^{(d)},
\ldots,f_i^{(d-(i-1))}\rangle[W], }$$ where $\eta^*_{d-i,d-(i-1)}$
is   finite,  $\alpha_{d-(i-1)}^*$ is surjective and
$f_i^{(d-(i-1))}\in {\mathcal G}^{(d-(i-1))}$ is transversal to
$\beta_{d-(i-1),d-i}: (V^{(d-(i-1))},x_{i-1})\to
(V^{(d-i)},x_i)$.

 \

According to  Theorem \cite[4.11]{hpositive} (see also \ref{eliminationalg}),
the inclusion
$$\eta^*_{d-i,d-(i-1)}({\mathcal G}^{(d-i)})\subset \alpha_{d-(i-1)}^*({\mathcal G}^{(d-(i-1))})$$ is
a finite extension of  graded algebras. Therefore,
\begin{equation}
\label{inte}
\delta^{*}_{d-(i-1)}\left(\eta^*_{d-i,d-(i-1)}({\mathcal
G}^{(d-i)})\right)\subset
\delta^{*}_{d-(i-1)}\left(\alpha_{d-(i-1)}^*({\mathcal
G}^{(d-(i-1))})\right), \end{equation} is a finite extension of
graded algebras in ${\mathcal O}_{V^{(d)},x}/\langle f_1^{(d)},
\ldots,f_i^{(d-(i-1))}\rangle[W]$.

\

Since the map
$${\mathcal O}_{V^{(d)},x}/\langle f_1^{(d)},
\ldots,f_i^{(d-(i-1))}\rangle \to {\mathcal O}_{V^{(d)},x}/\langle
f_1^{(d)}, \ldots,f_s^{(d-(s-1))}\rangle,$$ is surjective, the
resulting maps
$$\begin{array}{c}
\gamma^*_{i-1}: {\mathcal O}_{V^{(d-(i-1))},x_{i-1}}\to {\mathcal
O}_{V^{(d)},x}/\langle f_1^{(d)},\ldots,f_s^{(d-(s-1))}\rangle, \\
\gamma^*_i: {\mathcal O}_{V^{(d-i)},x_i}\to {\mathcal
O}_{V^{(d)},x}/\langle f_1^{(d)},\ldots,f_s^{(d-(s-1))}\rangle
\end{array}$$ are a composition of finite and surjective maps. Hence
$$\gamma^*_{i}({\mathcal G}^{(d-i)})\subset
\gamma^*_{i-1}({\mathcal G}^{(d-(i-1))}) \left(\subset {\mathcal
O}_{V^{(d)},x}/\langle f_1^{(d)},
\ldots,f_s^{(d-(s-1))}\rangle[W]\right)$$ is a finite extension.

\

Using  an inductive argument we conclude that there is a sequence of
 finite inclusions of Rees algebras
$$\gamma^*_{s}({\mathcal G}^{(d-s)})\subset \ldots \subset
\gamma^*_1({\mathcal G}^{(d-1)})\subset \gamma^*_0({\mathcal
G}^{(d)}) \left(\subset {\mathcal O}_{V^{(d)},x}/\langle f_1^{(d)},
\ldots,f_s^{(d-(s-1))}\rangle[W]\right).$$  }
\end{Paragraph}

\

\noindent {\bf The existence of nested sequences for
differential Rees algebras}

\

In the following we consider a differential Rees algebra on a smooth scheme
 together with an admissible local smooth projection. Our goal  is to show that there is a $\tau$-sequence which is, in addition, a nested sequence for this given admissible projection. This will be settled in Corollary \ref{PA}. In this procedure we will start from an arbitrary  $\tau$-sequence at a singular point.

\begin{Proposition}\label{614} Let ${\mathcal G}=\oplus_nI_nW^n$ be a differential Rees algebra on a
$d$-dimensional smooth scheme $V^{(d)}$   over a field $k$. Let
$x\in \mbox{Sing }{\mathcal G}\subset V^{(d)}$ be a simple closed
point and let $$f_1^{(d)}W^{n_1},\ldots, f_s^{(d)}W^{n_s}$$ be a
maximal $\tau_{{\mathcal G},x}$-sequence of length $\tau_{{\mathcal G},x}=s\geq 2$.  Fix
a  ${\mathcal G}$-admissible  projection in a neighborhood of $x$, $\beta_{d,d-1}:V^{(d)} \to V^{(d-1)}$  (see Definition \ref{def83}). Then:

{\rm A)}  For some index $i$, $1\leq i \leq s$,  $f_i^{(d)}$ is transversal to
$\beta_{d,d-1}$.

{\rm B)}   Set $i=1$ as in {\rm (A) } (after reordering the sequence if needed) and
construct an elimination algebra ${\mathcal R}_{{\mathcal
G},\beta_{d,d-1}}$ as described in \ref{eliminationalg}. Then:

\begin{enumerate}

\item[ i.] There is a $\tau_{{\mathcal
R}_{{\mathcal G},\beta_{d,d-1}}}$-sequence of length $(s-1)$,
$f_2^{(d-1)}W^{l_2},\ldots,f_s^{(d-1)}W^{l_s}\in {\mathcal
R}_{{\mathcal G},\beta_{d,d-1}}$.
\item[ii.] The previous $\tau_{{\mathcal
R}_{{\mathcal G},\beta_{d,d-1}}}$-sequence   can be
constructed so that  $$\langle
f_2^{(d-1)},\ldots,f_s^{(d-1)}\rangle \subset \langle
f_1^{(d)},\ldots, f_s^{(d)} \rangle \subset
{\mathcal
O}_{V^{(d)},x}$$ via the  inclusion ${\mathcal
O}_{V^{(d-1)},x_1}\subset {\mathcal O}_{V^{(d)},x}$.
\end{enumerate}
\end{Proposition}

\

\noindent{\em Proof:}
A) Our hypotheses are that  $ \mathcal G $ is a differential Rees algebra
and that $f_1^{(d)}W^{n_1},\ldots, f_s^{(d)}W^{n_s}$ is a
maximal $\tau_{{\mathcal G},x}$-sequence of length $s\geq 2$ at
$x\in \mbox{Sing }{\mathcal G}\subset V^{(d)}$. Recall that each
$n_i= p^{e_i}$, that each $f_i$ has order $p^{e_i}$ at
$ {\mathcal
O}_{V^{(d)},x}$, and that $\mbox{In}_x(f_i) \in \mbox{Gr}_{m_x}({\mathcal
O}_{V^{(d)},x})$ is homogeneous of degree $p^{e_i}$ (and moreover a
$p^{e_i}$-th power of a linear form). In addition,  the tangent cone defined by $ \mathcal G $ at $\mathbb{T}_{V^{(d)},x}=\mbox{Spec }(\mbox{Gr}_{m_x}({\mathcal
O}_{V^{(d)},x})) $ is the closed set defined by the ideal $\langle \mbox{In}_x(f_1), \dots , \mbox{In}_x(f_s) \rangle$. Since  we are assuming that the conditions in Definition \ref{def83} hold, there must be an index $i$ for which $f_i$ is transversal to
$\beta_{d,d-1}:V^{(d)} \to V^{(d-1)}$.

\

B) By (A) we can assume that $f_1$ is transversal to $\beta_{d,d-1}$ (here a reordering of the $\tau$-sequence may be needed).
Suppose that $\mbox{In}_x(f_1)=Y_1^{p^{e_1}}$ for some linear form $Y_1 \in \mbox{Gr}_{m_x}({\mathcal
O}_{V^{(d)},x})$. Let $\{z_2, \dots , z_{d}\}$ be a regular system of parameters in  $ {\mathcal
O}_{V^{(d-1)},x_1}$. Choose $y_1$ to be  an element of order one at
$ {\mathcal
O}_{V^{(d)},x}$, so that $\mbox{In}_x(y_1)=Y_1 \in \mbox{Gr}_{m_x}({\mathcal
O}_{V^{(d)},x})$. Then  $\{y_1, z_2, \dots , z_{d}\}$ is a regular system of parameters in $ {\mathcal
O}_{V^{(d)},x}$, and $\mbox{Gr}_{m_x}({\mathcal
O}_{V^{(d)},x})$ is a polynomial ring in variables $\{Y_1, Z_2, \dots , Z_{d}\}$, where $Z_i=\mbox{In}_x(z_i)$  for $i=2,\ldots,d$.

\

 Recall that the $\tau$-sequence $f_1^{(d)}W^{n_1},\ldots, f_s^{(d)}W^{n_s}$ is defined with $n_i=p^{e_i}$, which can be chosen so that
$e_1=e_2= \cdots =e_s=e$ (see (\ref{tausseq})).  Let $k'$ denote the residue field of ${\mathcal
O}_{V^{(d)},x}$. Then:
\begin{enumerate}
\item[-] $\mbox{Gr}_{m_x}({\mathcal
O}_{V^{(d)},x})=k'[Y_1, Z_2, \dots , Z_{d}]$;

\item[-]  $\mbox{In}_x(f^{(d)}_1)=Y^{e}_1$;

\item[-]   For $j=2,\ldots,s$,  $\mbox{In}_x(f^{(d)}_j)=\lambda_j  Y^{p^e}_1+(L_j)^{p^e}$, for some $\lambda_j \in k'$ and $L_j $ a linear form in $k'[ Z_2, \dots , Z_{d}]$;

\item[-]  The linear forms $\{ÊL_j , j=2, \dots s\}$ are independent in
$\mbox{Gr}_{m_{x_1}}({\mathcal
O}_{V^{(d-1)},x_1})=k'[ Z_2, \dots , Z_{d}]$.
\end{enumerate}

Assume, for simplicity, that $k'=k$ (by finite extension of base field),
set ${f^{\prime}}^{(d)}_1=f^{(d)}_1$ and let ${f^{\prime}}^{(d)}_j=\lambda_j f^{(d)}_1-f^{(d)}_j$ for $j=2,\ldots,s$.
Notice  that $\{ {f^{\prime}}^{(d)}_1W^{p^e}, \dots, {f^{\prime}}^{(d)}_sW^{p^e}\}$ is a
$\tau_{{\mathcal G},x}$-sequence, that $\langle {f^{\prime}}^{(d)}_1, \dots {f^{\prime}}^{(d)}_s\rangle=\langle f^{(d)}_1, \dots f^{(d)}_s\rangle$. A regular system of parameters
$\{ v_2, \dots , v_{d}\}$
can be chosen in
${\mathcal
O}_{V^{(d-1)},x_1}$ so that:

\

a) The set $\{y_1, v_2, \dots ,v_{d}\}$ is a regular system of parameters in ${\mathcal
O}_{V^{(d)},x}$. In particular  $\mbox{Gr}_{m_x}({\mathcal
O}_{V^{(d)},x})=k'[Y_1, V_2, \dots , V_{d}]$,  and $V_i=\mbox{In}_x(v_i)$ for $i=2,\ldots,d$.

\

b) $\mbox{In}_x(f^{\prime}_1)=Y_1^{p^{e}}$, and
$\mbox{In}_x(f^{\prime}_j)=V_{j} ^{p^{e}} $,  for $=2,\ldots,s$.

\

Under these assumptions, part B i) of the Proposition was  proven in \cite[5.12]{hpositive}.  We  briefly sketch
the argument here: The setting now is that  $f_1^{(d)}W^{n_1},\ldots, f_s^{(d)}W^{n_s}$ is  a
$\tau_{{\mathcal G},x}$-sequence of length $s$,  all $n_i=p^e$, and  there is a regular system of parameters
$\{y_2^{},\ldots,y_d^{}\}\subset
{\mathcal O}_{V^{(d-1)},x_1}$ which extends to  $\{y_1^{},y_2\ldots,y_d^{}\}\subset
{\mathcal O}_{V^{(d)},x}$
and  $\mbox{In}_x
f_i^{(d)}=\mbox{In}_xy_i^{p^{e}}\in \mbox{Gr}_{m_x}({{\mathcal
O}_{V^{(d)},x}})$ for  $i=1,\ldots,s$.

\

We assume that $f_1$ is a monic polynomial of degree $p^e$ in $y_1$ and coefficients in ${\mathcal
O}_{V^{(d-1)},x_1}$, so ${\mathcal O}_{V^{(d)},x}/\langle f_1^{(d)}\rangle$ is a free ${\mathcal
O}_{V^{(d-1)},x_1}$-module of rank $p^e$.

\

For each $i=2,\ldots,s$, let $\overline{f_i^{(d)}}$ be the image of
$f_i^{(d)}$ in ${\mathcal O}_{V^{(d)},x}/\langle f_1^{(d)}\rangle$.
 Multiplying   by $\overline{f_i^{(d)}}$ induces a map
 of free
${\mathcal O}_{V^{(d-1)},x_1}$-modules:
$$\Gamma_{\overline{f_i^{(d)}}}:{\mathcal O}_{V^{(d)},x}/\langle f_1^{(d)} \rangle \to {\mathcal O}_{V^{(d)},x}/\langle
f_1^{(d)} \rangle,$$
and similarly,  multiplying  by $\overline{f_i^{(d)}}W^{n_i}$ defines a map of free ${\mathcal O}_{V^{(d-1)},x_1}[W]$-modules:
$$\Gamma_{\overline{f_i^{(d)}}W^{n_i}}:{\mathcal O}_{V^{(d)},x}/\langle f_1^{(d)} \rangle [W]\to {\mathcal O}_{V^{(d)},x}/\langle
f_1^{(d)} \rangle [W].$$
Let $p_i(t)$ be the characteristic polynomial of 
$\Gamma_{\overline{f_i^{(d)}}W^{n_i}}$. Let $g_iW^{l_i}\in {\mathcal O}_{V^{(d-1)},x_1}[W]$ be the
determinant (i.e.,  $g_iW^{l_i}=p_i(0)$). Note that $g_i$ is the determinant of $\Gamma_{\overline{f_i^{(d)}}}:{\mathcal O}_{V^{(d)},x}/\langle f_1^{(d)} \rangle \to {\mathcal O}_{V^{(d)},x}/\langle
f_1^{(d)} \rangle$.

\

Under these conditions, it can be shown that  $g_i$ has order $p^{2e}$ in $ {\mathcal O}_{V^{(d-1)},x_1}$, that $g_iW^{p^{2e}}\in {\mathcal
R}_{{\mathcal G},\beta_{d,d-1}}$,
and that $\mbox{In}_{x_1}g_i=\mbox{In}_{x_1}(y_i^{})^{p^{2e}}$ for
$i=2,\dots ,s$, where, as indicated before,
$\{y_2,\ldots,y_d\}$ is a regular system of parameters  in ${\mathcal
O}_{V^{(d-1)},x_1}$ (see
\cite[5.12]{hpositive} for more  details on  this  proof).

\

To prove  B ii), observe that  the composition
$${\mathcal O}_{V^{(d)},x}/\langle f_1^{(d)} \rangle \stackrel{\Gamma_{\overline{f_i^{(d)}}}}{\longrightarrow}
{\mathcal O}_{V^{(d)},x}/\langle f_1^{(d)} \rangle \to {\mathcal
O}_{V^{(d)},x}/\langle f_1^{(d)},f_i^{(d)}\rangle,$$ (where the
last row is just the natural quotient morphism),  maps the image of
$\Gamma_{\overline{f_i^{(d)}}}$ to zero. Since  $g_i$ is the determinant of the first,  any sufficiently high
power of $g_i$ is zero in ${\mathcal
O}_{V^{(d)},x}/\langle f_1^{(d)},f_i^{(d)}\rangle$.
In particular, for $e^{\prime}$ large enough, $g_i^{p^{e^{\prime}}}\in \langle f_1^{(d)},f_i^{(d)}\rangle$, for $i=2, \dots s$.

\

Finally define $f_i^{(d-1)}=g_i^{p^{e^{\prime}}}$  and
$l_i=p^{2e+e^{\prime}}$ for $i=2, \dots s$. So:
$$f_2^{(d-1)}W^{l_2},\ldots,f_s^{(d-1)}W^{l_s}\in {\mathcal
R}_{{\mathcal G},\beta_{d,d-1}}$$
is a $\tau$-sequence, and
$$\langle
f_2^{(d-1)},\ldots,f_s^{(d-1)}\rangle \subset \langle
f_1^{(d)},\ldots, f_s^{(d)} \rangle \subset
{\mathcal
O}_{V^{(d)},x}.\qed$$

\begin{Corollary}\label{PA} Let ${\mathcal G}^{(d)}=\oplus_nI_nW^n$ be a differential Rees algebra over a
  $d$-dimensional
 smooth scheme $V^{(d)}$  over a field $k$.
 Let $x\in \mbox{Sing }{\mathcal G}^{(d)}$ be a simple point, and let
$f_1^{(d)}W^{n_1},\ldots, f_s^{(d)}W^{n_s}\in {\mathcal G}^{(d)}$
be a maximal $\tau_{{\mathcal G}^{(d)},x}$-sequence of length $s=\tau_{{\mathcal G}^{(d)},x}$.
Consider a ${\mathcal G}^{(d)}$-admissible local smooth projection to a
$(d-s)$-dimensional scheme,
\begin{equation}
\label{admissiblepreorder}
\begin{array}{rrcl} \beta_{d,d-s}:
& (V^{(d)},x) &
\longrightarrow &  (V^{(d-s)},x_s)\\
 & {\mathcal G}^{(d)} & & {\mathcal G}^{(d-s)},
\end{array}
\end{equation} and a
factorization of (\ref{admissiblepreorder}) as  a sequence of
${\mathcal G}^{(d-i)}$-admissible projections (see Definition \ref{def83},
\begin{equation}
\label{inductive}
\begin{array}{ccccccc} (V^{(d)},x) &
\stackrel{\beta_{d,d-1}}{\longrightarrow} &  \ldots &
\longrightarrow & (V^{(d-(s-1))},x_{s-1})&
\stackrel{\beta_{d-(s-1),d-s}}{\longrightarrow}&  (V^{(d-s)},x_s)\\
{\mathcal G}^{(d)} &    & &   & {\mathcal G}^{(d-(s-1))} &   &
{\mathcal G}^{(d-s)}.
\end{array}
\end{equation}
Then, after reordering $f_1^{(d)}W^{n_1},\ldots, f_s^{(d)}W^{n_s}$,
if needed,  for each  $i=1,\ldots,s-1$:
\begin{enumerate}
\item[i.] There is a $\tau_{\Gdi,x_i}$-sequence of length $(s-i)$,
 $$f_{i+1}^{(d-i)}W^{l_{i,i+1}},\ldots,f_{s}^{(d-i)}W^{l_{i,s}}\in \Gdi\subset {\mathcal
 O}_{V^{(d-i)},x_i}[W].$$
\item[ii.]  For $f_{i+1}^{(d-i)},\ldots, f_s^{(d-1)}$ as in (i) 
there is an inclusion of ideals
$$\langle f_{i+1}^{(d-i)},\ldots,f_s^{(d-i)}\rangle\subset \langle f_{i}^{(d-(i-1))},f_{i+1}^{(d-(i-1))},
\ldots, f_{s}^{(d-(i-1))} \rangle\subset {\mathcal
 O}_{V^{(d-(i-1))},x_{i-1}}$$ via the
inclusion ${\mathcal
 O}_{V^{(d-i)},x_i}\to  {\mathcal
 O}_{V^{(d-(i-1))},x_{i-1}}$;
 \item[iii.] Moreover there is a  ${\mathcal G}^{(d)}$-nested-sequence of length
$s$ in a neighborhood of $x$,
$$f_1^{(d)}W^{t_1},\ldots,f_s^{(d-(s-1))}W^{t_s}$$ relative  to
sequence (\ref{inductive}), that is also a $\tau_{{\mathcal
G}^{(d)},x}$-sequence,  and with
$$\langle f_1^{(d)},\ldots,f_s^{(d-(s-1))}\rangle \subset \langle f_1^{(d)},\ldots,f_s^{(d)}
\rangle \subset {\mathcal O}_{V^{(d)},x}$$ 
where $f^{(d-i)}$ is viewed  in ${\mathcal O}_{V^{(d)},x}$  via the inclusion 
${\mathcal O}_{V^{(d-i)},x_i}\to {\mathcal O}_{V^{(d)},x}$ for
$i=1,\ldots,s-1$.

\end{enumerate}

\end{Corollary}

\

\noindent{\em Proof:} After relabelling
$f_1^{(d)},\ldots,f_s^{(d)}\in {\mathcal G}^{(d)}$, we may assume
that $f_1^{(d)}$ is transversal to $\beta_{d,d-1}$. Now  the
corollary  follows from Proposition \ref{614} and an inductive
argument since the elimination algebra of a differential Rees algebra
is also a differential Rees algebra. \qed

\

\noindent{\bf Proof of Theorem \ref{ordertau}}

\begin{Paragraph} \label{prepreorder} {\rm By Proposition \ref{614}, (B) (ii), 
it is enough the proof the theorem in the case in which $m=\tau$.  The first part of the theorem will be proven in
two steps.}

\

\noindent {\rm {\bf Step 1.} Assume  that ${\mathcal G}^{(d)}$ is
a differential Rees algebra. Then, by Corollary
\ref{PA} (iii), we can assume that  there is a $\tau_{{\mathcal
G}^{(d)},x}$-sequence of length $m$,
\begin{equation}
\label{thesequence} f_1^{(d)}W^{n_1},\ldots,f_m^{(d)}W^{n_m}\in
{\mathcal G}^{(d)}
\end{equation}
 that is also ${\mathcal G}^{(d)}$-nested   relative  to
some sequence of ${\mathcal G}^{(d)}$-local admissible projections
\begin{equation}
\label{theprojections}
\begin{array}{ccccccc}
(V^{(d)},x)  & \to & (V^{(d-1)}_0,x_{0,1}) & \to  &  \ldots & \to & (V^{(d-m)}_0,x_{0,m})\\
{\mathcal G}^{(d)} & & {\mathcal G}^{(d-1)}_0  &  & \ldots & &
{\mathcal G}^{(d-m)}_0.
\end{array}
\end{equation}
 Hence the map of local rings
$${\mathcal O}_{V^{(d-m)}_0,x_{0,m}}\to {\mathcal O}_{V^{(d)},x}/\langle
f_1^{(d)},\ldots,f_m^{(d)}\rangle$$ is finite and flat and
therefore the quotient
$${\mathcal O}_{V^{(d)},x}/\langle
f_1^{(d)},\ldots,f_m^{(d)}\rangle$$ is a $(d-m)$-dimensional
Cohen-Macaulay ring (see \ref{PropNestedSeq}).  Let     $B={\mathcal O}_{V^{(d)},x}/\langle
f_1^{(d)},\ldots,f_m^{(d)}\rangle$, and denote by $m_B$ its maximal ideal.

\

\noindent {\bf Step 2.} Suppose that we are  given  {\bf an arbitrary}
${\mathcal G}^{(d)}$-admissible projection to some $(d-m)$-smooth
scheme (\ref{def838}), and an elimination algebra
\begin{equation}
\label{cualquiera}
\begin{array}{rrcl}
\beta_{{d,d-m}}: & V^{(d)}& \longrightarrow &  V^{(d-m)} \\
& ({\mathcal G}^{(d)},x) &   & ({\mathcal G}^{(d-m)},x_{m})
\end{array}
\end{equation}
then:

(a) Notice that there  is a natural map $${\mathcal O}_{V^{(d-m)},x_{m}} \to B=
{\mathcal O}_{V^{(d)},x}/\langle
f_1^{(d)},\ldots,f_m^{(d)}\rangle.$$

(b)   We claim that      the images of ${\mathcal G}^{(d-m)}$
and ${\mathcal G}^{(d)}$ in $B[W]={\mathcal O}_{V^{(d)},x}/\langle
f_1^{(d)},\ldots,f_m^{(d)}\rangle [W]$ have the same integral closure. Since the local
admissible  projection (\ref{cualquiera}) is arbitrary,  and the
sequence $f_1^{(d)},\ldots,f_m^{(d)}$ (and hence $B={\mathcal
O}_{V^{(d)},x}/\langle f_1^{(d)},\ldots,f_m^{(d)}\rangle$) is fixed,
 the first part of Theorem \ref{ordertau} follows from  Proposition \ref{preorder}.

\

The claim  in  (b) can be accomplished by finding  a nested
sequence relative to some factorization of (\ref{cualquiera}) into
locally admissible projections. This nested sequence will be
constructed using the $\tau$-sequence found in Step 1.

\

So, consider any ${\mathcal G}^{(d)}$-admissible projection   to
some $(d-m)$-smooth scheme in a neighborhood of $x$, with its corresponding elimination
algebra,
\begin{equation}
\label{anyadmissible}
\begin{array}{rrcl}
\beta_{{d,d-m}}: & (V^{(d)},x) & \longrightarrow &  (V^{(d-m)},x_{m}) \\
& {\mathcal G}^{(d)}&   & {\mathcal G}^{(d-m)},
\end{array}
\end{equation}
and  construct a ${\mathcal G}^{(d)}$-nested sequence   relative to
some factorization of (\ref{anyadmissible}), say 
\begin{equation}
\label{consSequence}
f_{1}^{(d)}W^{l_1},\ldots,f_{m}^{(d-(m-1))}W^{l_m}\in {\mathcal
G}^{(d)},
\end{equation}  using  the
$\tau_{{\mathcal G},x}$-sequence from step 1,
$$f_1^{(d)}W^{n_1},\ldots,f_m^{(d)}W^{n_m},$$ and following the strategy 
of  Corollary
\ref{PA} (ii) (here some relabeling may be needed). Notice that by construction,
\begin{equation}
\label{contenido} \langle
f_{1}^{(d)},\ldots,f_{m}^{(d-(m-1))}\rangle\subset \langle
f_1^{(d)},\ldots,f_m^{(d)}\rangle\subset {\mathcal O}_{V^{(d)},x}.
\end{equation}

\

According  to properties (1) and (2) in \ref{PropNestedSeq},
$${\mathcal O}_{V^{(d-m)},x_{m}}\to {\mathcal
O}_{V^{(d)},x}/\langle
f_{1}^{(d)},\ldots,f_{m}^{(d-(m-1))}\rangle$$ is a finite flat
local map  of $(d-m)$-dimensional Cohen-Macaulay rings (here an
\'etale change of base maybe needed) and by (\ref{contenido}) the
map
$${\mathcal
O}_{V^{(d)},x}/\langle f_{1}^{(d)},\ldots,f_{m}^{(d-(m-1))}\rangle
\to {\mathcal O}_{V^{(d)},x}/\langle
f_1^{(d)},\ldots,f_m^{(d)}\rangle$$ is surjective. Therefore
$${\mathcal O}_{V^{(d-m)},x_{m}}\to B={\mathcal
O}_{V^{(d)},x}/\langle f_1^{(d)},\ldots,f_m^{(d)}\rangle$$ is
finite. Since $B$ is a complete intersection, both  
${\mathcal O}_{V^{(d-m)},x_{m}}$ and $B$ 
have the same dimension, so it follows that the previous morphism is flat. 

\

Let $m_{x_m}$ be the maximal ideal in ${\mathcal
O}_{V^{(d-m)},x_{m}}$, and let $m_x$ be the maximal ideal in
${\mathcal O}_{V^{(d)},x}$. Then by \ref{PropNestedSeq} (2),
$$m_{x_m}{\mathcal O}_{V^{(d)},x}/\langle
f_{1}^{(d)},\ldots,f_{m}^{(d-(m-1))}\rangle$$ is a reduction 
of the maximal ideal of ${\mathcal O}_{V^{(d)},x}/\langle
f_{1}^{(d)},\ldots,f_{m}^{(d-(m-1))}\rangle$, say 
$$m_x\left({\mathcal O}_{V^{(d)},x}/\langle
f_{1}^{(d)},\ldots,f_{m}^{(d-(m-1))}\rangle\right).$$ Thus by
(\ref{contenido})
$$m_{x_m}\left({\mathcal
O}_{V^{(d)},x}/\langle f_1^{(d)},\ldots,f_m^{(d)}\rangle\right)$$ is a
reduction of the maximal ideal $m_B$ of $B$,
$$m_x \left({\mathcal O}_{V^{(d)},x}/\langle f_1^{(d)},\ldots,f_m^{(d)}\rangle\right).$$

\

Finally,  consider the  diagram:
$$
\xymatrix@R=0pc@C=0pc{    & &  &  &  & {\mathcal
G}^{(d)} \subset   {\mathcal O}_{V^{(d)},x}[W] \ar[dddd]^{\gamma^*}\\
\\
& & &   \\
\\
{\mathcal G}^{(d-m)}\subset {\mathcal O}_{V^{(d-m)} ,x_m}[W]
\ar[rrrrr]^{\gamma^*_m}  & & &  & &   {\mathcal
O}_{V^{(d)},x}/\langle f_1^{(d)},\ldots,f_m^{(d-(m-1))}\rangle[W]
}.$$

According to  \ref{PropNestedSeq}  (4) $\gamma^*_{m}({\mathcal
G}^{(d-m)})\subset \gamma^*({\mathcal G}^{(d)})$ is a finite
extension of graded algebras. Thus by (\ref{contenido}), their
images in $B$ are also a finite extension.

\

Since (\ref{anyadmissible}) was an arbitrary ${\mathcal
G}^{(d-m)}$-admissible projection, by Proposition \ref{preorder}
  $\mbox{ord}_{x_m}{\mathcal G}^{(d-m)}$ is independent on the choice of the
  projection. This proves the first part of the Theorem.

\

The second part of the Theorem will be accomplished in two steps.

\

\noindent{\bf Step 1$^{\prime}$.}  Fix the $\tau_{{\mathcal
G}^{(d)},x}$-sequence found in  (\ref{thesequence}) which is also a
${\mathcal G}^{(d)}$-nested sequence relative  to
(\ref{theprojections}). Now suppose that $V^{(d)} \leftarrow
V^{(d)^{\prime}}$ is a composition of permissible monoidal
transformations, and that  $x^{\prime}\in \mbox{Sing }{\mathcal
G}^{(d)^{\prime}}$   is a closed point dominating $x$. By Theorem
\ref{1theo}  there is a commutative diagram of permissible monoidal transformations and admissible
projections. Observe that the strict transforms of
$f_1^{(d)},\ldots,f_m^{(d)}\in {\mathcal G}^{(d)}$ in
$V^{(d)^{\prime}}$, say
$f_1^{(d)^{\prime}},\ldots,f_m^{(d)^{\prime}}$, form a ${\mathcal
G}^{(d)^{\prime}}$-nested sequence in a neigborhood of
$x^{\prime}$ (see \ref{PropNestedSeq}). Let  $B^{\prime}={\mathcal O}_{V^{(d)^{\prime}},x^{\prime}}/\langle
f_1^{(d)^{\prime}},\ldots,f_m^{(d)^{\prime}}\rangle $ and let $m_{B^{\prime}}$ denote its maximal ideal.

\

\noindent{\bf Step 2$^{\prime}$.}  Fix  {\bf an arbitrary}  ${\mathcal G}^{(d)}$-admissible local smooth projection as in
(\ref{anyadmissible}) and consider the composition of permissible monoidal transformations from step 1$^{\prime}$,  $V^{(d)} \leftarrow
V^{(d)^{\prime}}$. Again  by Theorem
\ref{1theo} there is a commutative diagram of elimination algebras
and admissible projections:
\begin{equation}
\label{explosion}
 \xymatrix@R=0pc@C=0pc{ (V^{(d)},x)  & & & & &
\quad (U\subset
V^{(d)^{\prime}}, x^{\prime})\ \ar[lllll]\\
{\mathcal G}^{(d)}\ar[dddd]_{\beta_{d,d-m}}    & &  &  & & {\mathcal G}^{{(d)}^{\prime}}\ar[dddd]^{\beta_{d,d-m}^{\prime}} \\
\\
& & & \circlearrowleft \\
\\
(V^{(d-m)},x_m) & &  &  & & (V^{(d-m)^{\prime}},x_m^{\prime}) \ar[lllll]\\
{\mathcal G}^{(d-m)}  & &  &  & & {\mathcal G}^{(d-m)^{\prime}}. }
\end{equation}
Here the vertical maps correspond to the arbitrary projection, and the horizontal arrows correspond to the blow-ups. 
\

Recall, as observed in step  1$^{\prime}$, that the strict transforms in $V^{(d)^{\prime}}$ of the
${\mathcal G}$-nested sequence
$f_{1}^{(d)},\ldots,f_{m}^{(d-(m-1))}$  given in
(\ref{consSequence}), say
$f_{1}^{{(d)}^{\prime}},\ldots,f_{m}^{(d-(m-1))^{\prime}}$, form a
${\mathcal G}^{(d)^{\prime}}$ nested sequence in a neighborhood of
$x^{\prime}$ relative to some factorization of the ${\mathcal
G}^{(d)^{\prime}}$-admissible local smooth projection (\ref{explosion}). Also
  there is an inclusion ideals
$$\langle f_{1}^{{(d)}^{\prime}},\ldots,f_{m}^{{(d-(m-1))}^{\prime}}\rangle\subset
\langle f_1^{(d)^{\prime}},\ldots,f_m^{(d)^{\prime}}\rangle,$$
(see \ref{PropNestedSeq} and (\ref{contenido})). Now the proof follows from a similar
argument  as the one given in step 2, using  $(B^{\prime}, m_{B^{\prime}})$ instead of $(B,m_B)$ (see also Theorem \ref{1theo} and its proof).
\qed }
\end{Paragraph}

\

\section{The non-simple case}\label{nonsimple}

Let ${\mathcal G}^{(d)}$ be a Rees algebra and  let  $x\in
\mbox{Sing }{\mathcal G}^{(d)}$ be a simple point with
$\tau_{{\mathcal G}^{(d)},x}\geq m$. Under the assumptions of
Theorem \ref{ordertau} there are well defined
upper semi-continuous functions in a neighborhood of $x$:
$$\begin{array}{rrcl}
\mbox{ord}^{(d-i)}:&
\mbox{Sing }{\mathcal G}^{(d)} & \longrightarrow & {\mathbb Q}\\
  & z & \to & \mbox{ord}_{z}^{(d-i)}{\mathcal G}^{(d)}= \mbox{ord}_{z_i}{\mathcal
  G}^{(d-i)},
\end{array}$$
where $z_i:=\beta_{d,d-i}(z)$,  and  $\beta_{d,d-i}: V^{(d)}\to
V^{(d-i)}$ is a ${\mathcal G}^{(d-i)}$-admissible local smooth projection
on to some $(d-i)$-dimensional smooth scheme $V^{(d-i)}$,  for
$i=0,1,\ldots,m$ (\ref{def838}). Since $\tau_{{\mathcal G}^{(d)},x}\geq m$,
$$\mbox{ord}_{x}^{(d)}{\mathcal G}^{(d)}=\ldots=\mbox{ord}_{x}^{(d-(m-1))}{\mathcal G}^{(d)}=1.$$

\

In the following  we  denote by $\mbox{max-ord}^{(d-i)}$  the
maximum value of the function $\mbox{ord}^{(d-i)}$ and we will use
${\mbox{\underline{Max}-ord}}^{(d-i)}$ to denote the   the closed
set
$$\{z\in \mbox{Sing }{\mathcal G}^{(d)}:
\mbox{ord}_z^{(d-i)}{\mathcal
G}^{(d)}=\mbox{max-ord}^{(d-i)}{\mathcal G}^{(d)}\}.$$

\

Now suppose that $\tau_{{\mathcal G}^{(d)},x}\geq m$. Fix a
${\mathcal G}^{(d)}$-admissible projection to some
$(d-m)$-dimensional smooth scheme, $\beta_{(d-m)}: V^{(d)}\to
V^{(d-m)}$, and let $x_m:=\beta_{(d-m)}(x)$. If $m=d$, then
$\mbox{Sing }{\mathcal G}^{(d)}=\{x\}$ in a neighborhood of $x$,
and a resolution of ${\mathcal G}^{(d)}$ is achieved by blowing up
this point (in the sense of  \ref{ragadprs3}). This is a particular case that will be discussed in Remark \ref{rkcdime}. 

\

On the other hand, if $m<d$, and if  $x_{m}\in \mbox{Sing }{\mathcal
G}^{(d-m)}$ is not a simple point, then it
would  be interesting to, somehow,``enlarge $\tau_{{\mathcal
G}^{(d)},x}$". This will be done by
  extending  ${\mathcal G}^{(d)}$ to a larger Rees
algebra $\widetilde{\mathcal G}^{(d)}\supset {\mathcal G}^{(d)}$
so that $\mbox{Sing }\widetilde{\mathcal
G}^{(d)}={\mbox{\underline{Max} ord}}^{(d-m)}{\mathcal G}^{(d)}$
and $\tau_{\widetilde{\mathcal G}^{(d)},x}\geq m+1$ (i.e.,
$x_{m}\in \mbox{Sing }\widetilde{\mathcal G}^{(d-m)}$ will be a
simple point). In this case, a stratification of $\mbox{Sing
}\widetilde{\mathcal G}^{(d)}$ will induce a stratification of
${\mbox{\underline{Max} ord}}^{(d-m)}\subset \mbox{Sing }{\mathcal
G}^{(d)}$ by descending induction on the value of  $\tau$.

\

The purpose of this section is to show how this enlargement can be
done in full generality, the main result  is the formulation of  Theorems
\ref{twisting} and \ref{twistingtransform} in which the main
properties of $\widetilde{\mathcal G}^{(d)}$ are discussed. It is
at this point where the notions of weak equivalence introduced in
\ref{Defweak} appear in  full strength. In fact, these two theorems
show that $\widetilde{\mathcal G}^{(d)}$ can be chosen so as to be
well defined up to integral closure (see also Remark \ref{kr212}).

\

We begin by recalling the notion of twisted algebras introduced in
\cite{EV}.

\begin{Definition}\label{mdeftwisting}
{\rm  Let ${\mathcal G}=\bigoplus_{n\geq 0}J_nW^n$ be a Rees
algebra on a smooth  $d$-dimensional scheme $V$ and let $\omega$
be a positive rational number. The {\em twisted algebra}
${\mathcal G}(\omega)$ is defined  as
$${\mathcal G}(\omega)=\bigoplus_{n\geq
0}J_{\frac{n}{\omega}}W^n$$ where it is assumed that
$J_{\frac{n}{\omega}}=0$ if $\frac{n}{\omega}$ is not an integer.}
\end{Definition}


\begin{Proposition}\label{propertiestwisting}
The twisted algebra of Definition \ref{mdeftwisting}   satisfies
the following properties:
\begin{enumerate}
\item[(i)] If ${\mathcal G}={\mathcal G}_{(J,b)}$ and if $w$ is a
positive rational number with $b \omega\in {\mathbb Z}$ then
${\mathcal G}(\omega)={\mathcal G}_{(J,\omega b)}$.
\item[(ii)] If ${\mathcal G}_1$ and ${\mathcal G}_2$ have the same
integral closure, then so do ${\mathcal G}_1(\omega)$ and
${\mathcal G}_2(\omega)$.
\item[(iii)] ${\mathcal G}(\omega)$ is a Rees algebra and
$\omega\cdot \mbox{ord}_x{\mathcal G}(\omega)=\mbox{ ord}_x{\mathcal
G}$. In particular if $\omega=\mbox{ord}_x{\mathcal G}$ then
$\mbox{ord}_x{\mathcal G}(\omega)=1$.
\item[(iv)] If $\omega=\mbox{max-ord }^{(d)}{\mathcal G}$, then
${\mathcal G}(\omega)$ is simple and $\mbox{Sing }{\mathcal
G}(\omega)=\mbox{\underline{Max}-ord }^{(d)}{\mathcal G}$.
\end{enumerate}
\end{Proposition}
For the proof we  refer the reader to   \cite[Propositions 6.4, 6.5, 6.7 and
Corollary 6.7]{EV}.
\begin{Remark}
\label{setwisting} {\rm Let ${\mathcal G}$ be a Rees algebra, and
let $\omega$ be the maximum of the function
$$\begin{array}{rrcl}
\mbox{ord}: & \mbox{Sing }{\mathcal G} & \longrightarrow & {\mathbb
Q}.
\end{array}$$
If  $x\in \mbox{Sing }{\mathcal G}$, then  ${\mathcal G}={\mathcal
G}(\omega)$ at $x$ if and only if $x$ is a simple point. If $x$ is
not a simple point for ${\mathcal G}$, then $\tau_{\mathcal
G,x}=0$, but  then  $x\in \mbox{Sing }{\mathcal G}(\omega)$ is a
simple point of ${\mathcal G}(\omega)$, so in particular $\tau_{\mathcal G(\omega),x}\geq
1$.}
\end{Remark}

\begin{Definition} \label{definetwisting} {\rm Let ${\mathcal G}$ be a Rees algebra on a smooth scheme $V$, and
let $x\in \mbox{Sing }{\mathcal G}$. Let
$\omega=\mbox{ord}_x{\mathcal G}$. If $x$ is not a simple point,
i.e., if $\omega>1$, then, define  $\widetilde{\mathcal
G}:={\mathbb D}\mbox{iff}{ (\mathcal G}(\omega))$ (\ref{extensiondif}).}
\end{Definition}

\begin{Remark}
{\rm Using  the same notation as in the previous definition, notice
that $\mbox{Sing }{\widetilde{\mathcal G}}=\{z\in V:
\mbox{ord}_z{\mathcal G}=\omega\}$ in a neighborhood of $x$.}
\end{Remark}

\begin{Definition}
{\em Given two algebras over \(V\), for instance  \(\mathcal{G}_{1}\) and
\(\mathcal{G}_{2}\), set $ \mathcal{G}_{1}\odot\mathcal{G}_{2} $
as the smallest subalgebra of \(\calo_V[W]\) containing both (as
in (\ref{eq876})). Let $U$ be an affine open set in $V$. If the
restriction of $\mathcal{G}_1$ to $U$ is
$\calo_V(U)[f_1W^{n_1},\dots ,f_sW^{n_s}]$, and that of
$\mathcal{G}_2$ is $\calo_V(U)[f_{s+1}W^{n_{s+1}},\dots
,f_tW^{n_t}]$, then the restriction of
$\mathcal{G}_{1}\odot\mathcal{G}_{2}$ to $U$ is
$$\calo_V(U)[f_1W^{n_1},\dots ,f_sW^{n_s}, f_{s+1}W^{n_{s+1}},\dots
,f_tW^{n_t}].$$

One can check that:
\begin{enumerate}
    \item \(\Sing(\mathcal{G}_{1}\odot\mathcal{G}_{2})=
    \Sing(\mathcal{G}_{1})\cap\Sing(\mathcal{G}_{2})\).  In
    particular, if \( V\leftarrow V^{\prime}\)   is a permissible transformation for
    \(\mathcal{G}_{1}\odot\mathcal{G}_{2}\), then it is also a
    permissible transformation
    for \(\mathcal{G}_{1}\) and for \( \mathcal{G}_{2}\).

    \item If  \( V\leftarrow V^{\prime}\) is a permissible transformation
    for \(\mathcal{G}_{1}\odot\mathcal{G}_{2}\), and if
    \((\mathcal{G}_{1}\odot\mathcal{G}_{2})'\), \(\mathcal{G}'_1\),
    and \(\mathcal{G}'_2\) denote their transforms in  \(V'\), then:
    \begin{equation*}
    (\mathcal{G}_{1}\odot\mathcal{G}_{2})'=
    \mathcal{G}'_1\odot\mathcal{G}'_2.
    \end{equation*}
\end{enumerate}
}
\end{Definition}

\begin{Theorem}
\label{deftwisting} Let ${\mathcal G}=\bigoplus_nI_nW^n$ be a
differential Rees algebra defined on a $d$-dimensional smooth  scheme
$V^{(d)}$   over a perfect field $k$, and  let $x\in \mbox{Sing }
{\mathcal G}$ be a simple point (i.e., $\tau_{{\mathcal G},x}\geq
1$). Assume that $x$ is not contained in any component of  codimension one
 of $\mbox{Sing }{\mathcal G}$,  and assume that locally at $x$, 
$$\omega:=\mbox{max-ord}^{(d-1)} {\mathcal G}>1.$$ Fix  two
${\mathcal G}$-admissible local smooth projections to some $(d-1)$-dimensional
smooth schemes, and consider the corresponding elimination
algebras and    twisted algebras as in Definition
\ref{definetwisting},
$$\begin{array}{lcclclccl}
\beta_1: & (V^{(d)},x) & \longrightarrow & (V^{(d-1)}_1,x_{1,1}) &
\hspace{0.5cm} &
\beta_2: & (V^{(d)},x) & \longrightarrow & (V^{(d-1)}_2,x_{1,2}) \\
 & {\mathcal G} & &  {\mathcal R}_{{\mathcal G},\beta_1} \subset
{\mathcal R}_{{\mathcal G},\beta_1}(\omega) & &  & {\mathcal G} &
& {\mathcal R}_{{\mathcal G},\beta_2}\subset {\mathcal
R}_{{\mathcal G},\beta_2}(\omega).
\end{array}$$
Then
$$\widetilde{\mathcal G}_1={\mathcal G}\odot {\mathcal R}_{{\mathcal
G},\beta_1}(\omega) \ \ \mbox{ and } \ \  \widetilde{\mathcal
G}_2={\mathcal G}\odot {\mathcal R}_{{\mathcal
G},\beta_2}(\omega)$$ are weakly equivalent  (see Definition
\ref{Defweak}).
\end{Theorem}

\begin{Remark} \label{kr212}  Observe that $\mbox{Sing }{\widetilde{\mathcal G}}_1=\mbox{Sing }{\widetilde{\mathcal G}}_2=\underline{Max }{\mbox ord}^{(d-1)}{\mathcal G}$. Since
$\widetilde{\mathcal G}_1$ and $ \widetilde{\mathcal G}_2$ are weakly equivalent, Theorem \ref{examplesweak} asserts  that there is a
canonical differential Rees algebra    $\widetilde{\mathcal
G}\subset {\mathcal O}_V[W]$,   such that
$$\mbox{\underline{Max}-ord}^{(d-1)}{{\mathcal G}}=\mbox{Sing }
\widetilde{\mathcal G} \  \ \mbox{ and } \ \
\tau_{\widetilde{\mathcal G},x}\geq 2$$ for all $z\in
\mbox{\underline{Max}-ord}^{(d-1)}{\mathcal G}$ in some
neighborhood of $x$.

\end{Remark}

\noindent{\em Proof of Theorem \ref{deftwisting}.} We   have to
show that $\widetilde{\mathcal G}_1$ and $\widetilde{\mathcal
G}_2$ define the same singular locus under any permissible
morphism in the sense of Definition
\ref{permissibletransformation}. This is straightforward   for
permissible morphisms as in \ref{specialsmooth}  (i) and
(ii), so we are only left with the case of permissible monoidal
transformations.

\

Let  $V   \leftarrow  V^{\prime}$ be a permissible monoidal
transformation with center $$Y\subset
\mbox{\underline{Max}-ord}^{(d-1)}{\mathcal G}=\mbox{Sing
}\widetilde{\mathcal G}_1=\mbox{Sing }\widetilde{\mathcal G}_2,$$
and let $x^{\prime}\in V^{\prime}$ be a closed point that dominates
$x$. Then by Theorem \ref{1theo} there is
 a commutative diagram of algebras and elimination algebras in a suitable open set of $V^{\prime}$:
$$\xymatrix@R=0pc@C=0pc{ (V,x)  & & & & & (V^{\prime}, x^{\prime})  \ar[lllll]\\
{\mathcal G}\subset \widetilde{\mathcal G}_i \ar[dddd]^{\beta_i} &
&  & & & {\mathcal
G}^{\prime} \subset \widetilde{\mathcal G}_i^{\prime} \ar[dddd]^{\beta_i^{\prime}} \\
\\
& & & \\
\\
(V^{(d-1)}_i,x_{1,i})  & &  &  & & (V^{{(d-1)}^{\prime}}_i,x_{1,i}^{\prime})  \ar[lllll]\\
{\mathcal R}_{{\mathcal G},\beta_i}\subset  {\mathcal
R}_{{\mathcal G},\beta_i}(\omega)   & &  &  & & {\mathcal
R}_{{\mathcal G},\beta_i}^{\prime}\subset {\mathcal R}_{{\mathcal
G},\beta_i}(\omega)^{\prime}  }
$$
for $i=1,2$.

\

Now, on the one hand  by Theorem  \ref{ordertau},
$$\beta_i^{\prime}(\mbox{\underline{Max}-{w-ord}}^{(d-1)}{\mathcal
G}^{\prime}) =\mbox{\underline{Max}-w-ord}{\mathcal R}_{{\mathcal
G}^{\prime},\beta_i}=\underline{\mbox{Max}} \mbox{ ord} {\mathcal R}_{{\mathcal G},\beta_i}(\omega)^{\prime}$$ in a neighborhood of $x_{1,i}^{\prime}$
(see \ref{lasat1} for the definition of $\mbox{w-ord}^{(d-i)}$). Here ${\mathcal R}_{{\mathcal G},\beta_i}(\omega)^{\prime}$ is the transform of the simple  Rees algebra ${\mathcal R}_{{\mathcal G},\beta_i}(\omega)$. 
Moreover,
$$\widetilde{\mathcal G}_i^{\prime}=\left({\mathcal
G}\odot{\mathcal R}_{{\mathcal
G},\beta_i}(\omega)\right)^{\prime}= {\mathcal G}^{\prime}\odot
{\mathcal R}_{{\mathcal G},\beta_i}(\omega)^{\prime},$$ and
$\mbox{\underline{Max}-w-ord}^{(d-1)}{\mathcal
G}^{\prime}=\mbox{Sing }\widetilde{\mathcal G}_i$ for $i=1,2$.
Therefore
$$\mbox{Sing }\widetilde{\mathcal G}_1^{\prime}=\mbox{Sing }\widetilde{\mathcal
G}_2^{\prime}$$ in a neighborhood of $x^{\prime}$. \qed

\begin{Theorem}\label{twisting}
Let ${\mathcal G}^{(d)}$ be a differential Rees algebra on a  smooth 
$d$-dimensional  scheme   over a perfect field $k$, let
$x\in \mbox{Sing }{\mathcal G}^{(d)}$ be a simple point and assume
that $\tau_{{\mathcal G}^{(d)},x}=m\geq 1$.  Then, in an open neighborhood of $x$,     there is a
differential Rees algebra $\widetilde{\mathcal G}^{(d)}$ containing ${\mathcal
G}^{(d)}$ with the following properties:
\begin{enumerate}
\item[(i)] $\tau_{\widetilde{\mathcal G}^{(d)},z}\geq
m+1$ for $z\in \mbox{Sing }\widetilde{\mathcal G}^{(d)}$.
\item[(ii)] There is an equality of  closed sets $${\mbox{Sing }}{\widetilde{\mathcal G}}^{(d)}=
\mbox{\underline{Max} ord}^{(d-m)}{\mathcal G}^{(d)}.$$
\item[(iii)] The differential Rees algebra $\widetilde{\mathcal G}^{(d)}$  is unique
up to weak equivalence. Furthermore, this differential Rees algebra is
unique up to integral closures of algebras (see Remark
\ref{kr212}).
\end{enumerate}
\end{Theorem}

\noindent{\em Proof:}  Consider a ${\mathcal G}^{(d)}$-admissible
local smooth projection to some $(d-m)$-smooth dimensional scheme,
$$(V^{(d)},x)\to (V^{(d-m)},x_m),$$ and a factorization as in the
  diagram,
$$
\xymatrix@R=0pc@C=0pc{
{\mathcal G}^{(d)}  \subset  {\mathcal O}_{V^{(d)}}[W]\\
  \\
  \\
  \\
  {\mathcal G}^{(d-(m-1))} \subset
  {\mathcal O}_{V^{(d-(m-1))}}[W] \ar[uuuu]_{\beta^*_{d,d-(m-1)}}\\
  \\
   \\
    \\
{\mathcal G}^{(d-m)} \subset {\mathcal O}_{V^{(d-m)}}[W]
\ar[uuuu]_{\beta^*_{d-(m-1),d-m}}.}$$

By Theorem \ref{deftwisting}   there is a differential Rees algebra,
$\widetilde{\mathcal G}^{(d-(m-1))}$, containing ${\mathcal
G}^{(d-(m-1))}$ with the following properties:
\begin{enumerate}
\item[(a)] Its $\tau$-invariant  at $x_{m-1}=\beta_{d,d-(m-1)}(x)$ is larger than that of ${\mathcal
G}^{(d-(m-1))}$, i.e.,
\begin{equation}\label{comparetau}\tau_{\widetilde{\mathcal
G}^{(d-(m-1))},x_{m-1}}\geq \tau_{{\mathcal
G}^{(d-(m-1))},x_{m-1}}+1, \end{equation} and therefore
$x_{m-1}\in \mbox{Sing }\widetilde{\mathcal G}^{(d-(m-1))}$ is a
simple point.
\item[(b)] By construction  $\widetilde{\mathcal G}^{(d-(m-1))}$ is
unique up to weak equivalence.
\item[(c)] There is an equality of closed sets (using the
identification between singular loci),
\begin{equation}
\label{singups} \mbox{Sing } \widetilde{\mathcal
G}^{(d-(m-1))}=\mbox{\underline{Max}-ord } {\mathcal G}^{(d-m)}=
\mbox{\underline{Max}-ord}^{(d-m)}{\mathcal G}^{(d)}
\end{equation}
\end{enumerate}
in a neighborhood of $x$.

Now set
$$\widetilde{\mathcal G}^{(d)}={\mathcal
G}^{(d)}\odot\beta_{d,d-(m-1)}^*\widetilde{\mathcal
G}^{(d-(m-1))}.$$

Next we check that this algebra satisfies the properties stated in
the Theorem:
\begin{enumerate}
\item[(i)] By construction
$$\tau_{\widetilde{\mathcal G}^{(d)},x}=\tau_{\widetilde{\mathcal
G}^{(d-(m-1))},x_{m-1}}+ (m-1)  \geq   \tau_{{\mathcal
G}^{(d-(m-1))},x_{m-1}}+1 + (m -1) \geq m +1.
$$
\item[(ii)] Notice that via the natural identification of the
singular loci and from  (\ref{singups}), locally, in a
neighborhood of $x$,
$$\mbox{Sing }\widetilde{\mathcal G}^{(d)}=\mbox{Sing }\widetilde{\mathcal
G}^{(d-(m-1))}= \mbox{Sing }\widetilde{\mathcal G}^{(d-m)}=
\mbox{\underline{Max}-ord } {\mathcal
G}^{(d-m)}=\mbox{\underline{Max}-ord}^{(d-m)}{\mathcal G}^{(d)}.$$

\item[(iii)] The argument to show this part is similar to the
proof of Theorem \ref{deftwisting} since by Theorem \ref{ordertau}
$\mbox{ord}^{(d-m)}$  does not depend on the choice of the
${\mathcal G}^{(d)}$-admissible projections and therefore
$\mbox{\underline{Max}-ord}^{(d-m)}{\mathcal
G}^{(d)}=\mbox{\underline{Max}-ord } {\mathcal G}^{(d-m)}$ is well defined in a
neighborhood of $x$. \qed

\end{enumerate}

\

Finally we state a similar result for permissible transforms of
differential Rees algebras.

\begin{Theorem}\label{twistingtransform}
Let ${\mathcal G}^{(d)}$ be a differential Rees algebra on a
$d$-dimensional smooth scheme of finite type over a perfect  field $k$, let
$x\in \mbox{Sing }{\mathcal G}^{(d)}$ be a simple point and assume
that $\tau_{{\mathcal G}^{(d)},x}=m\geq 1$. Let
$$V^{(d)}\leftarrow V^{(d)^{\prime}}$$
be a composition of permissible monoidal transformations, let
${\mathcal G}^{(d)^{\prime}}$ be the weighted transform of ${\mathcal
G}^{(d)}$ and let $x^{\prime}\in \mbox{Sing }{\mathcal
G}^{(d)^{\prime}}$ be a closed point that dominates  $x$. Then there
exists an algebra, $\widetilde{\mathcal G}^{(d)^{\prime}}$
containing ${\mathcal G}^{(d)^{\prime}}$ with the following
properties:
\begin{enumerate}
\item[(i)] $\tau_{\widetilde{\mathcal G}^{(d)^{\prime}},x^{\prime}}\geq
m+1.$
\item[(ii)] Locally at $x^{\prime}$ there is an equality of  closed sets
$${\mbox{Sing }}{\widetilde{\mathcal G}}^{(d)^{\prime}}=
\mbox{\underline{Max} w-ord}^{(d-m)}{\mathcal G}^{(d)^{\prime}}.$$
\item[(iii)] The algebra $\widetilde{\mathcal G}^{(d)^{\prime}}$  is unique
up to integral closure of algebras.
\end{enumerate}
\end{Theorem}

\noindent {\em Proof:} Under the assumptions of the theorem, consider a local 
${\mathcal G}$-admissible projection to some $(d-m)$-dimensional 
smooth scheme, $\beta_{d,d-m}:V^{(d)}\to V^{(d-m)}$. Set $\omega:=\mbox{max ord}^{(d-m)}{\mathcal G}^{(d)}=\mbox{max ord }{\mathcal G}^{(d-m)}$, and define $ \widetilde{\mathcal G}^{(d)}={\mathcal G}^{(d)}\odot 
\beta_{d,d-m}^*{\mathcal  G}^{(d-m)}(\omega)$ as in the proof of   Theorem
\ref{twisting}. After considering the sequence of blow ups, $V^{(d)}\leftarrow V^{(d)^{\prime}}$, 
the proof proceeds in the same manner as that of Theorem 
\ref{twisting}, since,  as in that case, by Theorem \ref{ordertau},
the functions $\mbox{ord}^{(d-m)}$ are well defined for ${\mathcal
G}^{(d)^{\prime}}$ in a neighborhood of $x^{\prime}$. \qed

\section{Stratification of the singular locus by smooth
strata}\label{stratificationtheorem} The purpose of this section
is to prove the following theorem:
\begin{Theorem}
\label{stratification} Let ${\mathcal G}^{(d)}$ be a differential
 algebra on a smooth $d$-dimensional scheme  $V^{(d)}$
  over a perfect field $k$. Let ${\mathbb Q}^*={\mathbb Q}\cup
\{\infty\}$ and let
$$I_d=\underbrace{{\mathbb Q}^*\times {\mathbb Q}^* \times \ldots \times{\mathbb
Q}^*}_{d-\mbox{\tiny{times}}}$$ ordered lexicographically.  Then
there is an upper semi-continuous function,
$$\gamma_{{\mathcal G}^{(d)}}: \mbox{Sing }{\mathcal G}^{(d)} \to I_d$$
such that:
\begin{enumerate}
\item[(i)] The level sets of $\gamma_{{\mathcal G}^{(d)}}$ stratify
$\mbox{Sing }{\mathcal G}^{(d)}$ in smooth  locally closed strata.
\item[(ii)] If $k$ is a field of characteristic zero then
$\gamma_{{\mathcal G}^{(d)}}$ coincides with the resolution
function used for resolution of singularities in characteristic
zero.
\end{enumerate}
\end{Theorem}

The proof of the Theorem is given in  \ref{function}. First we
need the following lemma.

\begin{Lemma}\label{codimensionone} Let ${\mathcal G}$ be a Rees
algebra on a smooth $d$-dimensional scheme $V$  over a field $k$,
and assume that $x\in \mbox{Sing }{\mathcal G}$ is a simple point.
If  $x$ is contained in a component of
  codimension one, $Y$,  of
$\mbox{Sing }{\mathcal G}$, then $Y$ is smooth.
\end{Lemma}

\noindent{\em Proof:} We may assume that, up to integral closure
${\mathcal G}=\bigoplus_nJ^nW^{k_0n}$ for some sheaf of ideals
$J\subset {\mathcal O}_V$. The hypothesis of the lemma asserts
that $Y\subset \mbox{Sing }{\mathcal G}.$ Since ${\mathcal G}$
is simple at $x$, the hypothesis means that locally in a suitable
neighborhood of $x$, $Y\subset V(J)$ where $V(J)$
 denotes the closed set determined by $J$. By restricting to a smaller
 neighborhood $U$ of $x$ if needed, we may assume that ${\mathcal
I}(Y)=\langle f\rangle$, for some reduced element $f\in {\mathcal O}_V(U)$ and
that $J=\langle f^s\rangle$ for some positive integer $s$. Since
$x$ is a
 simple point and $Y\cap U\subset \mbox{Sing }{\mathcal G}\cap U$, $f$ has to be smooth at $x$. \qed

\begin{Remark}\label{rkcdime}
{\rm As a consequence of the previous lemma an inductive argument, using elimination algebras, shows that  if $x\in\mbox{Sing }{\mathcal G}$ is contained in a
component of codimension $\tau_{{\mathcal G},x}$, then the
component is smooth in a neighborhood of $x$. One can check that 
in this case a resolution of ${\mathcal G}$ is achieved by blowing up such components. 
}
\end{Remark}

\begin{Paragraph} \label{function} {\it Proof of Theorem
\ref{stratification}.} {\rm  We start by defining the function
$\gamma_{{\mathcal G}^{(d)}}$.  Let $x\in \mbox{Sing }{\mathcal
G}^{(d)}$. To associate a value to $\gamma_{{\mathcal G}^{(d)}}$
at $x$, we will argue by induction on the dimension of $V^{(d)}$.

\

Suppose   that $V^{(1)}$ is a one-dimensional smooth scheme over a
field $k$, that ${\mathcal G}^{(1)}$ is a non-zero differential
algebra and that $x\in \mbox{Sing }{\mathcal G}^{(1)}$ is a closed
point. Set
$$\gamma_{{\mathcal G}^{(1)}}(x)=(\mbox{ord}^{(1)}_x{\mathcal
G}^{(1)}).$$

\

Suppose that the function $\gamma$ can be defined   for any
differential Rees algebra ${\mathcal G}^{(n)}$ on a $n$-dimensional
smooth scheme
 over a perfect field $k$, $V^{(n)}$,  with $n<d$. We will show that then the function can be defined  for any
 non-zero differential Rees algebra ${\mathcal G}^{(d)}$ on  a $d$-dimensional smooth scheme
 over $k$, $V^{(d)}$.

\

\noindent $\bullet$ First assume that $x\in \mbox{Sing }{\mathcal
G}^{(d)}$ is a simple closed point. We now distinguish between two
cases:

\

\noindent {\bf Case 1.} If $x$ is  contained in a component of
codimension one, $Y$, of  $\mbox{Sing }{\mathcal
G}^{(d)}$, then   set
$$\gamma_{{\mathcal G}^{(d)}}(x)=
(\mbox{ord}^{(d)}_x{\mathcal
G}^{(d)},\underbrace{\infty,\ldots,\infty}_{d-1-\mbox{\tiny{times}}})=
(1,\underbrace{\infty,\ldots,\infty}_{d-1-\mbox{\tiny{times}}}).$$
Note that by Lemma \ref{codimensionone} the closed subscheme  $Y$
is smooth locally at $x$.

\

\noindent {\bf Case 2.}  If $x$  is not contained in any component of
codimension one of $\mbox{Sing }{\mathcal G}^{(d)}$ then
construct a ${\mathcal G}^{(d)}$-admissible local smooth projection to a
$(d-1)$-dimensional scheme $V^{(d-1)}$, and an elimination algebra
as in \ref{localproj} and \ref{eliminationalg},
$$\begin{array}{rrcl}
\beta_{d,d-1}: & (V^{(d)},x) & \longrightarrow &  (V^{(d-1)},x_1)\\
  & {\mathcal G}^{(d)} & & {\mathcal G}^{(d-1)}.
 \end{array}$$
By the induction hypothesis, $\gamma_{{\mathcal G}^{(d-1)}}(x_1)$
is defined. Now set
$$\gamma_{{\mathcal G}^{(d)}}(x)=(\mbox{ord}_{x}^{(d)}{\mathcal
G}^{(d)},\gamma_{{\mathcal G}^{(d-1)}}(x_1)).$$

\noindent $\bullet$ If $x\in \mbox{Sing }{\mathcal G}$ is not a
simple point, then let $\tilde{\mathcal G}^{(d)}$ be the twisted
algebra as in Definition \ref{definetwisting} with
$\omega=\mbox{ord}_x{\mathcal G}$. Then $x\in \mbox{Sing
}\tilde{\mathcal G}^{(d)}$ is a simple point and cases 1 and 2 can
be applied to $\tilde{\mathcal G}^{(d)}$. Now define
$$\gamma_{{\mathcal G}^{(d)}}(x)=(\mbox{ord}^{(d)}_x{\mathcal
G}^{(d)},\gamma_{\tilde{\mathcal G}^{(d-1)}}(x)),$$ where
$\gamma_{\tilde{\mathcal G}^{(d-1)}}(x)$ are the last
$(d-1)$-coordinates of the function $\gamma_{\tilde{\mathcal
G}^{(d)}}(x)$.

\

We will see next that $\gamma_{{\mathcal G}^{(d)}}$ is
upper semi-continuous and that it  stratifies $\mbox{Sing
}{\mathcal G}^{(d)}$ in  smooth strata.

\

The fact that this function takes only a finite number of values
follows  by induction. Thus it  only  remains  to show that   for
any value $(a_1,a_2,\ldots,a_d)\in ({\mathbb Q}^*)^d$, the set
$\{x\in V^{(d)}: \gamma_{\mathcal G}(x)\geq (a_1,a_2,\ldots,a_d)\}$ is
(locally) closed and  smooth (if it is non-empty). Observe that it is
enough to prove this fact in the case when $(a_1,a_2,\ldots,a_d)$
is the maximum value of the  function.  As in the
previous discussion we will use induction on the dimension of
$V^{(d)}$.

\

First suppose that $V^{(1)}$ is a one-dimensional-scheme   and let
$a_1\in {\mathbb Q}^*$ be any value such that $\{z\in V^{(1)}:
\gamma_{\mathcal G}(z)\geq a_1\}$ is non-empty. Since ${\mathcal
G}^{(1)}\neq 0$, $\{z\in V^{(1)}: \sigma_{\mathcal G}(z)\geq a_1\}$
consists of a finite number of closed points which is clearly a
smooth closed subscheme of $V^{(1)}$.

\

Assume now that part (i) of the theorem holds   for differential
algebras in any $n$-dimensional smooth scheme $V^{(n)}$ of finite
type over a field $k$ with $n<d$. We will show that it   also
holds for differential Rees algebras over a $d$-dimensional scheme
$V^{(d)}$ of finite type over a field $k$.

\

Let $(a_1,\ldots,a_d)\in {\mathbb Q}^*\times\ldots\times{\mathbb
Q}^*$ be the maximum value of $\gamma_{{\mathcal G}^{(d)}}^{(d)}$,
and let $x\in \mbox{Sing }{\mathcal G}^{(d)}$ with
$\gamma_{{\mathcal G}^{(d)}}^{(d)}(x)=(a_1,\ldots,a_d)$. We will
prove that there is an open subset $U^{(d)}\subset V^{(d)}$
containing $x$ such that $U^{(d)}\cap\{z\in V^{(d)}: \gamma_{{\mathcal
G}^{(d)}}^{(d)}(z)=(a_1,\ldots,a_d)\}$ is closed and smooth. We
distinguish three cases.

\

\noindent{\bf Case 1.} If $a_2=\ldots=a_d=\infty$ then  $x$ is
contained in a component of  codimension one  of $\mbox{Sing }{\mathcal
G}^{(d)}$. In this case by Lemma \ref{codimensionone} (applied to
${\mathcal G}^{(d)}$ if $x$ is a simple point, or to some
twisting, $\widetilde{\mathcal G}^{(d)}$ of ${\mathcal G}^{(d)}$
otherwise), there is an open neighborhood $U^{(d)}$ of $x$
satisfying the required property.

\

\noindent{\bf Case 2.}  If $a_2\neq \infty$ and $x$ is a simple
point  consider a ${\mathcal G}^{(d)}$-admissible local smooth projection
and an elimination algebra as in \ref{localproj} and
\ref{eliminationalg} in an open neighborhood $U^{(d)}$ of $x$:
$$\begin{array}{ccc}
(V^{(d)},x) & \longrightarrow & (V^{(d-1)},x_1)\\
 {\mathcal G}^{(d)} & &  {\mathcal G}^{(d-1)}.
 \end{array}$$
 Notice that then $\{z\in
U^{(d)}:\mbox{ord}_z^d{\mathcal G}^{(d)}= a_1\} $ can be
identified with $\mbox{Sing }{\mathcal G}^{(d-1)}$ in some open
neighborhood $U^{(d-1)}$ containing $x_1$. Therefore via this
identification
$$\{z\in U^{(d)}: \gamma_{{\mathcal G}^{(d)}}(z)=(a_1,\ldots,a_d)\}=$$
$$=\mbox{Sing }{\mathcal G}^{(d-1)}\cap \{z\in
U^{(d-1)}:\gamma_{{\mathcal G}^{(d-1)}}(z)=(a_2,\ldots,a_d)\}=$$
$$=\{z\in H^{(d-1)}:\gamma_{{\mathcal G}^{(d-1)}}(z)=(a_2,\ldots,a_d)\}$$
for some open subset $ H^{(d-1)}\subset V^{(d-1)}$.  Restricting
 $U^{(d-1)}$  if  necessary we may assume that $(a_2,\ldots,a_d)$
is actually the maximum of $\gamma_{{\mathcal G}^{(d-1)}}$. According to our
inductive hypothesis there is an open neighborhood of $x_1$ where
$\{z\in U^{d-1}:\gamma_{{\mathcal G}^{(d-1)}}(z)=(a_2,\ldots,a_d)\}$
is locally closed and  smooth. Again, via the identification
$\mbox{Sing }{\mathcal G}^{(d)}\cap U^{(d)}$ with $\mbox{Sing
}{\mathcal G}^{(d-1)}\cap U^{(d-1)}$ we conclude that  there is an
open neighborhood of $x$ where the stratum $\{z\in V^{(d)}:
\gamma_{\mathcal G^{(d)}}(z)=(a_1,\ldots,a_d)\}$ is closed and smooth.

\

\noindent{\bf Case 3.}  If $a_2\neq \infty$ and $x$ is not a
simple point, then, in a suitable neighborhood of $x$,  replace
${\mathcal G}^{(d)}$ by $\widetilde{\mathcal G}^{(d)}$ as in
Theorem \ref{twisting}. By restricting to a smaller neighborhood
if needed, it can be assumed that in    addition,
$\{z:\mbox{ord}_z{\mathcal G}^{(d)}=a_1\}=\mbox{Sing
}\widetilde{\mathcal G}^{(d)}$. Now the
 argument in Case 2 can  be applied to $\widetilde{\mathcal G}^{(d)}$. } \qed
\end{Paragraph}

\part{Epilogue and example}\label{partEpilogue}
\label{epilogue} Let ${\mathcal G}^{(d)}$ be a Rees algebra on a
smooth scheme $d$-dimensional smooth scheme $V^{(d)}$ over a perfect  field
$k$. The study of a stratification on $\Sing {\mathcal G}^{(d)} $
achieved by means of an upper semi-continuous function is one
example of an application of Main Theorem \ref{ordertau}. However, this stratification is mainly interesting   due to the following
fact:  as in Part \ref{Char0} (via the dictionary between
Rees algebras and pairs) similar satellite functions can be
defined thanks to Theorems \ref{1theo}, \ref{ordertau} and
\ref{twisting}. In this way an upper semi-continuous function is
constructed whose maximum value determines permissible centers,
and the blow-up along these centers produces a simplification of
the singularities. To be precise, once we blow-up at the smooth
center defined by the function (on the worst points), a new
upper semi-continuous function is defined, which provides a new
stratification and a new closed and smooth stratum of worst
singularities. We then consider the blow-up at such center and so on.

\

A number of exceptional hypersurfaces arise in this process of
monoidal transforms, and it is important that these hypersurfaces
have normal crossings. So we have  to define a procedure so that the
maximum stratum (center of the monoidal transform) have normal
crossings with the exceptional hypersurfaces introduced in the
previous steps. Here is where the second satellite functions play an important
role (see \ref{lasat2}).

\

On the other hand, the notion of codimensional type in Definition
\ref{cdtype} provides a natural form of induction used in resolution
problems:

\

-  Observe first that when ${\mathcal G}^{(d)}$ is of codimensional
type $d$, then $\Sing({\mathcal G}^{(d)})$ is a zero dimensional
closed set, and a resolution is achieved by blowing up these closed
points.

\

- When ${\mathcal G}^{(d)}$ is of codimensional type $\geq m$, then
by  Theorem \ref{twisting}
a new Rees algebra ${\widetilde{\mathcal G}}^{(d)}$ of
codimensional type $\geq m+1$ can be attached to ${\mathcal
G}^{(d)}$. Theorem \ref{twistingtransform} says that this Rees
algebra is determined up to integral closure of algebras. This is
what allows us to define the upper semi-continuous functions after
successive monoidal transforms, and it also leads to the reduction to
the monomial case. In fact, if we assume by induction an algorithm
of resolution for algebras of codimensional type $\geq m+1$, then
 the following theorem holds:

\

\noindent{\bf Theorem (Reduction to the monomial case). \cite[Corollary 6.15]{positive}} 
{\em  Let $V^{(d)}$ be a $d$-dimensional smooth scheme over a perfect field $k$,   
let ${\mathcal G}^{(d)}$ be a Rees algebra of codimensional type $m\geq 0$, and let $E$ be a set of 
smooth hypersurfaces in $V^{(d)}$ with normal crossings. 
Assume that there is an algorithm of resolution of Rees algebras of codimensional type $\geq m+1$. Then it is possible  to define a sequence of permissible transformations, 
\begin{equation}\label{Atryo2}
\begin{array}{ccccccc}
(V^{(d)},\mathcal{G}^{(d)},E^{(d)}) & \leftarrow &
(V^{(d)}_1,\mathcal{G}_1^{(d)},E_1^{(d)}) & \leftarrow & \cdots &
\leftarrow & (V^{(d)}_s,\mathcal{G}_s^{(d)},E_s^{(d)}) \\
\downarrow^{\beta} & & \downarrow^{\beta_1} & &  \cdots & & \downarrow^{\beta_s}   \\
(V^{(d-m)},\mathcal{G}^{(d-m)},E^{(d-m)}) & \leftarrow &
(V^{(d-m)}_1,\mathcal{G}_1^{(d-m)},E_1^{(d-m)}) & \leftarrow &
\cdots & \leftarrow &
(V^{(d-m)}_s,\mathcal{G}_s^{(d-m)},E_s^{(d-m)})
\end{array}
\end{equation}
so that  
$$ \max   \ord^{(d-m)}_{{\mathcal G}^{(d)}}  \geq \max  \vword^{(d-m)}_{{\mathcal G}_1^{(d)}} \geq \dots \geq \max
 \vword^{(d-m)}_{{\mathcal G}^{(d)}_s},$$
where in addition   $ \max \vword^{(d-m)}_{{\mathcal G}_s}=0$. 
 In other words, the sequence of transformations can be
 defined so that  $(V^{(d)}_s,\mathcal{G}^{(d)}_s,E^{(d)}_s)$ is in the
   monomial case  as described in \ref{rgn}.}

\

So Theorem \ref{twistingtransform} provides our form of induction.
However, to be precise, the invariant dealt with in that Theorem
is essentially that in  (\ref{eqstdeg1}) of \ref{ragadprs3} (or
say, the first satellite function in equation (\ref{eqdword}) of
\ref{lasat1}). As indicated above, after the first monoidal
transformation, the upper semi-continuous function that defines the
reduction to the monomial case makes use of the second satellite
function as in (\ref{ect}). The reader can find the formulation of
Theorem \ref{twistingtransform} in terms of the second satellite
function and other technical aspects  in \cite{positive}.

\


 For  the case of characteristic zero it is simple to extend
 sequence (\ref{Atryo2}) (i.e., the monomial case) to a resolution of $(V^{(d)},\mathcal{G}^{(d)},E^{(d)})$
 (see  Steps A and B in Section \ref{algoritmicpairs}). When the
 characteristic is positive, the containment  $\beta_s(\mbox{Sing }\mathcal{G}_s^{(d)})\subset
 \mbox{Sing }\mathcal{G}_s^{(d-m)}$ may be strict, and then a
  resolution of $(V^{(d-m)}_s,\mathcal{G}_s^{(d-m)},E_s^{(d-m)})$ may not lift to a  resolution of
 $(V^{(d)}_s,\mathcal{G}_s^{(d)},E_s^{(d)})$ (see \cite{BVV} for concrete  examples). However the condition of $(V^{(d-m)}_s,\mathcal{G}_s^{(d-m)},E_s^{(d-m)})$ being monomial in positive characteristic
 (i.e. the elimination algebra being monomial) opens the way to new invariants, as those treated in  \cite{BVV}.
 We hope to be
 able to address the monomial case in arbitrary characteristic in the
 future.

\

We conclude this section with an example to illustrate the computation of an elimination algebra and its use in stratification. We also indicate how  our resolution functions  define a sequence of blow-ups that lead to  the monomial case. This example has been treated in \cite{HauserBul} and \cite{Hauser} to show some of the pathologies raising in positive characteristic. This is one of the cases where the natural resolution invariant (defined as a  generalization of the one  used in characteristic zero) grows after a finite number of blow-ups. 

\begin{Example}\label{canguro}
{\rm Assume that $k$ is a field of characteristic 2,  let
$V^{(3)}$ be the affine 3-dimensional space $\mbox{Spec}
(k[X,Y,Z])$, and let
$$S:=\{f=Z^2+(Y^7+YX^4)\in k[X,Y,Z]=0\}.$$
Clearly, the maximum order at points of $S$ is two. This maximum
is reached at the points of the curve $\{Z=0, Y^3+X^2=0\}$, but
this is not a smooth closed subscheme, and hence  we are
forced to look for other invariants that refine the order
function.

\

Let ${\mathcal G}$ be the differential Rees algebra generated by $f$ in
degree two: $${\mathcal G}={\mathcal
O}_{V^{(3)}}[Z^2+(Y^7+YX^4)W^2,  (Y^3+X^2)^2W].$$ Notice
that
  $\mbox{Max-ord}^{(3)}_{\mathcal G}=1$,  and that $$\mbox{Sing }{\mathcal
G}= \mbox{\underline{Max}-ord}^{(3)}_{{\mathcal
G}}=\{Z=0,Y^3+X^2=0\}.$$ Again, observe that the function
$\mbox{ord}^{(3)}_{{\mathcal G}}$ is too coarse: its singular
locus is not even smooth.

\

Let $V^{(2)}=\mbox{Spec}(k[X,Y])$.  We choose the ${\mathcal
G}$-admissible projection $$ \beta_{3,2}: V^{(3)}\to V^{(2)} $$ and
compute the corresponding elimination algebra, ${\mathcal
R}_{\mathcal G}$:
$${\mathcal R}_{\mathcal G}={\mathcal O}_{V^{(2)}}[(Y^3+X^2)^2W].$$
Now $\mbox{Max-ord}^{(2)}_{{\mathcal G}}=
\mbox{Max-ord}_{{\mathcal R}_{\mathcal G}}= 4$,
 and
$\mbox{\underline{Max}-ord}^{(2)}_{{\mathcal G}}=\{(0,0,0)\}$. The
procedure now involves associating a simple   differential
algebra to  $\mbox{\underline{Max}-ord}^{(2)}_{{\mathcal G}}$ and
then projecting onto some smooth scheme of dimension 1.

\

The first monoidal transformation is the blow-up at the origin,
$V^{(3)}\leftarrow  V^{(3)}_1$, which induces a blow-up at
$V^{(2)}$ with center $\{(0,0)\}$, $V^{(2)}\leftarrow  V^{(2)}_1$.
Recall that by Theorem \ref{1theo}, there are commutative diagrams
of monoidal transformations, restrictions and elimination
algebras.

\

Consider  the affine charts $U_{1,Y}^{(3)}=\mbox{Spec
}\left(k\left[\frac{X}{Y},Y,\frac{Z}{Y}\right]\right)\subset
V^{(3)}_1$, and $U_{1,Y}^{(2)}=\mbox{Spec
}\left(k\left[\frac{X}{Y},Y\right]\right)\subset V^{(2)}_1$. To
simplify notation, set again $X=\frac{X}{Y}$ and $Z=\frac{Z}{Y}$:
$$\begin{array}{ccc}
V^{(3)} & \stackrel{\pi_1}{\longleftarrow} & V^{(3)}_1  \\
 & & \cup  \\
  & & U_{1,Y} \\
 & & f_1=Z^2+Y^3\cdot(Y+X^2)^2.
\end{array}$$
In $U_{1,Y}^{(3)}$ consider the weighted transform of ${\mathcal G}$,
$${\mathcal G}_1= {\mathcal O}_{U_{1,Y}^{(3)}}[(Z^2+Y^3\cdot(Y+X^2)^2)W^2,  Y^3(Y+X^2)^2W],$$  and the weak
transform of ${\mathcal R}_{\mathcal G}$ in $U_{1,Y}^{(2)}$,
${\mathcal R}_{{\mathcal G}_1}= {\mathcal
O}_{U_{1,Y}^{(2)}}[Y^3(Y+X^2)^2W].$

\

Notice that $\mbox{Max-ord}^{(3)}_{\mathcal G_1}=1$,  and that
$\mbox{Max-ord}^{(2)}_{{\mathcal G}_1}=\mbox{Max-w-ord}_{{\mathcal
R}_{{\mathcal G}_1}}=2$, so this invariant has dropped, and hence
the second satellite function plays a role counting exceptional
divisors (see \ref{lasat2}). Now   the same procedure that works
for algorithmic resolution in characteristic zero  applies here
(see Section \ref{Char0}), and after two more blow-ups
  at closed points (the centers that are determined using the upper semi-continuous functions derived from Theorem  \ref{ordertau}), the monomial case is
achieved.

}

\end{Example}

\end{document}